\documentclass[a4paper]{article}
\usepackage{authblk}
\usepackage{hyperref}
\usepackage{url}
\usepackage{amssymb,amsthm,latexsym,amsmath}
\usepackage{msc}
\usepackage{graphicx}
\usepackage{subcaption}
\usepackage{textcomp} 
\begin{document}
\title{Maximal Arrangement of Dominos in the Diamond}
\author[1]{Dominique D\'{e}s\'{e}rable}
\author[2]{Rolf Hoffmann}
\author[3]{Franciszek Seredyński}

\affil[1]{Institut National des Sciences Appliqu\'{e}es, Rennes, France\newline
\url{domidese@gmail.com}
}
\affil[2]{Technische Universit{\"a}t Darmstadt, Darmstadt, Germany\newline
\url{hoffmann@ra.informatik.tu-darmstadt.de}
}
\affil[3]{Department of Mathematics and Natural Sciences, Cardinal Stefan Wyszyński University, Warsaw, Poland\newline
\url{f.seredynski@uksw.edu.pl}
}
\maketitle
%
%
\begin{abstract}
\noindent ``Dominos'' are special entities consisting of a hard dimer--like kernel surrounded by a soft hull and governed by local interactions.
``Soft hull'' and ``hard kernel'' mean that the hulls can overlap while the kernel acts under a repulsive potential.
Unlike the dimer problem in statistical physics, which lists the number of all possible configurations for a given $ n \times n $ lattice, the more modest goal herein is to provide lower and upper bounds for the maximum allowed number of dominos in the diamond.
In this NP problem, a deterministic construction rule is proposed and leads to a sub--optimal solution $ \psi_n $ as a  {\em lower} bound.
A certain disorder is then injected and leads to an {\em upper} bound  $  \overline{\psi}_n $ reachable or not.
In some cases, the lower and upper bounds coincide, so 
$  
          \overline{\psi}_n  =  \psi_n
$ 
becomes the {\em  exact} number of dominos for a maximum configuration. 
\end{abstract}
%
%
%
\newpage
\tableofcontents
\newpage
%
%
\section{Introduction}
\label{Section: Introduction}
%
%
The problem of maximal arrangement of dominos fits into the broad subject of discrete optimization, tiling and covering in two--dimensional spaces and polyominos.
A  ``domino'' is a special entity consisting of a hard dimer--like kernel surrounded by a soft hull and governed by local interactions.
``Soft hull'' and ``hard kernel'' mean that the hulls can overlap while the kernel acts under a repulsive potential.
In other words, a domino is a $ 4 \times 3 $ rectangle with one $ 2 \times 1 $ kernel in the center and with two possible layouts: horizontal or vertical.
Hull--hull overlap must be favored in a {\em max} problem and must be avoided in a {\em min} problem while hull--kernel contact is precluded.
In this paper, the problem is to arrange a maximum of dominos in the rhombic {\em diamond}, a 
$
           \frac{\pi} {4}\mbox{--tilted} 
$
square.  

This work is part of a research activity started more than five years ago in the area of discrete optimization and centered around a common topic: 
a model of cellular automata (CA) as a tool for simulating populations of small polyomino--like entities with local interaction.
Their interaction is modeled as an idealized multicellular ``tile'' with a hard kernel surrounded by a soft hull. 
The objective function is either to maximize or to minimize the population.

In
\cite{Hoffmann:Deserable-2017,Hoffmann:Deserable-2019}
the tile was a spin--like left--up--right--down domino defined from four $ 3 \times 3 $ Moore templates in a multi--agent system where the agent was controlled by a finite state machine evolved by a genetic algorithm forming domino patterns and the objective was to arrange a maximum of dominos in the square.
In
\cite{Hoffmann:Deserable-Seredynski-2019}
the domino was redefined as a rectangle with a spin--like kernel surrounded by a 10--cell hull, close to its current definition, and the genetic evolution was replaced by a probabilistic CA.
 In
\cite{Hoffmann:Deserable-Seredynski-2021a}
the robustness of the CA rule with respect to the geometry of the host shape was tested and extended to the rhombic {\em diamond}.
 In
\cite{Hoffmann:Deserable-Seredynski-2021b}
the robustness of the CA rule with respect to the objective function was tested and the minimization of the population of dominos in the square was carried out without too much difficulty.
 In
\cite{Hoffmann:Seredynski-2021,Hoffmann:Deserable-Seredynski-2022}
the robustness of the CA rule with respect to the entity --\,including its interacting field\,-- was tested:
the practical object was a sensor with its sensing area, the kernel a monomino as ``sensor point'', the hull a von Neumann neighborhood of range 2 as ``sensing area'' and the objective function was the minimization of the population of sensors.

This paper is a continuation of
\cite{Hoffmann:Deserable-Seredynski-2021a}
in the sense that the theoretical framework presented therein to support the results of the simulation had just been initiated.
Moreover, the $ n $--sample was only suitable for small sizes; for larger sizes, the model would have exhibited a divergence.
Here, the deterministic construction rule which is fixed has the advantage of offering the greatest possible symmetry and simplicity and, on the other hand, it leads to a quasi--optimal configuration before giving an upper bound for the requested maximum.

The following section defines the geometric and physical frameworks of this study. At first, the ``domino'' entity is disambiguated in order to remove the troublesome homonymy with the popular domino--dimer in statistical physics.
Two density indices are suggested, that may be useful for further examination.
Then we present the deterministic rule --\,one could even say the {\em axiom}\,-- of construction which leads to the quasi--optimal configuration.
Finally, disorder injection in the configuration leads to the required upper bound.
%
%
\section{General Statements}
\label{Section: General Statements}
%
%
%
%
\subsection{The Domino Entity}                            
%
%
%
%
\begin{figure}
\centering
\includegraphics[width=7cm]{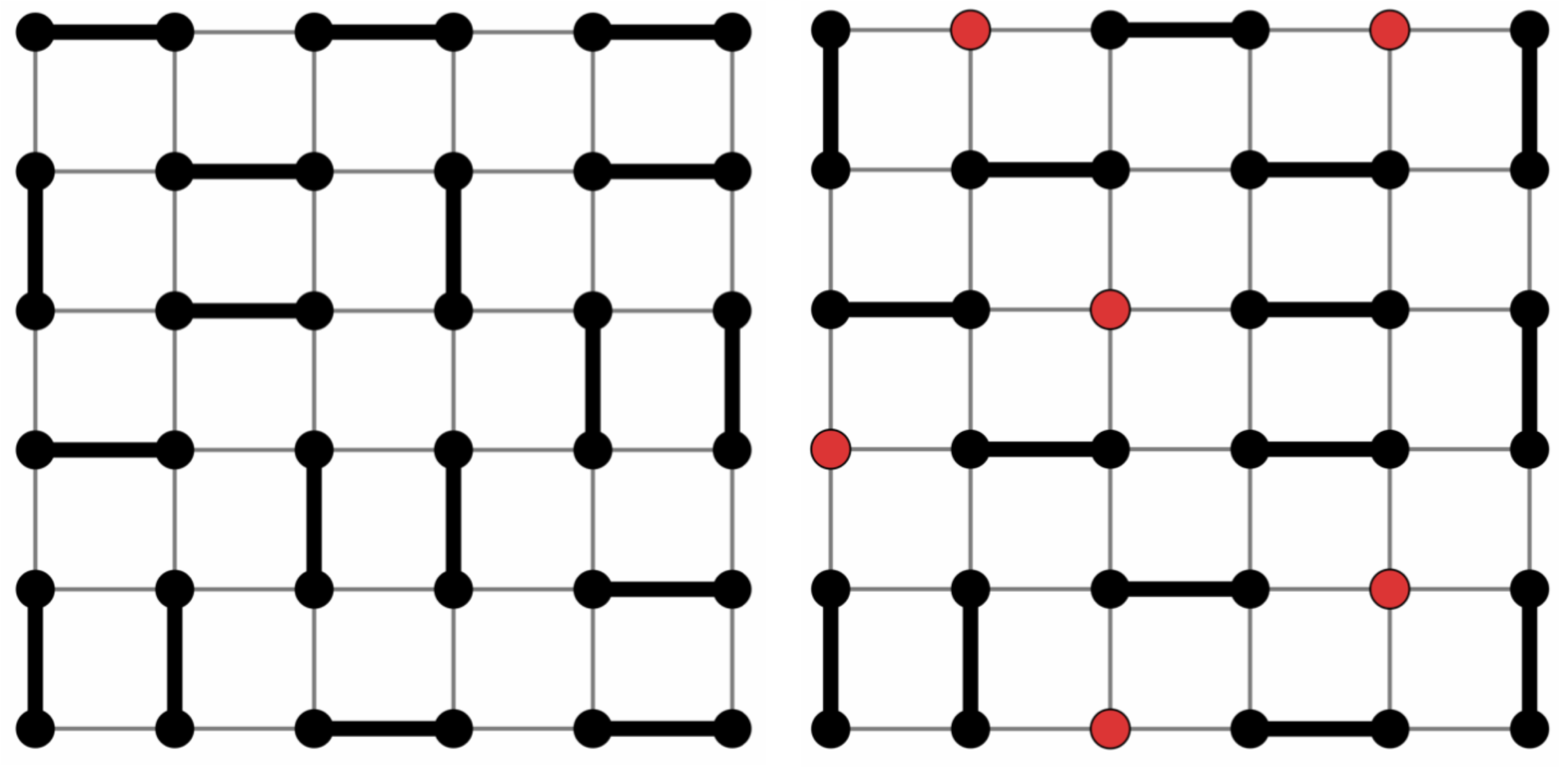}  
\caption{Dimer configurations in a $ 6 \times 6 $ square lattice
\cite{Alegra:Fortin-2014}.
({\em Left}) Perfect matching without monomer.
({\em Right}) Dimers with monomers (red dots).
} 
\label{Figure: DimerMonomer }
\end{figure}
%
%
%
The term ``domino'' can have several meanings.
Basically, a domino is a twofold object belonging to the polyomino family.
It is encountered in statistical physics in the form of a ``dimer''
\cite{Temperley:Fisher-1961,Kasteleyn-1961}
as a pure tiling problem or closely related to short--range interaction coupling in spin systems
\cite{Niss-2005}
(Fig.\,\ref{Figure: DimerMonomer }).
It is also correlated with the alternating--sign matrices
\cite{Elkies:Kuperberg-Larsen:Propp-1992}
and the {\em square--ice} model of statistical mechanics
\cite{DiMarzio:Stillinger-1963,Lieb-1967}:
namely, the six possible patterns for a vertex in a 4--regular digraph with two incoming arcs and two outgoing arcs.
But the most amazing result is the existence of a circular area at the thermodynamic limit --the ``Arctic Circle''-- induced by random domino tilings in the Aztec diamond
\cite{Jockusch:Propp:Shor-1998}.

\begin{figure}
\centering
\includegraphics[width=9cm]{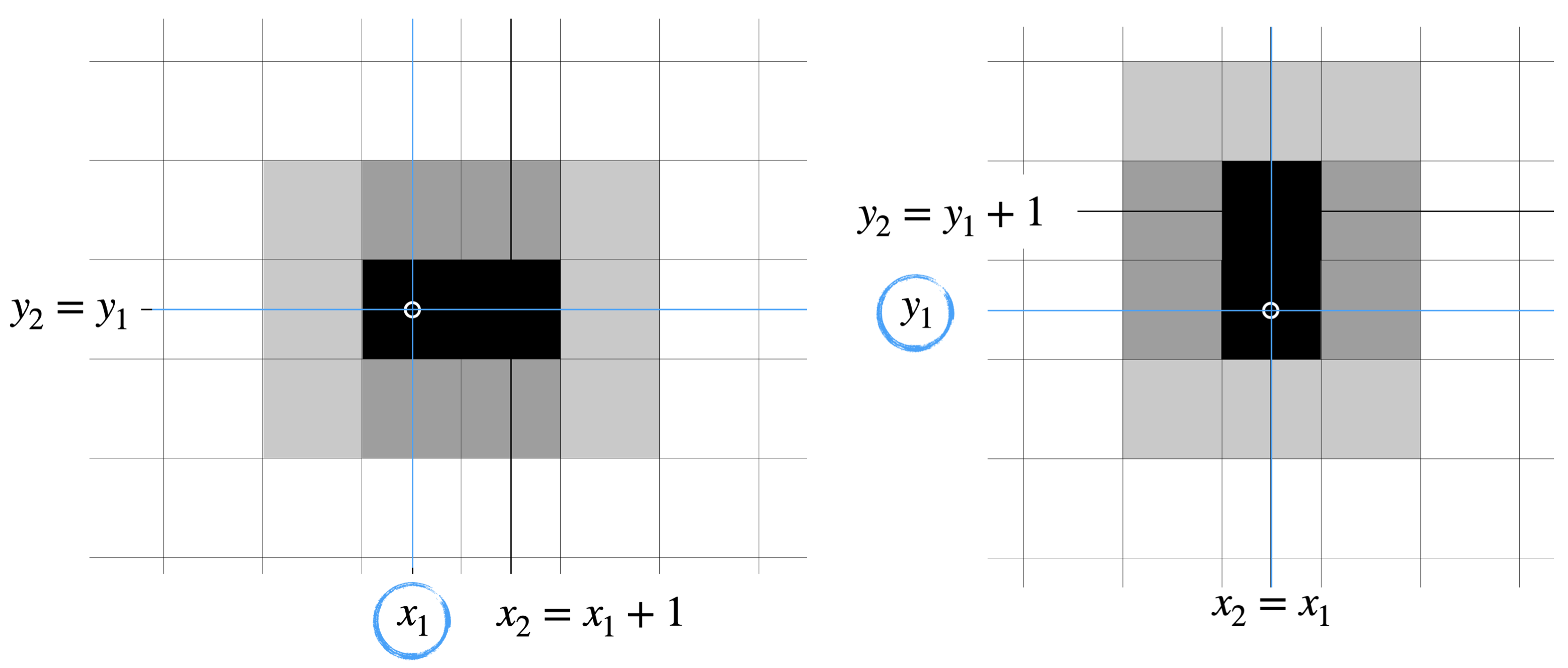}  
\caption{Horizontal and vertical dominos as the union of two Moore neighborhoods 
$
          \mathcal{ M }( \omega_1 \,)   \cup   \mathcal{ M }( \omega_2 \,) 
$ 
anchored at point $  \omega_1  $ in the
$  
            \mathbb{Z}^2
$ 
lattice.
} 
\label{Figure: Domino-DhDv }
\end{figure}

Our domino entity herein is quite different. 
Let 
$
             \mathcal{ B }  =  \{ e_x, e_y \}    =   \{ ( 1 , 0 ), ( 0 , 1 ) \} 
$
be the orthonormal basis in the 
$  
            \mathbb{Z}  \times   \mathbb{Z}
$ 
lattice and let 
\[
         \mathcal{ M }\,(( x_0 , y_0 \,))  =  
         \mathcal{ M }_{ r }\,(( x_0 , y_0 \,))  =  \{ ( x , y ) : | x - x_0 | \le r ,  | y - y_0 | \le r \}  \ \ ( r = 1 )
\]
be the Moore neighborhood of range $ r = 1 $ surrounding a given point $ ( x_0, y_0 ) $.
Without loss of generality, let us choose the two neighboring points 
$ 
           \omega_1 =  ( x_1 , y_1 \,) 
$ 
and 
$ 
	  \omega_2 = ( x_2 , y_2 ) =  \omega_1 +   e \ (  e \in  \mathcal{ B } )
$ 
and the vector
$
          \omega_1  \omega_2 
$
anchored at $ \omega_1 $. A  {\em domino} 
$
        (  \omega_1,  \omega_2 )
$
is the set of the 12 points
\[            
              \mathcal{ D }_{ \omega_1 \, \omega_2 }  =  \mathcal{ M }( \omega_1 \,)   \cup   \mathcal{ M }( \omega_2 \,)  =
                                                                                    \mathcal{ M }( \omega_1 \,)   \cup   \mathcal{ M }( \omega_1 + e \,)  
\]
and where
$           
              \mathcal{ K }_{ \omega_1 \, \omega_2 }  =    \{   \omega_1 , \omega_2   \} 
$
is the {\em kernel} while the set of the 10 points
\[            
              \mathcal{ H }_{ \omega_1 \, \omega_2 }  =    \mathcal{ D }_{ \omega_1 \, \omega_2 }  -   \mathcal{ K }_{ \omega_1 \, \omega_2 } 
\]
surrounding the kernel is the {\em hull}.
Points in
$  
            \mathbb{Z}^2
$ 
and cells in the cellular space are associated by duality as shown in 
Fig.\,\ref{Figure: Domino-DhDv }.
%
%
%
%
%
{\em Horizontal} and {\em vertical} dominos are such that $ e = e_x $ and $ e = e_y $ respectively.

Perhaps the best way to remove the troublesome homonymy from its namesake in statistical physics would be to rename our composite domino herein.
A name like ``$\{4,3\}$--do--deca--mino'' ({\em do--} for kernel, {\em deca--} for hull) might then be more relevant.

A layout of dominos, whether maximum or minimum, is governed by the hull--kernel repulsive {\em exclusion} rule  
\[            
               \forall   \omega_1   \in   \mathbb{Z}^2   \ ; \
               \forall   \omega'_1   \in   \mathbb{Z}^2  \ : \
               \mathcal{ H }_{ \omega'_1 \, \omega'_2 }  \cap  \mathcal{ K }_{ \omega_1 \, \omega_2 }  =  \varnothing 
\]
for any pair 
$
             (  \mathcal{ D }_{ \omega_1 \, \omega_2 }  ,   \mathcal{ D }_{ \omega'_1 \, \omega'_2 }  ).
$
This means, on the contrary, that the hulls of local neighbors are allowed to {\em overlap}.
Each cell  
$ 
             \omega    \in    \mathbb{Z}^2  
$  
is assigned a {\em cover} (or overlap) level 
$
           v_{ \omega }  \  ( 0 \le   v_{ \omega }  \le 4 )
$  
which measures the number of dominos covering it.
A {\em lacunar} void is an empty cell  $ \omega $ such that $ v_{ \omega }  = 0 $.
The hull--kernel exclusion implies
$
               v_{ \omega_1 }  =     v_{ \omega_2 }  =  1
$
for any 
$
              \mathcal{ K }_{ \omega_1 \, \omega_2 }.
$
Figure \ref{Figure: Four-layouts }
shows various scenarios of overlapping. 

The kernel of a domino is allowed to evolve strictly within a given convex shape.
It follows that the hull will necessarily exceed the envelope of the convex.
In the simplest case, if we consider a square array of size $ n \times n $ in which a population of dominos will evolve, then the relevant convex will be the square $ \mathcal{S}_{n} $ of size  
$
       ( n + 2 ) \times ( n+ 2 ).
$

One of the best down--to--earth metaphor illustrating our domino entity is the arrangement of cars in a car–park. The {\em kernel} is the car itself while the {\em hull} is the free space left to open the doors, the tailgate or possibly the hood. 
Incidentally, the {\em diamond} might represent a complex idealized host shape --\,the parking lot\,-- that could exist in old city centers or in steep regions.
Similar physical distancing situations are found in structures such as a set of tables in a classroom or in any public space.
Our problem of maximal arrangement of dominos fits into the broad subject of tiling and covering in two-dimensional spaces
\cite{Zong-2014}
and can be brought back to that of ellipsoid packing
\cite{Donev:Stillinger:Chaikin:Torquato-2004,Borzsonyi:Stannarius-2013}.
%
%
%
\subsection{Density Measurement}                            
\label{Subsection: Density Measurement-- }
%
%
%
For a given domino 
$
        (  \omega_1,  \omega_2 )
$
we define the two following measures:
%
%
\begin{itemize}
%
%
\item          
%
%
the {\em overlap} index 
\vspace{ - 2mm }
\[
          \nu_{  \omega_1,  \omega_2 } =  \sum_{ \mathcal{ D }_{ \omega_1 \, \omega_2 } }^{ }   v_{ \omega } = 
           \sum_{ \mathcal{ K }_{ \omega_1 \, \omega_2 } }^{ }   v_{ \omega }   +  
           \sum_{ \mathcal{ H }_{ \omega_1 \, \omega_2 } }^{ }   v_{ \omega }   =
          2  +   \sum_{ \mathcal{ H }_{ \omega_1 \, \omega_2 } }^{ }   v_{ \omega } 
\]
as the sum of all cover levels $ v_{ \omega } $ on  
$ 
             \omega    \in    \mathcal{ D }_{ \omega_1 \, \omega_2 }
$  
and
%
\item         
%
%
the {\em occupancy} index 
\vspace{ - 2mm }
\[
          \rho_{  \omega_1,  \omega_2 } =  \sum_{ \mathcal{ D }_{ \omega_1 \, \omega_2 } }^{ }   \tau_{ \omega } = 
           \sum_{ \mathcal{ K }_{ \omega_1 \, \omega_2 } }^{ }   \tau_{ \omega }   +  
           \sum_{ \mathcal{ H }_{ \omega_1 \, \omega_2 } }^{ }   \tau_{ \omega }   =
          2  +   \sum_{ \mathcal{ H }_{ \omega_1 \, \omega_2 } }^{ }   \tau_{ \omega } 
\]
where $  \tau_{ \omega }  $  is the occupancy ratio (the inverse $ 1 / v_{ \omega } $ of the cover level) of cell
$ 
             \omega    \in    \mathcal{ D }_{ \omega_1 \, \omega_2 }.
$  
\end{itemize}
%
%
%
\begin{figure}
\centering
\includegraphics[width=9cm]{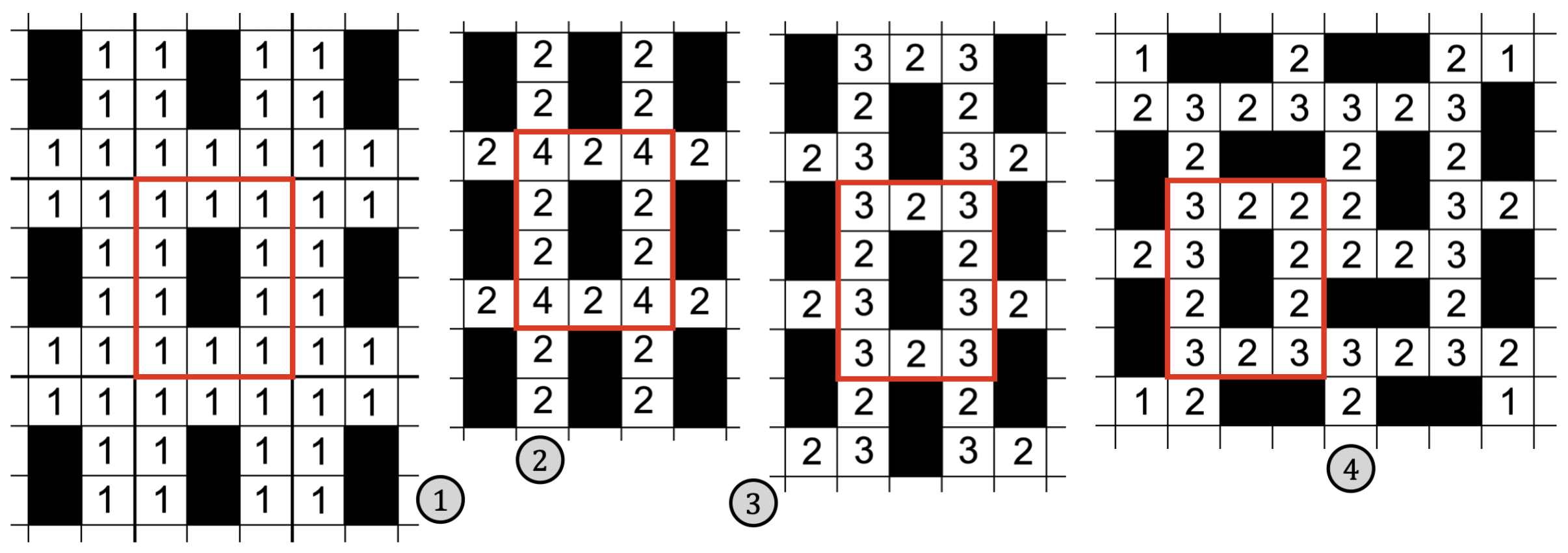}  
\caption{Four typical arrangements of dominos:
(1) loosely coupled physical distancing configuration
(2) tightly coupled orthotropic configuration
(3) tightly coupled staggered configuration
(4) tightly coupled isotropic configuration. 
} 
\label{Figure: Four-layouts }
\end{figure}
%
%
For example, the two measures above applied to the red domino in the four cases of  
Fig.\,\ref{Figure: Four-layouts }
give:
%
%
\begin{enumerate}
%
%
\item   $  \nu_{  \omega_1,  \omega_2 }  =   2 + 10 \times 1  = 12  $ \ ;    \hspace{13mm}  
           $  \rho_{  \omega_1,  \omega_2 }  =   2 + 10 \times 1  = 12  $
\item   $   \nu_{  \omega_1,  \omega_2 }  =   2 + 6 \times 2 + 4 \times 4  =  30  $ \ ;   \hspace{2.5mm} 
           $  \rho_{  \omega_1,  \omega_2 }  =   2 + 6 \times 1/2 + 4 \times 1/4 = 6  $	
\item   $  \nu_{  \omega_1,  \omega_2 } =   2 + 4 \times 2 + 6 \times 3  =  28  $ \ ;    \hspace{2.5mm} 
           $ \rho_{  \omega_1,  \omega_2 }  =   2 + 4 \times 1/2 + 6 \times 1/3 = 6  $	
\item   $   \nu_{  \omega_1,  \omega_2 } =   2 + 6 \times 2 + 4 \times 3  =  26  $ \ ;    \hspace{5mm} 
           $  \rho_{  \omega_1,  \omega_2 }  =   2 + 6 \times 1/2 + 4 \times 1/3 = 19 / 3  $	
\end{enumerate}
%
%
In the first case, the overlap index is minimal while the occupancy index is maximal for a {\em minimal} layout.
Depending of their objective function, the first and the second configurations are optimal, but an overall configuration is constrained by the boundary conditions.
Second and third cases show better accuracy for the overlap index.
Finally, the occupancy index provides a relationship between the number of dominos and the surface of the shape. Thus, if we reconsider the square
$ 
          \mathcal{S}_{n} 
$
it follows that
%
%
\begin{equation}
\label{equation: The Square : Occupancy vs. Surface }
          \sum_{ k = 1 }^{ \xi_n }   \rho \, ( k )   +   \varpi_n  \, = \
          \mid \mathcal{S}_{n}  \mid  \ = ( n + 2 )^2
\end{equation}
%
%
where $ \xi_n $ denotes the number of dominos for a given population in the square while
$
             \rho \, ( k )
$
and
$
             \varpi_n 
$
represent respectively the occupancy index of domino $ k $ and the number of lacunar voids.
Here $ k $ defines an arbitrary numbering of the population of dominos.
This relation holds whatever the population of dominos (either minimal or maximal, or neither minimal nor maximal).

%
%
%
\subsection{Geometry of the Diamond}                            
%
%
%
%
%
\begin{figure}
\centering
\includegraphics[width=8cm]{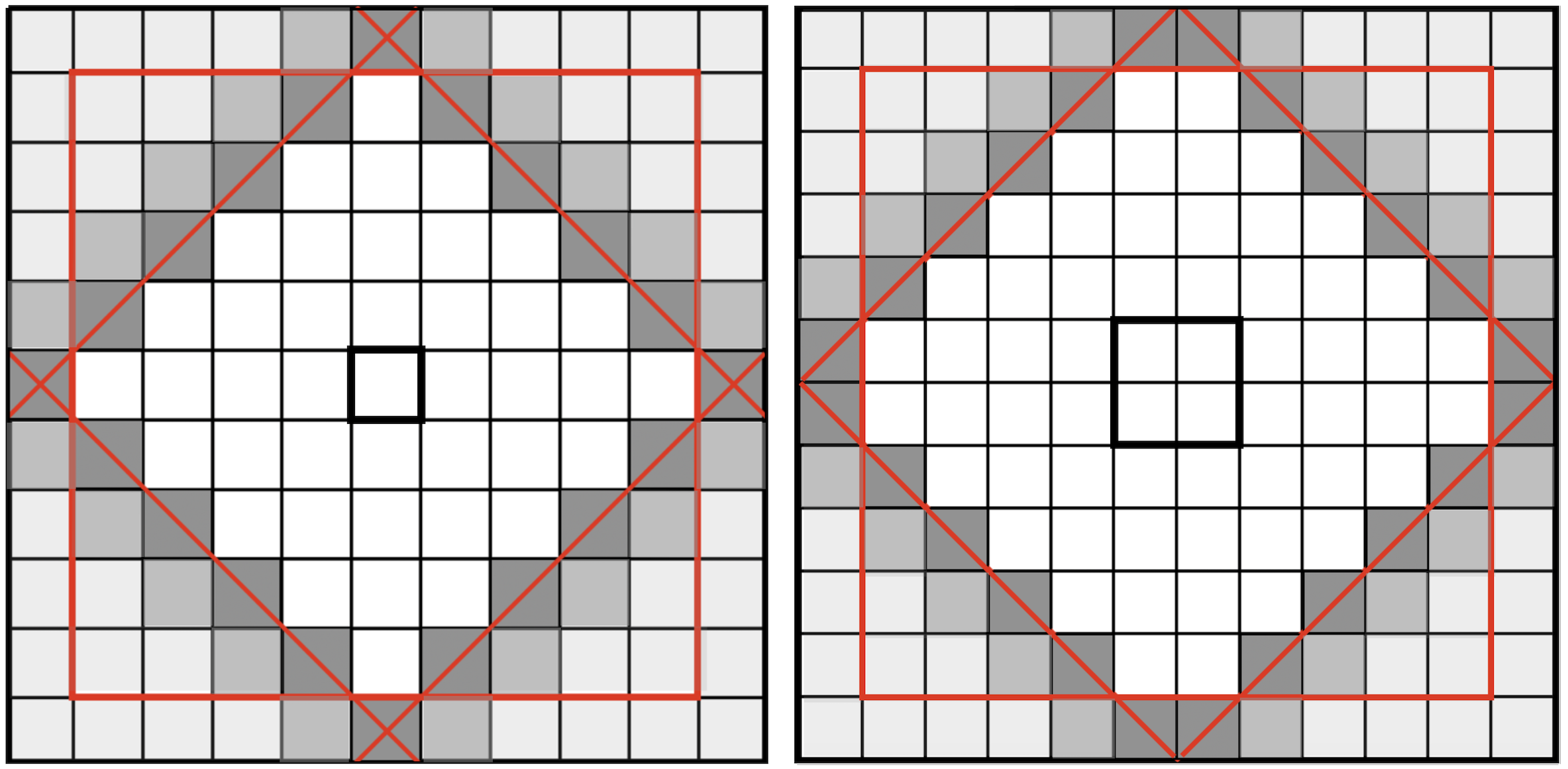}  
\caption{The (cellular) diamond
$
       \mathcal{D}_n
$
as a $\pi/4$--rotated (cellular) square.\newline
%
%
-- $n$ odd ({\em left}):   
$
       \mathcal{D}_{ 9 }  
$
as a von Neumann neighborhood of range $ r = 5 $ with $ 2 r = n + 1 $ and
$
        \mid \mathcal{D}_{ n }\mid \   =  n_{ r } =   61.
$
For the white diamond
$
        \mid \mathcal{D}_{ n - 2  }\mid \   =  n_{ r  - 1 } =   41
$
and for the full domain
$
       \mid \mathcal{ S }_{ n }  \cap   \mathcal{D}_{ n + 2 }\mid  \, =  
       n_{ r  +  1 }  - 4   =   81.
$
\newline 
%
%
-- $n$ even ({\em right}):
$
       \mathcal{D}_{ 10 }  
$
as an Aztec diamond of range $ r = 6 $ with $ 2 r = n + 2 $ and 
$
        \mid \mathcal{D}_{ n }\mid \   =  n_{ r }^{ ( a ) } =   84.
$
For the white diamond
$
        \mid \mathcal{D}_{ n - 2  }\mid \   =  n_{ r - 1 }^{ ( a ) }  =   60
$
and for the full domain
$
       \mid \mathcal{ S }_{ n }  \cap   \mathcal{D}_{ n + 2 }\mid  \, =  
       n_{ r  +  1 }^{ ( a ) }  -  8   =   104.
$
} 
\label{Figure: D9-D10 }
\end{figure}
%
%
\noindent The (cellular) {\em diamond}
$
       \mathcal{D}_n
$
is defined herein as a $\pi$/4--rotated square inscribed in a square array
$
       \mathcal{S}_n
$
of
$
       ( n + 2 ) \times ( n+ 2 )
$
cells including a border with perimeter
$
     4 n + 4
$
enclosing a 
$
     n \times n
$
square (see
$
       \mathcal{ D }_{ 9 }
$
and
$
       \mathcal{ D }_{ 10 }
$
in Fig.\,\ref{Figure: D9-D10 }).

The kernel of a domino is allowed to evolve within the white convex shape
$
       \mathcal{ D }_{ n - 2 }.
$
It follows that the hull will be able to occupy the domain defined by the tipless diamond
$
        \mathcal{ S }_{ n }  \cap   \mathcal{D}_{ n + 2 }
$
where the four tips of $  \mathcal{ D }_{ n + 2 } $ (1--cell tip for $ n $ odd, 2--cell tip for $ n $ even) lying outside $ \mathcal{S}_{n} $ are truncated.
%
%
%
\begin{itemize}
%
%
\item  $ n $ odd :   
%
%
$ 
          \mathcal{D}_n 
$
has two perpendicular diagonals of length $ n + 2 $ and a center cell. 
It is the von Neumann neighborhood of range $ r $ at point $ ( 0, 0 ) $ 
(the center of $ \mathcal{D}_n $) and such that 
$ 
          2 r = n + 1.
$
$ \mathcal{D}_{ n  } $ is then the set of cells whose centers $ ( x, y ) $ are such that
$
          \{  ( x, y ) :  | \,  x  \, |  +   |  \, y  \ |  \le r \, \}  \  ( r \in \mathbb{N}  ) 
$
and so that both $ x $ and $ y $ are integers.
The cardinality of $ \mathcal{D}_{ n } $ is the centered square number 
$
          n_{ r  } =  2 r^2  + 2 r  + 1   
$
whence
%
%
\begin{equation}
\label{equation:  | D_{ n } |  -- n odd -- }
          |  \mathcal{D}_{ n } |   =   ( n^2 + 4 \,n + 5 ) / 2                                       
\end{equation}
%
%
thus
$
        \mid \mathcal{D}_{ n - 2 }\mid \   =     2 r^2 - 2 r + 1
$
and  
$
       \mid \mathcal{ S }_{ n }  \cap   \mathcal{D}_{ n + 2 }\mid  \, =   2 r^2  +  6  r  +  1 \, .
$
%
%
\item  $ n $ even :     
%
%
$ 
          \mathcal{D}_n 
$
has two perpendicular double diagonals of length $ n + 2 $ and a 4--cell center centered at point $ (0,0) $.
It is the Aztec diamond 
$ 
              \mathcal{AD}_{ r } 
$ 
of range $ r $ and such that
$ 
          2 r  =  n + 2.  
$
$ \mathcal{D}_{ n  }$ is then the set of cells whose centers $ ( x, y ) $ are such that
$
          \{  ( x, y ) :  | \,  x  \, |  +   |  \, y  \ |  \le r \, \}  \  ( r \in \mathbb{N}^\ast  ) 
$
and so that both $ x $ and $ y $ are half--integers.
It is the concatenation of 4 contiguous staircase quadrants, each containing ($ 1 + 2 + \ldots + r $) cells, then the cardinality of $ \mathcal{D}_{ n } $ is
$ 
            | \mathcal{AD}_{ r } |  =  n_{ r }^{ ( a ) }  =  2 r \,( r + 1 ) 
$ 
whence
%
%
\begin{equation}
\label{equation:  | D_{ n } |  -- n even -- }
          |  \mathcal{D}_{ n } |   =   ( n^2 + 6 \,n + 8 ) / 2                                       
\end{equation}
%
%
thus
$
        \mid \mathcal{D}_{ n - 2 }\mid \   =   2 r \, ( r - 1 )
$
and
$
       \mid \mathcal{ S }_{ n }  \cap   \mathcal{D}_{ n + 2 }\mid  \, =   2 r^2 + 6 r -4 \, .
$
\end{itemize}
%
%
The cardinality of the diamond satisfies the following relation
%
%
\begin{equation}
\label{equation:  | D_{ n } |  --  | D_{ n - 1 } | -- }
          |   \mathcal{D}_{ n } |  -  | \mathcal{D}_{ n - 1 } |   =   \left\{
                        \begin{array}{ll}
                                  1		& 	\mbox{ if }	n  \mbox{ odd } 	\\
                                 2 n + 3 	& 	\mbox{ if }	n  \mbox{ even }		 
                        \end{array}
                                                                                                \right.
\end{equation}
%
%
and it suffices to notice that $ n - 1 $ is even when $ n $ is odd --and vice versa-- then to adjust the cardinality of $ \mathcal{D}_{ n - 1 } $ accordingly.

 %
%
\begin{figure}
\centering
\includegraphics[width=8cm]{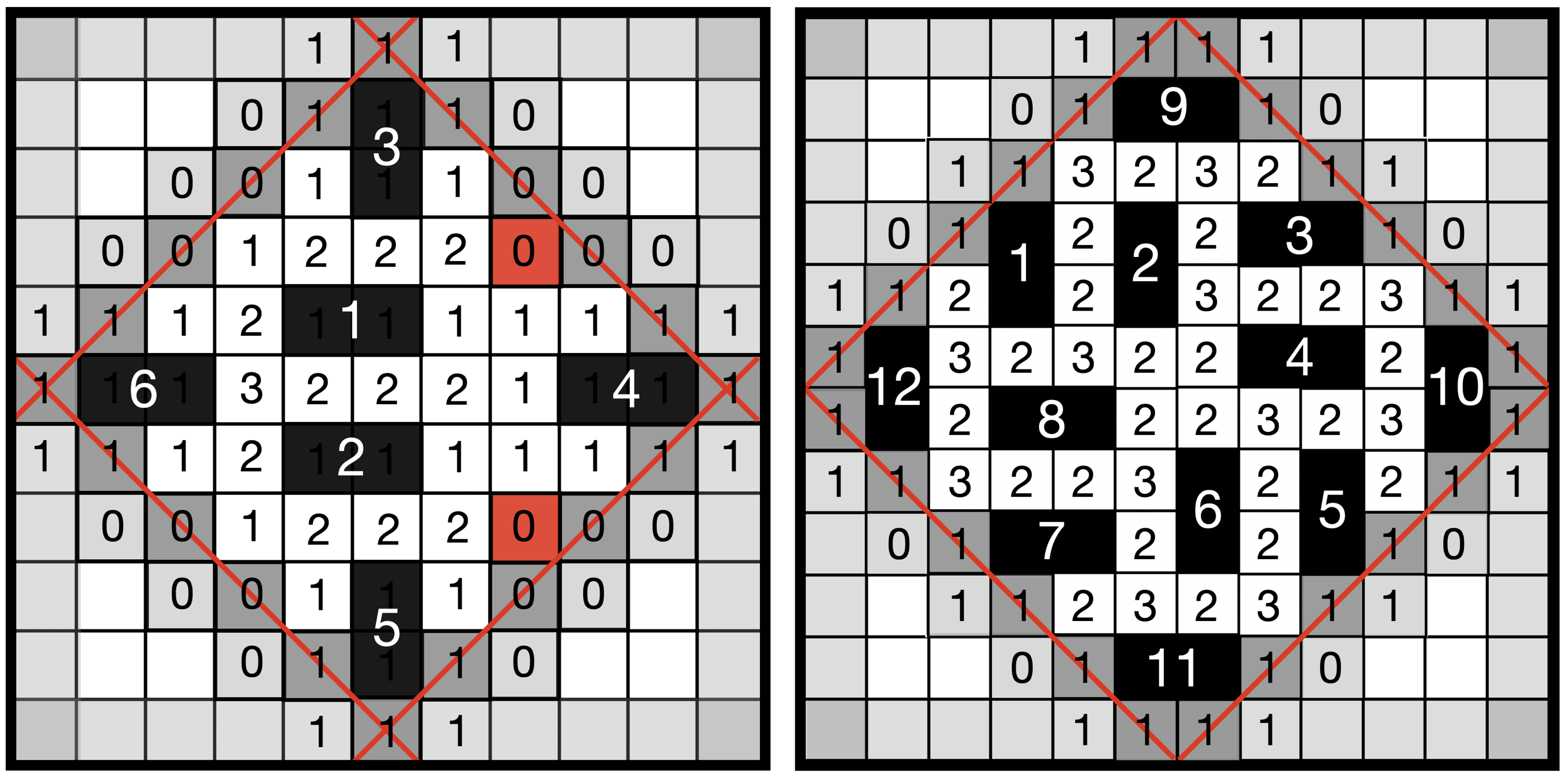}  
\caption{Connecting density and surface for a given population of dominos in $  \mathcal{ D}_{ 9 } $ and $  \mathcal{ D}_{ 10 } $.
Dominos are numbered in the range 
$
          k  \in   [  1, \psi_n  ] .
$
} 
\label{Figure: DOR-D9-D10 }
\end{figure}
%
%
%
Finally, it is useful to revisit the relationship in
(\ref{equation: The Square : Occupancy vs. Surface }) 
connecting a number of dominos and the surface of a shape from the occupancy index.
Applied to the diamond it will follow that
%
%
\begin{equation}
\label{equation: The Diamond : Occupancy vs. Surface }
          \sum_{ k = 1 }^{ \psi_n }   \rho \, ( k )    +   \varpi_n  \, = \
          \mid \mathcal{ S }_{ n }  \cap   \mathcal{D}_{ n + 2 }\mid   
\end{equation}
%
%
where $ \psi_n $ denotes the number of dominos for a given population in
$ 
          \mathcal{ D}_{n} 
$
 while
$
             \rho \, ( k )
$
and
$
             \varpi_n 
$
represent respectively the occupancy index of domino $ k $ and the number of lacunar voids.
As displayed in
Fig.\,\ref{Figure: DOR-D9-D10 }
%
%
%
\begin{itemize}
%
%
\item  in $  \mathcal{ D}_{ 9 } $ :   
$
           \rho ( 1 )  =  \rho ( 2 )  = 47 / 6 ; \ \rho ( 3 )  =  \rho ( 5 ) = 63 / 6 ;  \ \rho ( 4 )  =  72 / 6 ;  \ \rho ( 6 )  =  62 / 6 \\
$
\[
        \Rightarrow \sum_{ k = 1 }^{ 6 }   \rho \, ( k )  = 59 \ ; \hspace{4mm}      \varpi_{ 9 } =  2  +  4\, ( 2 + 3 )  =  22 \hspace{2mm}   \mbox{ and }
         \mid \mathcal{ S }_{ 9 }  \cap   \mathcal{D}_{ 11 }\mid  \, =  81
\]
%
%
%
\item  in $  \mathcal{ D}_{ 10 } $ :   
note the rotational symmetry for
$
           k \in (1,3,5,7)
$
and
$
           k \in (2, 4, 6, 8)
$
and
$
            k \in (9, 10, 11, 12).
$
Now 
$
           \rho ( 1 )  = 24 / 3 ; \ \rho ( 2 )  =  19 / 3 ;  \ \rho ( 9 )  =  29 / 3  \\
$
\[
        \Rightarrow \sum_{ k = 1 }^{ 12 }   \rho \, ( k )  = 96 \ ; \hspace{4mm}      \varpi_{ 10 } =  4 \times 2  =  8 \hspace{2mm}   \mbox{ and }
         \mid \mathcal{ S }_{ 10 }  \cap   \mathcal{D}_{ 12 }\mid  \, =  104
\]
\end{itemize}
%
%
%
and it will be observed in the sequel that among these two configurations, that of 
$  
           \mathcal{ D}_{ 10 } 
 $ 
 is maximal while that of 
 $  
           \mathcal{ D}_{ 9 } 
 $ 
 is not.
%
%
%
\subsection{The Construction Rule}                            
%
%
%
Unlike the dimer problem in statistical physics, which lists the number of all possible configurations for a given size shape, the more modest goal of this paper is to provide an {\em exact} value for the maximum number of dominos that can be contained in the diamond, or at least to provide {\em lower} and {\em upper} bounds.

We give here a heuristic approach. The analysis proceeds either according to an inductive formulation or according to a direct formulation. The inductive formulation makes it possible to underline certain recurrence relations. It will be observed that the geometry of the domino will involve a distribution of the patterns divided into  {\em six} classes according to the size $ n $ of the sample. In order to simplify the description, it will therefore be possible to proceed by “expansion” within the same class. The theoretical measures will split into six tables in the appendix, with a sample of up to $ n = 200 $. 

\begin{figure}
\centering
\includegraphics[width=9.5cm]{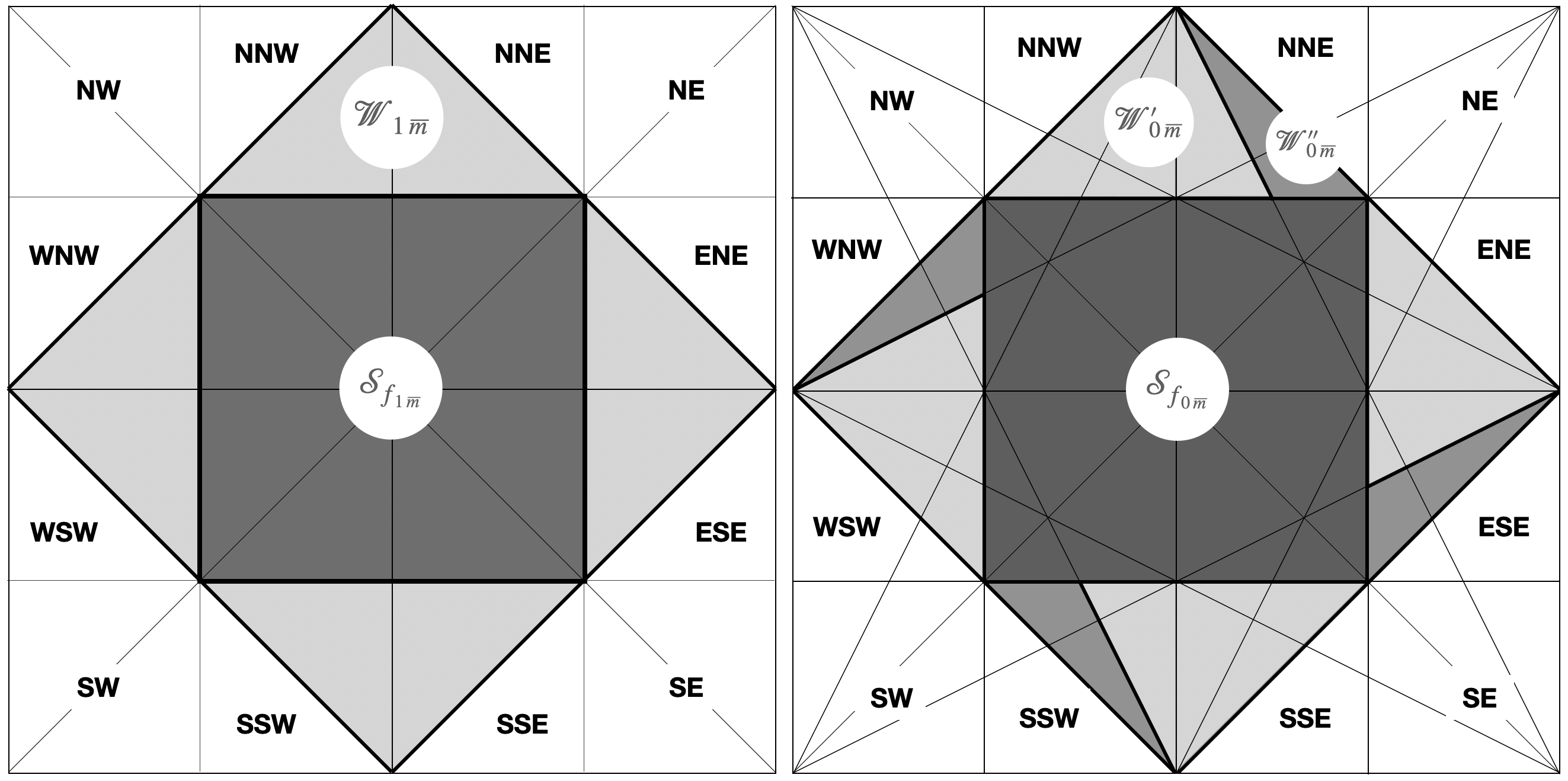}  
\caption{The (cellular) diamond
$
       \mathcal{D}_n
$
divided into regions.\newline
-- $ n $ odd ({\em left}):  a square central core
$
          \mathcal{S}_{ f_{ 1  \overline{m} }} 
$
of size
$
         f_{ 1  \overline{m} } \approx  m  
$
and four symmetric wedges
$
       \mathcal{W}_{ \, 1 \,  \overline{m}  }
$
in  \texttt{N}--\texttt{S}--\texttt{E}--\texttt{W} directions.\newline
-- $ n $ even ({\em right}):  a square central core
$
          \mathcal{S}_{ f_{ 0  \overline{m} }} 
$
of size
$
         f_{ 0  \overline{m} } \approx  m 
$
and four symmetric wedges
$
       \mathcal{W}_{ \, 0 \,  \overline{m}  }
$
in  \texttt{N}--\texttt{S}--\texttt{E}--\texttt{W} directions. Region
$
       \mathcal{W}_{ \, 0 \,  \overline{m}  }
$
splits into two subregions:
$
       \mathcal{W}'_{ \, 0 \,  \overline{m}  }
$
and
$
       \mathcal{W}''_{ \, 0 \,  \overline{m}  }.
$
} 
\label{Figure: Core-Wedges }
\end{figure}

Our domino arrangement follows a recursive scheme, organized around a square central ``{\em core}'' and four \texttt{N}--\texttt{S}--\texttt{E}--\texttt{W} triangular ``{\em wedges}'', symmetric by rotation as shown in 
Fig.\,\ref{Figure: Core-Wedges }.
The case $ n $ even is more intricate and the wedge must be divided into two subregions.

%
%
\begin{figure}
\centering
\includegraphics[width=11cm]{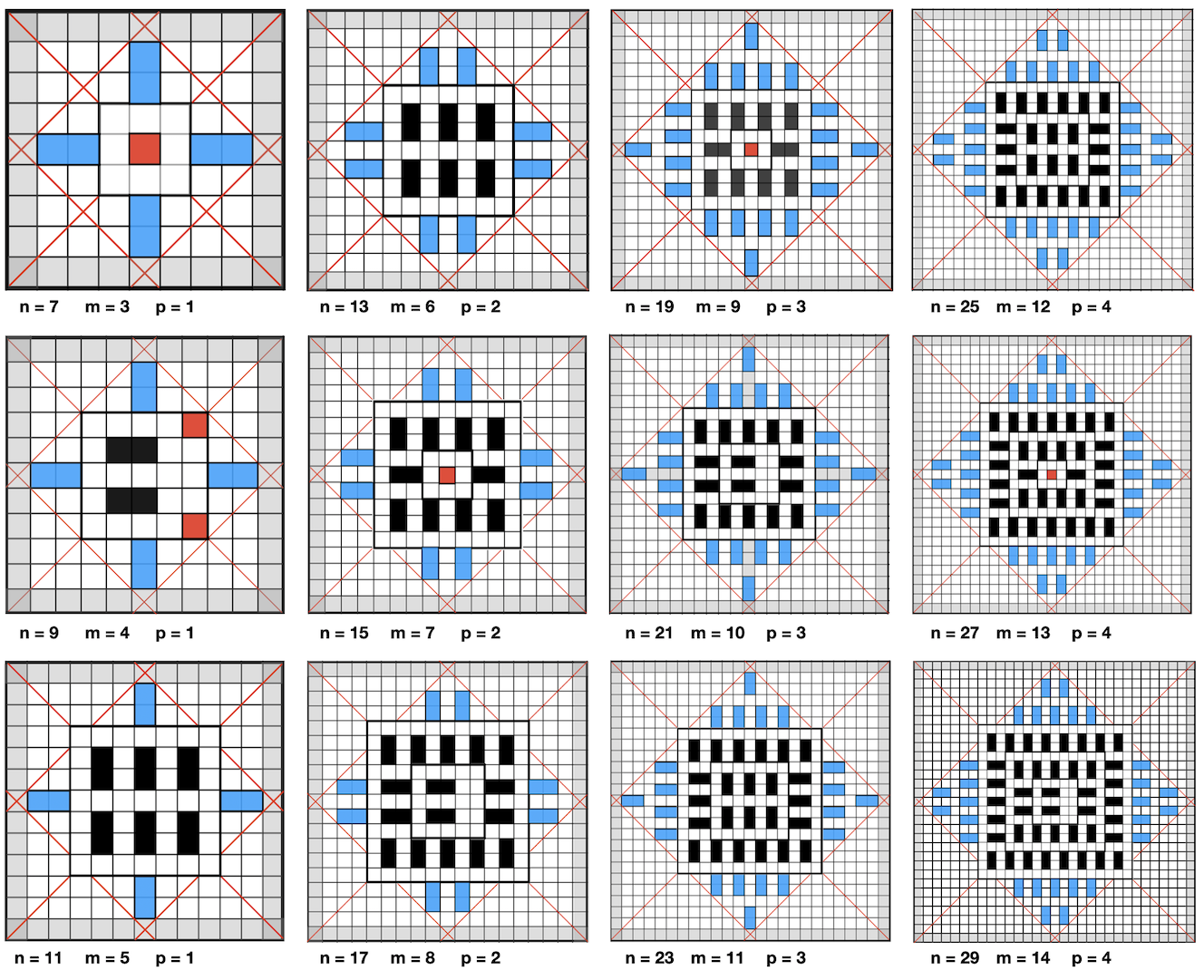}  
\caption{-- $n$ odd. The first configurations with small values of $ p \hspace{1mm} (p > 0) $. From top to bottom:
Class $ { 10 } $, 
Class $ { 11 } $,
Class $ { 12 } $.
Some lacunar voids are integral part of the core.\newline
} 
\label{Figure: C1-n7-n29 -- n odd } 
%
%
%
%
%
\centering
\includegraphics[width=11cm]{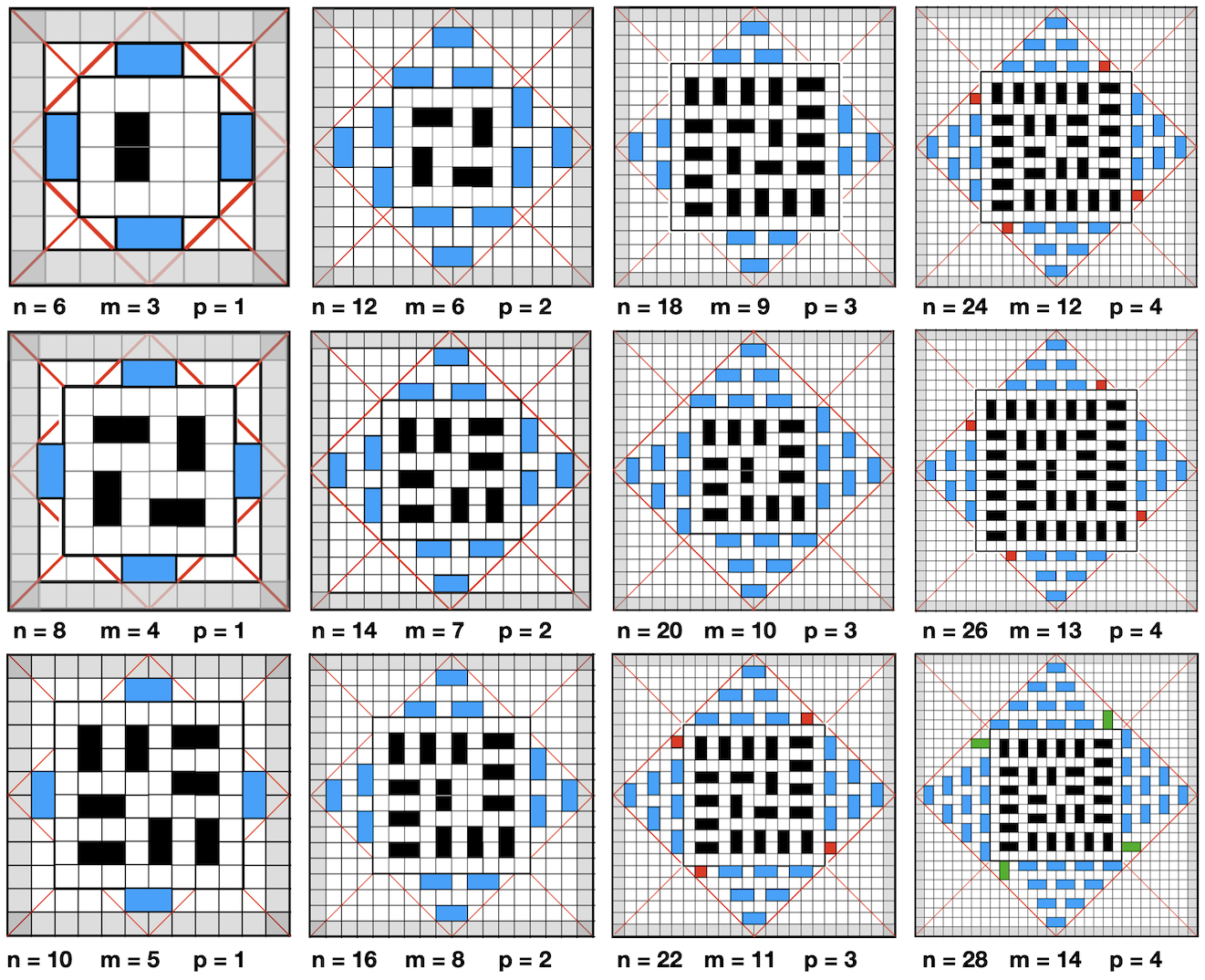}  
\caption{-- $n$ even. The first configurations with small values of $ p \hspace{1mm} (p > 0) $. From top to bottom:
Class $ { 00 } $, 
Class $ { 01 } $,
Class $ { 02 } $.
Some lacunar voids emerge as isolated cells from $ n > 20 $.
} 
\label{Figure: C0-n6-n28 -- n even } 
\end{figure}

%
%
The capacity, in number of dominos, of the {\em  square} shape and the configurations for different sizes of the core are evaluated in 
Sect.\,\ref{Section: Dominos in the Square}.
Setting 
$ 
             m = \lfloor n/2 \rfloor 
$ 
and 
$ 
             p = \lfloor m/3 \rfloor 
$ 
this hierarchy of configurations can then be divided into {\em  six} equivalence classes.
In the sequel, any function, say
$
        \mathcal{F},
$
 is indexed by the pair
 $
         { \,  \overline{n} \,  \overline{m}  } 
 $
 --\,often denoted as 
 $
        ({  \overline{n \, m}  }) 
 $
 or
  $
         n \, m 
 $
 for short\,-- where
 $
         \overline{n}
 $
stands for the class of $ n $  modulo 2 (its parity) and where
 $
         \overline{m}
 $
stands for the class of $ m $  modulo 3.
Whence the six families of configurations, namely
 $
          (00), (01), (02)
 $
for $ n $ {\em  even} and 
 $
          (10), (11), (12)
 $
for n {\em  odd}.
The first configurations for small values of $ n $ are displayed in
Figs. \ref{Figure: C1-n7-n29 -- n odd }--\ref{Figure: C0-n6-n28 -- n even }.
%
%
%
%
%

In this NP problem of domino arrangement, this arbitrary choice of construction around a central square whose side length
$
         f_{ \overline{n} \, \overline{m} } \approx  m 
$
is determined and where several solutions can coexist, the exact definition of this function must be considered as the {\em  axiom}, from which the placement of dominos follows a logical approach.

Besides, the case $ n $ even argues for a ``horizontal'' arrangement of the dominos with respect to the  \texttt{North} wedge (parallel to the adjacent  \texttt{North}  side of the core) while the case $ n $ odd argues for a ``vertical'' arrangement (perpendicular to the adjacent core’s side). In addition, a rotational symmetry of the wedges is required. 
We could dare the analogy, in a rather relative measure, with the observation of the authors\footnote{{\em in the four outer sub-regions, every tile lines up with nearby tiles, while in the fifth, central sub-region, differently-oriented tiles co-exist side by side.
-- \texttt{In}} ``Random Domino Tilings and the Arctic Circle Theorem'' 
\cite{Jockusch:Propp:Shor-1998}.
} except that theirs results from a physical phenomenon whereas ours results only from a rule of construction.

Finally, the placement of the dominos is organized from the left part of the (\texttt{North}) wedge, even if it leads to the completion of the inner right border with lacunar voids. Based on these criteria, our theoretical arrangement of dominos is {\em  unique}. The analysis is developed in
Sect.\,\ref{Section:Dominos in the Diamond -- n odd}--\ref{Section:Dominos in the Diamond -- n even}.
The penultimate column of the last six tables in the appendix
(Tabs.\,\ref{Table: Domino enumeration -- Class 10 }--\ref{Table: Domino enumeration -- Class 02 })
 lists the value $ \psi_n $ representing the {\em  lower} bound reached at the end of this construction.
%
%
%
\subsection{Injection of Disorder}                            
%
%
%
Our theoretical arrangement is almost optimal, not totally optimal. For example, it would be easy to observe in
Fig.\,\ref{Figure: C1-n7-n29 -- n odd }    
and in more detail in 
Fig.\,\ref{Figure: DOR-D9-D10 }    
that the diamond could accept at least $ 7 $ dominos for $ n = 9 $.
Nevertheless, by relaxing the constraint, it is possible to inject {\em  disorder} into the model in order to increase its population density. 
The counting of lacunar voids --\,emerging in 
Fig.\,\ref{Figure: C0-n6-n28 -- n even }
from  $ n > 20 $ in a systematic way\,-- will also be considered: in a disordered system, lacunar voids could indeed join together and thus form a new domino. In
Sect.\,\ref{Section: Injection of Disorder}, 
different scenarios of injection are examined on a case--by--case basis.
The aim is to achieve {\em  tight} upper bounds whenever possible.
The last column of
Tabs.\,\ref{Table: Domino enumeration -- Class 10 }--\ref{Table: Domino enumeration -- Class 02 }
 lists the value $  \overline{\psi}_n $ representing the {\em  upper} bound reached at the end of this injection.
In some cases, the lower and upper bounds coincide, so 
$  
          \overline{\psi}_n  =  \psi_n
$ 
and $ \psi_n $ becomes the {\em  exact} number of dominos for a maximum configuration.
%
%
%
\section{Dominos in the Square}
\label{Section: Dominos in the Square}
%

Given a square array 
$
       \mathcal{S}_n
$
of
$
       ( n + 2 ) \times ( n+ 2 )
$
cells including a border with perimeter
$
     4 n + 4
$
enclosing the 
$
     n \times n
$
field  
$
       \mathcal{ S }_{ n - 2 }  
$
of order $ n^2 $. 

\begin{figure}
\centering
\includegraphics[width=8.6cm]{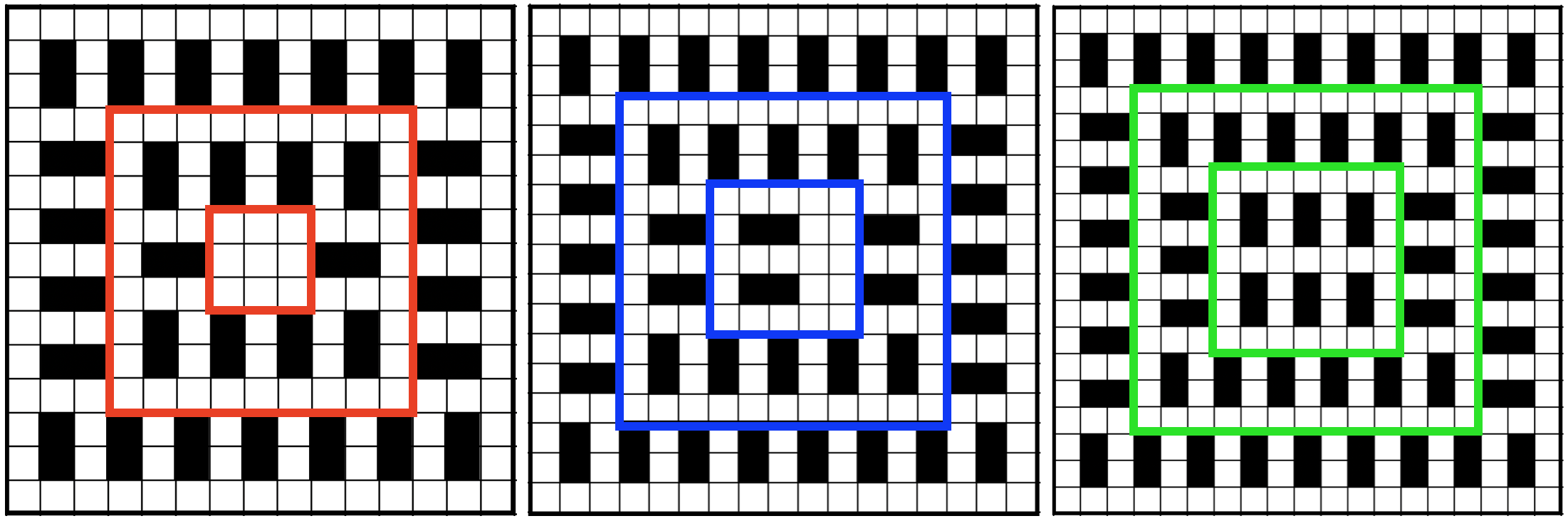}  
\includegraphics[width=8.6cm]{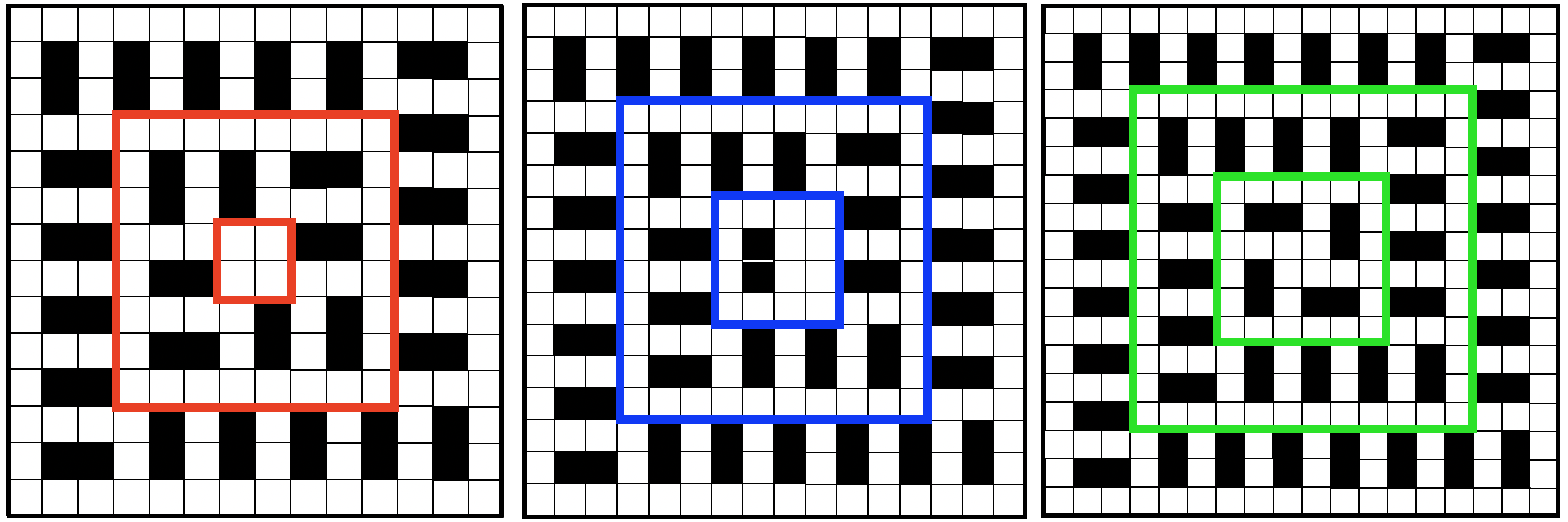}  
\caption{Dominos in the square $ \mathcal{S}_{n} $.
\newline $(\uparrow)$ The first configurations for $ n $ odd, from left to right:\newline
Class $ ( \dot{ 10 } ) $:
$
     ( \mathcal{S}_{1}, \mathcal{S}_{7}, \mathcal{S}_{13} ),
$
Class $ ( \dot{ 11 } ) $:
$
     ( \mathcal{S}_{3}, \mathcal{S}_{9}, \mathcal{S}_{15} ),
$
Class $ ( \dot{ 12 } ) $:
$
     ( \mathcal{S}_{5}, \mathcal{S}_{11}, \mathcal{S}_{17} ).
$
\newline $(\downarrow)$ The first configurations for $ n $ even, from left to right:\newline
Class $ ( \dot{ 00 } ) $:
$
     ( \mathcal{S}_{0}, \mathcal{S}_{6}, \mathcal{S}_{12} ),
$
Class $ ( \dot{ 01 } ) $:
$
     ( \mathcal{S}_{2}, \mathcal{S}_{8}, \mathcal{S}_{14} ),
$
Class $ ( \dot{ 02 } ) $:
$
     ( \mathcal{S}_{4}, \mathcal{S}_{10}, \mathcal{S}_{16} ).
$
} 
\label{Figure: Dominos in the Square -- n odd -- n even-- }
\end{figure}

The chosen arrangement is organized in concentric crowns surrounding a basic motif of minimal size. By construction, induced by the morphology of the domino, the successive crowns will be of thickness 3. It follows that the square shapes will be divided into {\em six} equivalence classes
 $
               \dot{ n } \,  \dot{ m }  
 $
--or shortly
$
           ( \dot{ nm })
$
--where
 $
               \dot{ n }  
 $
stands for the class of $ n $  modulo 2 (its parity) while
 $
               \dot{ m }  
 $
stands for the class of $ m $  modulo 3.
Whence the six families of configurations, namely
$
          ( \dot{ 10 } ), (  \dot{ 11 } ), (  \dot{ 12 } )
$
for $ n $ {\em  odd} and
$
          ( \dot{ 00 } ), (  \dot{ 01 } ), (  \dot{ 02 } )
$
for $ n $ {\em  even}. 
The first configurations for small values of $ n $ are displayed in
Fig.\,\ref{Figure: Dominos in the Square -- n odd -- n even-- }.
%
%
%
%
%

The maximal domino number
$
          \xi_n (\mathcal{S}_n)
$
covering 
$
       \mathcal{S}_n
$
is given by the inductive formula
\vspace{-3mm}
%
%
\begin{equation}
\label{equation: Induction Xi -- n odd }
     n  \hspace{2mm}  \mbox{odd:}   \hspace{5mm} 
      \xi_1 = 0,   \,\xi_3 = 2,  \, \xi_5 = 6   \hspace{5mm}   \mbox{and}   \hspace{5mm}     \xi_n  =  \xi_{n-6} + 2 ( n - 2 )
\vspace{-2mm}
\end{equation}
%
%
\vspace{-4mm}
%
%
\begin{equation}
\label{equation: Induction Xi -- n even }
     n  \hspace{2mm}  \mbox{even:}   \hspace{4mm} 
    \xi_0 = 0, \, \xi_2 = 1,  \, \xi_4 = 4   \hspace{5mm}   \mbox{and}   \hspace{5mm}     \xi_n  =  \xi_{n-6} + 2 ( n - 2 ) 
\vspace{-1mm}
\end{equation}
%
%
%
for
$
      n \ge 6
$
and where
$
         2 ( n - 2 )
$
denotes the number of dominos in the crown surrounding the inner subgrid 
$
    \mathcal{S}_{n-6}.
$
 A proof is straightforward: the crown 
 $
     \widehat{ \mathcal{S}}_{n}
$
of size
 $
    \mid \widehat{ \mathcal{S}}_{n}\mid = 4((n+1)+(n-1)+(n-3)) = 12(n-1)
$
yields a domino area (excluding border) with size 
$
     \mid \widehat{ \mathcal{S}}_{n}\mid \hspace{-1mm} - (4n+4) = 4(2n-4)
$
whence the number 
$
         ( 2n - 4 )
$
of black--white tetraminos.
 %
%

The exact expressions of $ \xi_n $ depend on the class of  $ \mathcal{S}_{ n } $ and are given hereafter.
%
%
\subsection{$n$ odd --- $\xi_n$}                               
%
%
Depending on the class of  
$
         m\bmod 3,
$
it follows from 
(\ref{equation: Induction Xi -- n odd })
that in
%
%
%
%
\begin{itemize}
%
%
%
\item  
$
           \mbox{ \texttt{Class} $ ( \dot{ 10 } ) $ : } m \equiv 0  \hspace{1mm}  \Rightarrow  \hspace{1mm}  n = 6\, p + 1  \hspace{2mm}  \mbox{and then}  
$
%
%
\begin{equation}
\label{equation: Xi_n -- n odd -- m = 3*p}
             \xi_n   =   \xi_1  + 2  \sum_{k=1}^{p} ( 6\, k - 1 )  =  2p \, ( 3p + 2 )=  \frac{ ( n - 1 ) ( n + 3 )}{6} 
\vspace{-3mm}
\end{equation}
%
%
%
%
\item  
$
           \mbox{ \texttt{Class} $ ( \dot{ 11 } ) $ : } m \equiv 1  \hspace{1mm}  \Rightarrow  \hspace{1mm}  n = 6\, p + 3  \hspace{2mm}  \mbox{and then}  
$
%
%
\begin{equation}
\label{equation: Xi_n -- n odd -- m = 3*p+1}
            \xi_n  =  \xi_3 + 2 \sum_{k=1}^{p} ( 6\, k + 1 )  = 2 + 2 p \, ( 3 p + 4 )  =  \frac{ ( n - 1 ) ( n + 3 )}{6} 
\vspace{-3mm}
\end{equation}
%
%
%
%
\item  
$
            \mbox{ \texttt{Class} $ ( \dot{ 12 } ) $ : } m \equiv 2  \hspace{1mm}  \Rightarrow  \hspace{1mm}  n = 6\, p + 5  \hspace{2mm}  \mbox{and then}  
$
%
%
\begin{equation}
\label{equation: Xi_n -- n odd -- m = 3*p+2}
              \xi_n  =   \xi_5  + 2  \sum_{k=1}^{p} ( 6\, k + 3 )  = 6 + 2 p \, ( 3 p + 6 )  =  \frac{ ( n + 1 )^2 }{6} 
\vspace{-3mm}
\end{equation}
%
%
%
\end{itemize}
%
%
%
%
\subsection{$n$ even --- $\xi_n$}                            
%
%
Depending on the class of  
$
         m\bmod 3,
$
it follows from 
(\ref{equation: Induction Xi -- n even }) 
that in
%
%
%
%
\begin{itemize}
%
%
%
\item  
$
          \mbox{ \texttt{Class} $ ( \dot{ 00 } ) $ : } m \equiv 0  \hspace{1mm}  \Rightarrow  \hspace{1mm}  n = 6\, p  \hspace{2mm}  \mbox{and then}  
$
%
%
\begin{equation}
\label{equation: Xi_n -- n even -- m = 3*p}
             \xi_n   =   \xi_0  + 2  \sum_{k=1}^{p} ( 6\, k - 2 )  =  2 p \, ( 3p + 1 )=  \frac{ n \, ( n + 2 )}{6} 
\vspace{-3mm}
\end{equation}
%
%
%
%
\item  
$
            \mbox{ \texttt{Class} $ ( \dot{ 01 } ) $ : } m \equiv 1  \hspace{1mm}  \Rightarrow  \hspace{1mm}  n = 6\, p + 2  \hspace{2mm}  \mbox{and then}  
$
%
%
\begin{equation}
\label{equation: Xi_n -- n even -- m = 3*p+1}
            \xi_n  =  \xi_2  + 2 \sum_{k=1}^{p}  6\, k  = 1 + 2 p \, ( 3 p + 3 ) =  \frac{  n ( n + 2 ) - 2 }{6} 
\vspace{-3mm}
\end{equation}
%
%
%
%
\item  
$
            \mbox{ \texttt{Class} $ ( \dot{ 02 } ) $ : } m \equiv 2  \hspace{1mm}  \Rightarrow  \hspace{1mm}  n = 6\, p + 4  \hspace{2mm}  \mbox{and then}  
$
%
%
\begin{equation}
\label{equation: Xi_n -- n even -- m = 3*p+2}
              \xi_n  =   \xi_4  + 2  \sum_{k=1}^{p} ( 6\, k + 2 )  = 4 + 2 p \, ( 3 p + 5 )  =  \frac{  n ( n + 2 ) }{6} 
\vspace{-2mm}
\end{equation}
%
%
%
\end{itemize}
%
%
%
The 100 first values of $ \xi_n $ are displayed in 
Tabs.\,\ref{Table: Dominos in the Square --  n odd  }--\ref{Table: Dominos in the Square --  n even  }
in the appendix.
%
%
%
%
\section{Dominos in the Diamond --- $n$ odd}
\label{Section:Dominos in the Diamond -- n odd}
\begin{figure}
\centering
\includegraphics[width=7cm]{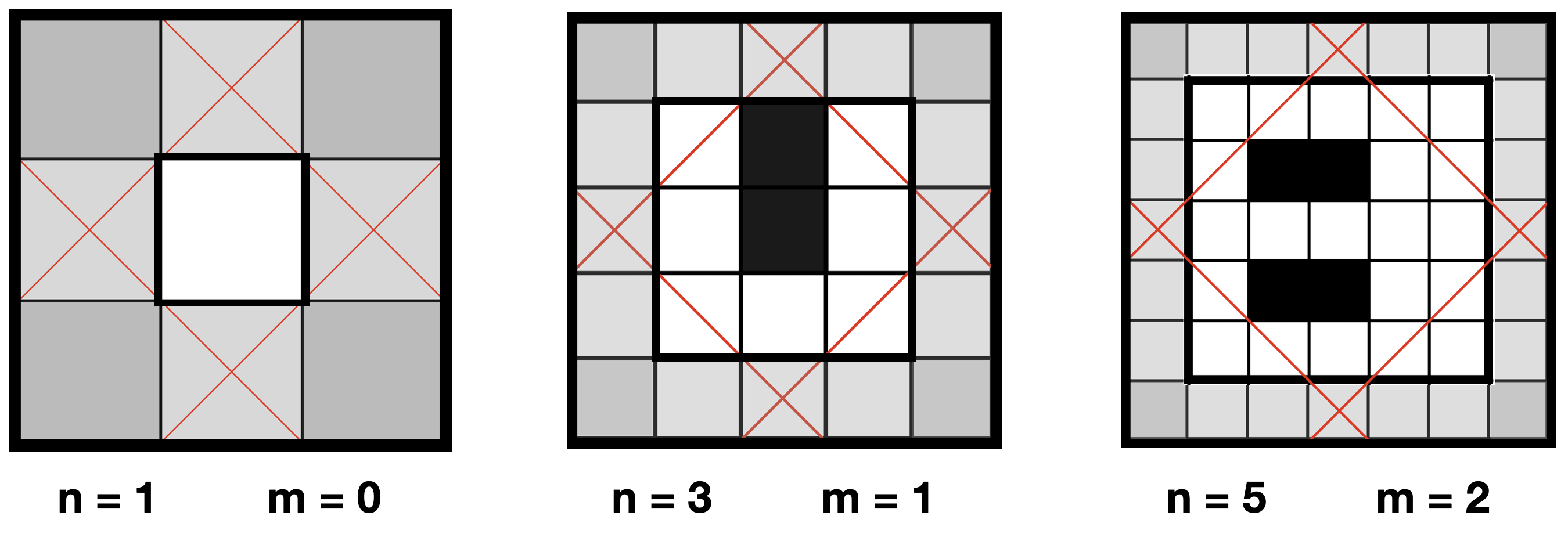}  
\caption{-- $n$ odd. The first configurations with $ p = 0 $. From left to right:
Class $ { 10 } $ ($  \psi_{ 1 }  = 0 $),
Class $ { 11 } $ ($  \psi_{ 3 }  = 1 $),
Class $ { 12 } $ ($  \psi_{ 5 }  = 2 $). \newline
} 
\label{Figure: C1m--P0 } 
\end{figure}

%
%
The following constructions for 
$
       \overline{ 1 \,m }
 $
 classes (often denoted as ``$ 1m $'' from now on) are valid only when the core--wedge structure actually appears from 
Fig.\,\ref{Figure: C1-n7-n29 -- n odd },
that is, for $ p \ge 1 $.
When growth formulas are used for the expansion of different regions, the strict condition $ p > 1 $  is required.
In the sequel, these conditions will be assumed everywhere, unless mentioned otherwise.
The first configurations with $ p = 0 $ are shown in
Fig.\,\ref{Figure: C1m--P0 }.
%
%
%
%
%

The case $ n $ odd argues for a ``vertical'' layout of the dominos (perpendicular to the adjacent core’s side). 
Referring back to
Fig.\,\ref{Figure: Core-Wedges } ({\em left})
the exact size 
$
         f_{ 1 m } \approx  m 
$
of the core will be adjusted so that the configuration of the wedges is entirely fixed by the $p$ parameter.
Let
$
             f_{ 1 m }^{+} ( m ) =   f_{ 1 m } ( m )  + 2   
$
be the extended length of the side of 
$
       \mathcal{S}_{  f_{ 1 m }(m)} 
$
including borders, and let  $ h_{1 m } $ such that
\[
             f_{ 1 m }^{+} ( m )  = n -  2  \cdot  h_{1 m } ( p )  \hspace{10mm}  ( n = 2 m + 1 )
\]
where $ h_{1 m } $ is the height of the vertical median column of Region
$
       \mathcal{W}_{ 1 m } .
$
By observing that this wedge is an isosceles right triangle, we fix its height so that
%
%
\begin{equation}
\label{equation:  n odd  --  2h_{ 1 m }( m )  }
           2  \cdot  h_{1 m } ( p )  =  3 p +  \mu_{ p } \hspace{2mm} \mbox{ with } \hspace{2mm}  \mu_{ p } = \left\{
                        \begin{array}{ll}
                                 0	& 	 \mbox{ for } p \mbox{ even }	\\	 
                                 1	& 	 \mbox{ for } p \mbox{ odd }	
                        \end{array}
                                                                                                \right.
\end{equation}
%
%
and it follows that
%
%
\[        
            f_{ 1 m }^{+} ( m ) =  \left\{
                        \begin{array}{ll}
                                 ( 6 p + 1 )  - ( 3 p +  \mu_{ p } ) = 3 p + 1 - \mu_{ p } 	& 	 \mbox{ in  \texttt{Class} $ 10 $ }  \\
                                 ( 6 p + 3 )  - ( 3 p +  \mu_{ p } ) = 3 p + 3 - \mu_{ p }	& 	 \mbox{ in  \texttt{Class} $ 11 $ }  \\
                                 ( 6 p + 5 )  - ( 3 p +  \mu_{ p } ) = 3 p + 5 - \mu_{ p } 	& 	 \mbox{ in  \texttt{Class} $ 12 $ }  
                        \end{array}
                                                                                                \right.
\]
%
%
whence finally the strict length of the side in the square core
%
%
\begin{equation}
\label{equation:  n odd  --  f_{ 1 m }( m ) -- C10 C11 C12 }
            f_{ 1 m }( m ) =   f_{ 1 m }^{+} ( m )  - 2 =   \left\{
                        \begin{array}{ll}
                                 m - 1 - \mu_{ p }		& 	 \mbox{ in  \texttt{Class} $ 10 $ }  \\
                                 m     - \mu_{ p }		& 	 \mbox{ in  \texttt{Class} $ 11 $ }  \\
                                 m + 1 - \mu_{ p }	& 	 \mbox{ in  \texttt{Class} $ 12 $ }  
                        \end{array}
                                                                                                \right.
\end{equation}
%
%
now expressed in terms of $ m $. 
Or, by using a general notation unifying the three classes
%
%
\begin{equation}
\label{equation:  n odd  --  f_{ 1 m }( m ) }
          f_{ 1 m }( m )   =  ( m - 1 +  m  \bmod 3  )  - \mu_{ p }	
\end{equation}
%
%
for any $ m > 0 $. Moreover, with
$
            \mu_{ p - 1 }   =  1 -  \mu_{ p } 	    
$
we get
%
%
\begin{equation}
\label{equation:  n odd  --  Square Core Elongation }
           f_{ 1 m }( m )  -  f_{ 1 m }( m - 3 ) =   3  -  \mu_{ p }  +   \mu_{ p - 1 }   =   4  -  2  \,  \mu_{ p }  =  \left\{
                        \begin{array}{ll}
                                  + 4	& 	 \mbox{ for } p \mbox{ even }	\\
                                  + 2	& 	 \mbox{ for } p \mbox{ odd }		 
                        \end{array}
                                                                                                \right.
\end{equation}
%
%
giving the elongation of the side of the square core within a
$
          p - 1  \rightarrow  p  
$
expansion.
%
%
\subsection{ \texttt{Class} $ 10 $ --- $ n $ odd \, \&  \, $ m \equiv  0  \pmod 3 $}
\label{Subsection: n  odd --- m = 3p -- }
%
%
%
The basic parameters of \texttt{Class} $ 10 $ verify: 
\[
		p \, \in  \, \mathbb{N} \,  ;  \hspace{3mm} m  \, = 3 \, p  \, ;  \hspace{3mm} n \, =  2 \, m \, + 1 \, ;
\]
%
%
\subsubsection{Capacity and expansion rate of the square core $  \mathcal{S}_{  f_{10}  } $ }                 
%
%
From
(\ref{equation:  n odd  --  f_{ 1 m }( m ) -- C10 C11 C12 }--\ref{equation:  n odd  --  f_{ 1 m }( m )  })
$
            f_{ 10 }( m ) =   m - 1 - \mu_{ p }	
$
and we get the following two cases
%
%
\begin{itemize}
%
%
\item  $  \mu_{ p }  = 0      \Rightarrow    f_{10} ( m ) =    m  - 1  $      
%
%
and 
$ 
          ( m - 1 ) / 2  \in \dot{2}   \hspace{1mm} \mbox{in} \hspace{1mm}    \mathbb{Z} / 3 \, \mathbb{Z}
$
then  $ \xi_{ m - 1 } $ follows from
(\ref{equation: Xi_n -- n odd -- m = 3*p+2})
whence
%
%
\begin{equation}
\label{equation: n odd  -- Xi-- m=3p -- mu=0 }
         \xi_{ m - 1 }   =  \frac{ m^2 }{6}   
\end{equation}
%
%
%
\item  $  \mu_{ p }  = 1      \Rightarrow    f_{10} ( m ) =    m  -  2  $       
%
%
and 
$ 
          ( m - 2 ) / 2  \in \dot{ 0 }   \hspace{1mm} \mbox{in} \hspace{1mm}    \mathbb{Z} / 3 \, \mathbb{Z}  
$
then  $ \xi_{ m - 2 } $ follows from
(\ref{equation: Xi_n -- n odd -- m = 3*p})
whence
%
%
\begin{equation}
\label{equation: n odd  -- Xi-- m=3p -- mu=1 }
         \xi_{ m - 2 }   =  \frac{ ( m - 3 ) ( m + 1 ) }{6} 
\end{equation}
%
%
\end{itemize}
%
%
giving the capacity of the square core. 
In the same way, from
(\ref{equation:  n odd  --  Square Core Elongation })
we get the following two cases
%
%
\begin{itemize}
%
%
\item  $  \mu_{ p }  = 0      \Rightarrow    f_{10} ( m - 3 ) =  f_{10} ( m ) - 4  =  m - 5    $     
%
%
and 
$ 
          (m - 5 )/2  \in \dot{0}  
$
then  $ \xi_{ m - 5 } $ follows from
(\ref{equation: Xi_n -- n odd -- m = 3*p})
whence
%
%
\begin{equation}
\label{equation: n odd  -- DXi-- m=3p -- mu=0 }
         \xi_{ m - 1 }  -  \xi_{ m - 5 }  =  \frac{ m^2 }{6}  -  \frac{ ( m - 6 ) ( m - 2 ) }{6}  = \frac{ 4 \, m - 6 }{3}   =  4 \, p  - 2  
\end{equation}
%
%
%
\item  $  \mu_{ p }  = 1      \Rightarrow     f_{10} ( m - 3 ) =  f_{10} ( m ) - 2  =  m - 4  $    
%
%
and 
$ 
          (m - 4 )/2  \in \dot{2}  
$
then  $ \xi_{ m - 4 } $ follows from
(\ref{equation: Xi_n -- n odd -- m = 3*p+2})
whence
%
%
\begin{equation}
\label{equation: n odd  -- DXi-- m=3p -- mu=1 }
         \xi_{ m - 2 }  -  \xi_{ m - 4 }  =  \frac{ ( m - 3 ) ( m + 1 ) }{6}  -  \frac{ ( m - 3 )^2 }{6}  = \frac{ 2 \, m - 6 }{3}   =  2 \, p  - 2  
\end{equation}
%
%
\end{itemize}
%
%
giving the expansion rate for the capacity in the square core.
%
%
\subsubsection{Capacity and expansion rate of the wedge $  \mathcal{W}_{ 10 } $ }                 
%
%
Let $ W_ { 10 }( p ) $ be the capacity of the wedge for any given $ p $. From
(\ref{equation:  n odd  --  2h_{ 1 m }( m )  })
comes
\[
              2 \, ( \,  h_{ 10 } ( p ) +    \mu_{ p }  )  =  3 \, (  \, p +  \mu_{ p }  )  
\]
and thus we get the capacity of the vertical median of Region
$
       \mathcal{W}_{10} 
$
\[
              \overline{h}_{10} ( p ) =  (  h_{10} ( p ) +  \mu_{ p } ) / 3  = ( p +  \mu_{ p }  ) / 2
\]
in number of vertical dominos --namely, in number of horizontal rows of dominos.
The baserow of 
$
       \mathcal{W}_{10}
$
has the length
\[
              b_{10} ( p ) =  2 \cdot h_{10} ( p ) - 1 =   3 p  +  \mu_{ p } - 1  \hspace{3mm}    
\]
always odd. This baserow can accommodate a binary string of the form
$
         (01)^*0
$
of length $  b_{10} ( p ) $ with 
$
          h_{10} ( p )
$
zeros and thereby a baserow with capacity
\[
            \overline{b}_{10} ( p )  = ( b_{10} ( p ) - 1 ) / 2  =  ( 3 p  +  \mu_{ p } - 2 ) / 2   
\]
in number of vertical dominos. Now
$
            \overline{h}_{10} ( 1 ) =  1  \hspace{1mm}  \mbox{;} \hspace{2mm}  \overline{b}_{10} ( 1 ) =  1
$
and
\[
            \overline{h}_{10} ( p )  -  \overline{h}_{10} ( p - 1 ) =  \mu_{ p }
             \hspace{2mm}  \mbox{and} \hspace{2mm}  
             \overline{b}_{10} ( p )  -  \overline{b}_{10} ( p - 1 ) =  1 + \mu_{ p }    \hspace{3mm}   
\]
whence the two following cases:
%
%
\begin{itemize}
%
%
\item  $  \mu_{ p }  = 0  \Rightarrow  $
%
%
The number of rows remains unchanged during expansion, and the capacity of the baserow is increased by 1, as well as, by extension, the capacity of the $ p / 2 $ rows. Then
$
             W_{10}( p ) - W_{10}( p - 1 ) = p / 2 .
$
%
%
\item  $  \mu_{ p }  = 1  \Rightarrow  $
%
%
The capacity of the baserow is increased by 2, as well as, by extension, the capacity of the $ ( p - 1 ) / 2 $ rows  $ ( p > 1 ) $.The capacity of the median is increased by 1, by adding a new row with a single domino at the top of the wedge. Then
$
             W_{10}( p ) - W_{10}( p - 1 ) = 2 \, ( p - 1 ) / 2 + 1 = p  .
$
\end{itemize}
%
%
It follows that 
%
%
\begin{equation}
\label{equation: n odd  -- W_10 expansion -- m=3p }
               W_{10}( 0 )  = 0  \hspace{1mm} \mbox{;} \hspace{2mm}  W_{10}( p ) - W_{10}( p - 1 ) =  \frac{ p \, ( 1 + \mu_{ p } ) }{ 2 } 
\end{equation}
%
%
yields the 
$
          p - 1  \rightarrow  p  
$
expansion of Region
$
         \mathcal{W}_{10}.
$
Nevertheless, since 
$
        \mu_{ p } + \mu_{ p - 1 }  =  1
$
the following expansion
\[
           \overline{h}_{10} ( p )  -  \overline{h}_{10} ( p - 2 )  =  1
            \hspace{2mm}  \mbox{and} \hspace{2mm}  
            \overline{b}_{10} ( p )  -  \overline{b}_{10} ( p - 2 ) =   3
\]
is relevant and does not depend on the parity of $ p $. As a result, the vertical median can accommodate {\em one} additional baserow extended with {\em three} additional dominos within any
$
          p - 2  \rightarrow  p 
$
expansion. For clarity's sake, we distinguish the two cases
%
%
\begin{itemize}
%
%
\item  $  \mu_{ p }  = 0 :    W_{10}( 2 ) =  \overline{b}_{10} ( 2 )  = 2 $  ;  
%
%
$
              W_{10}( p ) =   W_{10}( p - 2 ) +  ( 3 \, p - 2 ) / 2
$
and after the appropriate change of variable $ p = 2 k $
            %
              %
\[
             W_{10}( p )  =   \sum_{ k = 1 }^{ p / 2 } ( 3 \, k - 1 )      
            =    \frac{ p \, ( 3 \, p + 2 ) }{ 8 } 
\]
%
%
\item  $  \mu_{ p }  = 1 :    W_{10}( 1 ) =  \overline{b}_{10} ( 1 )  = 1 $  ;
%
%
$
              W_{10}( p ) =   W_{10}( p - 2 ) +  ( 3 \, p - 1 ) / 2
$
and after the appropriate change of variable $ p = 2 k - 1 $
            %
             %
\[
             W_{10}( p )  =   \sum_{ k = 1 }^{ ( p + 1 ) / 2 } ( 3 \, k - 2 )      
             =   \frac{ ( p + 1 ) ( 3 \, p + 1 ) }{ 8 } 
\]

\end{itemize}
%
%
and finally
%
%
\begin{equation}
\label{equation: n odd  -- expression of W_10 -- m=3p }
                W_{10}( p )  = \frac{ ( p +  \mu_{ p } ) \, ( 3 \, p + 2 -  \mu_{ p } ) }{ 8 } 
\end{equation}
%
%
or even in the form
%
%
\begin{equation}
\label{equation: n even  -- W_10 -- m=3p -- analogy with triangle area }
                W_{10}( p )  = \frac{   \overline{h}_{10} ( p ) \cdot ( \overline{b}_{10} ( p ) + ( 2 - \mu_{ p } ) ) }{ 2 } 
\end{equation}
%
%
in order to emphasize that  $ W_{10} $ grows with the area of Triangle
$
         \mathcal{W}_{10}.
$
%
%
\begin{figure}
\centering
\includegraphics[width=8cm]{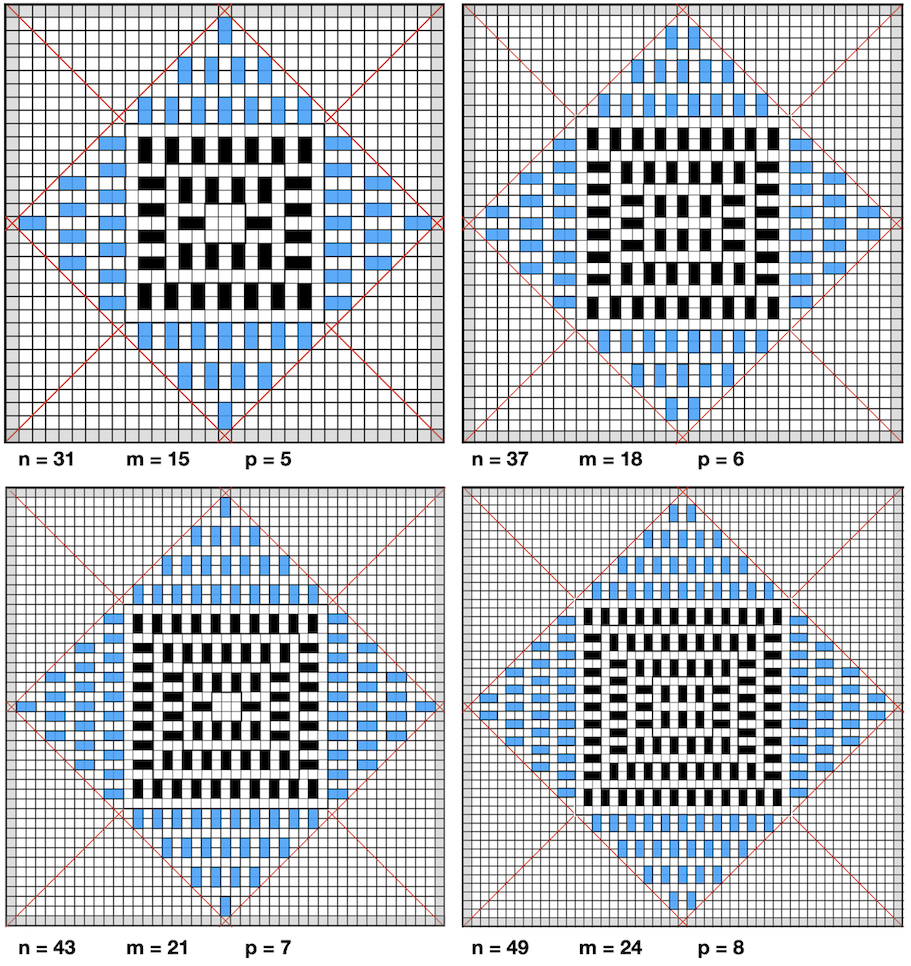}  
\caption{A 4--sequence of  $ \mathcal{D}_{ n } $ expansion in \texttt{Class} $ 10 $.
The class of the square core 
$
       \mathcal{S}_{ f_{10}( m )} 
$
alternates from $ ( \dot{ 10 } ) $ to $ ( \dot{ 12 } ) $ 
(refer to Fig.\,\ref{Figure: Dominos in the Square -- n odd -- n even-- })
whereas a new row is added in
$
         \mathcal{W}_{10}
$
for $ p $ odd.
} 
\label{Figure: C10-N31N49 }
\end{figure}
%
%
%
\subsubsection{Capacity and expansion  of  $ \mathcal{D}_{ n } $ in \texttt{Class} $ 10 $ }                 
%
%
The capacity $ \psi_{ n } $ of  $ \mathcal{D}_{ n } $ in \texttt{Class} $ 10 $ is simply given by
%
%
\begin{equation}
\label{equation: Psi_n in Class 10 }
      \psi_{n}  =   \xi_{ f_{10}( m ) }  + 4 \, W_{10}( p ) 
\vspace{-4mm}
\end{equation}
%
%
and then
%
%
\begin{itemize}
%
%
\item  $  \mu_{ p }  = 0      \Rightarrow   $ from 
%
%
(\ref{equation: n odd  -- Xi-- m=3p -- mu=0 })
and from
(\ref{equation: n odd  -- expression of W_10 -- m=3p })
%
%
\[        
             \psi_{n}  =  \frac{ m^2 }{6}  +  \frac{  p  \, ( 3 \, p + 2 ) }{ 2 }  
\]
%
%
\vspace{-5mm}
%
%
\item  $  \mu_{ p }  = 1      \Rightarrow    $ from   
%
%
(\ref{equation: n odd  -- Xi-- m=3p -- mu=1 })
and from
(\ref{equation: n odd  -- expression of W_10 -- m=3p })
%
%
%
\[ 
         \psi_{n}  =   \frac{ ( m - 3 ) ( m + 1 ) }{6}  +  \frac{   ( p + 1 )  \, ( 3 \, p + 1 ) }{ 2 }   
\]
%
%
\vspace{-4mm}
\end{itemize}
%
%
\vspace{-4mm}
and finally
%
%
\begin{equation}
\label{equation: n odd  -- PSi in C10 }
        \psi_{n}  =   p \, ( 3p + 1 )  =   \frac{ ( n - 1 ) ( n + 1 ) }{ 12 } 
\end{equation}
%
%
whatever the parity of $ p $. 

The expansion rate of  $ \psi_{ n } $ can also be obtained stepwise:
it follows from 
(\ref{equation: n odd  -- DXi-- m=3p -- mu=0 }--\ref{equation: n odd  -- DXi-- m=3p -- mu=1 })
and
from (\ref{equation: n odd  -- W_10 expansion -- m=3p })
that is, with $ p > 0 $  
\[
           \psi_{n} - \psi_{ n - 6 } =  ( \xi_{ f_{10}( m )}  - \xi_{ f_{10}( m - 3 )} )  +  4 \, ( W_{10}( p ) - W_{10}( p - 1 )  ) 
\]
\vspace{-5mm}
\[
            =    \mu_{ p }  \, ( 2 \, p - 2 )  +  ( 1 -  \mu_{ p } ) \, ( 4 \, p - 2 )   + 2 \, p \, ( 1 +  \mu_{ p } )  =  6 \, p - 2         
\]
now 
$ 
           n = 2\,m + 1 
$ 
and 
$
          m = 3\,p
$
then
%
%
\begin{equation}
\label{equation: n odd -- m == 0  -- Psi(n) - Psi( n - 6 ) }
          \psi_{ 1 }  = 0  \ ; \hspace{2mm}  \psi_{n} - \psi_{ n - 6 } =  n - 3  \, .
\end{equation}
%
%
%
The values of $ \psi_{n} $ in \texttt{Class} $ 10 $ are displayed in 
Table \ref{Table: Domino enumeration -- Class 10 }.
A 4--sequence of expansions is displayed in 
Fig.\,\ref{Figure: C10-N31N49 }.
%
%
%
%
%
%
%
%
%
%
\subsection{ \texttt{Class} $ 11 $ --- $ n $ odd \, \&  \, $ m \equiv  1  \pmod 3 $}
\label{Subsection: n  odd --- m = 3p + 1 -- }
%
%
%
The basic parameters of \texttt{Class} $ 11 $ verify: 
\[
		p \, \in  \, \mathbb{N} \,  ;  \hspace{3mm} m  \, = 3 \, p + 1  \, ;  \hspace{3mm} n \, =  2 \, m \, + 1 \, ;
\]
%
%
\subsubsection{Capacity and expansion rate of the square core $  \mathcal{S}_{  f_{11}  } $ }                 
%
%
From
(\ref{equation:  n odd  --  f_{ 1 m }( m ) -- C10 C11 C12 }--\ref{equation:  n odd  --  f_{ 1 m }( m )  })
$
            f_{ 11 }( m ) =   m  - \mu_{ p }	
$
and we get the following two cases
%
%
\begin{itemize}
%
%
\item  $  \mu_{ p }  = 0      \Rightarrow    f_{ 11 } ( m ) =    m    $      
%
%
and 
$ 
         m / 2  \in \dot{ 0 }   \hspace{1mm} \mbox{in} \hspace{1mm}    \mathbb{Z} / 3 \, \mathbb{Z}
$
then  $ \xi_{ m } $ follows from
(\ref{equation: Xi_n -- n odd -- m = 3*p})
whence
%
%
\begin{equation}
\label{equation: n odd  -- Xi-- m=3p+1 -- mu=0 }
         \xi_{ m }   =  \frac{ ( m - 1 ) ( m + 3 ) }{6}   
\end{equation}
%
%
%
\item  $  \mu_{ p }  = 1      \Rightarrow    f_{ 11 } ( m ) =    m  -  1  $       
%
%
and 
$ 
          ( m - 1 ) / 2  \in \dot{ 1 }   \hspace{1mm} \mbox{in} \hspace{1mm}    \mathbb{Z} / 3 \, \mathbb{Z}  
$
then  $ \xi_{ m - 2 } $ follows from
(\ref{equation: Xi_n -- n odd -- m = 3*p+1})
whence
%
%
\begin{equation}
\label{equation: n odd  -- Xi-- m=3p+1 -- mu=1 }
         \xi_{ m - 1 }   =  \frac{ ( m - 2 ) ( m + 2 ) }{6} 
\end{equation}
%
%
\end{itemize}
%
%
giving the capacity of the square core. 
In the same way, from
(\ref{equation:  n odd  --  Square Core Elongation })
we get the following two cases
%
%
\begin{itemize}
%
%
\item  $  \mu_{ p }  = 0      \Rightarrow    f_{ 11 } ( m - 3 ) =  f_{ 11 } ( m ) - 4  =  m - 4    $     
%
%
and 
$ 
          (m - 4 )/2  \in \dot{ 1 }  
$
then  $ \xi_{ m - 4 } $ follows from
(\ref{equation: Xi_n -- n odd -- m = 3*p+1})
whence
%
%
\begin{equation}
\label{equation: n odd  -- DXi-- m=3p+1 -- mu=0 }
         \xi_{ m }  -  \xi_{ m - 4 }  =  \frac{ ( m - 1 ) ( m + 3 ) }{6}  -  \frac{ ( m - 5 ) ( m - 1 ) }{6}  = \frac{ 4 \, m - 4 }{ 3 }   =  4 \, p    
\end{equation}
%
%
%
\item  $  \mu_{ p }  = 1      \Rightarrow     f_{ 11 } ( m - 3 ) =  f_{ 11 } ( m ) - 2  =  m - 3  $     
%
%
and 
$ 
          (m - 3 )/2  \in \dot{ 0 }  
$
then  $ \xi_{ m - 3 } $ follows from
(\ref{equation: Xi_n -- n odd -- m = 3*p})
whence
%
%
\begin{equation}
\label{equation: n odd  -- DXi-- m=3p+1 -- mu=1 }
         \xi_{ m - 1 }  -  \xi_{ m - 3 }  =  \frac{ ( m - 2 ) ( m + 2 ) }{6}  -  \frac{  m ( m - 4 ) }{6}  = \frac{ 2 \, m - 3 }{3}   =  2 \, p    
\end{equation}
%
%
\end{itemize}
%
%
giving the expansion rate for the capacity in the square core.
%
%
%
\begin{figure}
\centering
\includegraphics[width=8cm]{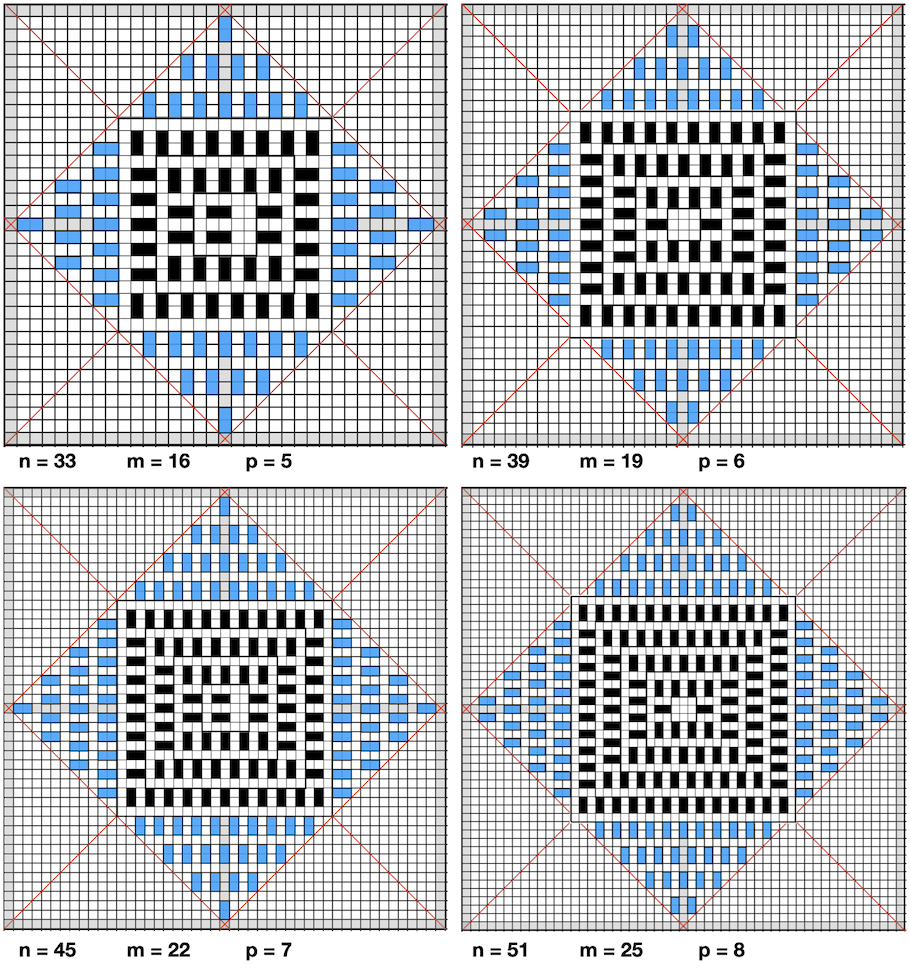}  
\caption{A 4--sequence of  $ \mathcal{D}_{ n } $ expansion in \texttt{Class} $ 11 $.
The class of the square core 
$
       \mathcal{S}_{ f_{11}( m )} 
$
alternates from $ ( \dot{ 11 } ) $ to $ ( \dot{ 10 } ) $ 
(refer to Fig.\,\ref{Figure: Dominos in the Square -- n odd -- n even-- })
whereas a new row is added in
$
         \mathcal{W}_{ 11 }
$
for $ p $ odd.
} 
\label{Figure: C11-N33N51 }
\end{figure}
%
%
%
%
\subsubsection{Capacity and expansion  of  $ \mathcal{D}_{ n } $ in \texttt{Class} $ 11 $ }                 
%
%
The capacity $ \psi_{ n } $ of  $ \mathcal{D}_{ n } $ in \texttt{Class} $ 11 $ is given by
%
%
\begin{equation}
\label{equation: Psi_n in Class 11 }
      \psi_{n}  =   \xi_{ f_{ 11 }( m ) }  + 4 \, W_{ 11 }( p ) 
\end{equation}
%
%
and with
$
          W_{ 11 } \equiv W_{ 10 }
$
we get
%
%
\begin{itemize}
%
%
\item  $  \mu_{ p }  = 0      \Rightarrow   $ from 
%
%
(\ref{equation: n odd  -- Xi-- m=3p+1 -- mu=0 })
and from
(\ref{equation: n odd  -- expression of W_10 -- m=3p })
%
%
\[        
             \psi_{n}  =  \frac{ ( m - 1 ) ( m + 3 ) }{6}   +  \frac{  p  \, ( 3 \, p + 2 ) }{ 2 }  
\]
%
%
%
\item  $  \mu_{ p }  = 1      \Rightarrow    $ from   
%
%
(\ref{equation: n odd  -- Xi-- m=3p+1 -- mu=1 })
and from
(\ref{equation: n odd  -- expression of W_10 -- m=3p })
%
%
%
\[ 
         \psi_{n}  =   \frac{ ( m - 2 ) ( m + 2 ) }{6}  +  \frac{   ( p + 1 )  \, ( 3 \, p + 1 ) }{ 2 }   
\]
%
%
\end{itemize}
%
%
and finally
%
%
\begin{equation}
\label{equation: n odd  -- PSi in C11 }
        \psi_{n}  =   3 p \, ( p + 1 )  =   \frac{ ( n - 3 ) ( n + 3 ) }{ 12 }  
\end{equation}
%
%
whatever the parity of $ p > 0 $.

The expansion rate of  $ \psi_{ n } $ can also be obtained stepwise:
it follows from 
(\ref{equation: n odd  -- DXi-- m=3p+1 -- mu=0 }--\ref{equation: n odd  -- DXi-- m=3p+1 -- mu=1 })
and
from (\ref{equation: n odd  -- W_10 expansion -- m=3p })
that is, with $ p > 0 $  
\[
           \psi_{n} - \psi_{ n - 6 } =  ( \xi_{ f_{ 11 }( m )}  - \xi_{ f_{ 11 }( m - 3 )} )  +  4 \, ( W_{ 11 }( p ) - W_{ 11 }( p - 1 )  ) 
\]
\[
            =    2 \, p  \cdot  \mu_{ p }  \,  +  4 \, p \, ( 1 -  \mu_{ p } ) \,  + 2 \, p \, ( 1 +  \mu_{ p } )  =  6 \, p         
\]
now 
$ 
           n = 2\,m + 1 
$ 
and 
$
          m = 3\,p + 1
$
then
%
%
\begin{equation}
\label{equation: n odd -- m == 1  -- Psi(n) - Psi( n - 6 ) }
          \psi_{ 1 }  = 1  \ ; \hspace{2mm}  \psi_{n} - \psi_{ n - 6 } =  n - 3  \, .
\end{equation}
%
%
%
The values of $ \psi_{n} $ in \texttt{Class} $ 11 $ are displayed in 
Table \ref{Table: Domino enumeration -- Class 11 }.
A 4--sequence of expansions is displayed in 
Fig.\,\ref{Figure: C11-N33N51 }.
%
%
%
%
%
%
%
%
%
%
%
\subsection{ \texttt{Class} $ 12 $ --- $ n $ odd \, \&  \, $ m \equiv  2  \pmod 3 $}
\label{Subsection: n  odd --- m = 3p + 2 -- }
%
%
%
The basic parameters of \texttt{Class} $ 12 $ verify: 
\[
		p \, \in  \, \mathbb{N} \,  ;  \hspace{3mm} m  \, = 3 \, p + 2  \, ;  \hspace{3mm} n \, =  2 \, m \, + 1 \, ;
\]
%
%
\subsubsection{Capacity and expansion rate of the square core $  \mathcal{S}_{  f_{ 12 }  } $ }                 
%
%
From
(\ref{equation:  n odd  --  f_{ 1 m }( m ) -- C10 C11 C12 }--\ref{equation:  n odd  --  f_{ 1 m }( m )  })
$
            f_{ 12 }( m ) =   m  + 1 - \mu_{ p }	
$
and we get the following two cases
%
%
\begin{itemize}
%
%
\item  $  \mu_{ p }  = 0      \Rightarrow    f_{ 12 } ( m ) =    m  + 1   $      
%
%
and 
$ 
         ( m + 1 ) / 2  \in \dot{ 1 }   \hspace{1mm} \mbox{in} \hspace{1mm}    \mathbb{Z} / 3 \, \mathbb{Z}
$
then  $ \xi_{ m + 1 } $ follows from
(\ref{equation: Xi_n -- n odd -- m = 3*p+1})
whence
%
%
\begin{equation}
\label{equation: n odd  -- Xi-- m=3p+2 -- mu=0 }
         \xi_{ m + 1 }   =  \frac{ m ( m + 4 ) }{6}   
\end{equation}
%
%
%
\item  $  \mu_{ p }  = 1      \Rightarrow    f_{ 12 } ( m ) =    m    $       
%
%
and 
$ 
          m / 2  \in \dot{ 2 }   \hspace{1mm} \mbox{in} \hspace{1mm}    \mathbb{Z} / 3 \, \mathbb{Z}  
$
then  $ \xi_{ m } $ follows from
(\ref{equation: Xi_n -- n odd -- m = 3*p+2})
whence
%
%
\begin{equation}
\label{equation: n odd  -- Xi-- m=3p+2 -- mu=1 }
         \xi_{ m }   =  \frac{ ( m + 1 )^2 }{6} 
\end{equation}
%
%
\end{itemize}
%
%
giving the capacity of the square core. 
In the same way, from
(\ref{equation:  n odd  --  Square Core Elongation })
we get the following two cases
%
%
\begin{itemize}
%
%
\item  $  \mu_{ p }  = 0      \Rightarrow    f_{ 12 } ( m - 3 ) = f_{ 12 } ( m ) - 4  =  m - 3    $     
%
%
and 
$ 
          (m - 3 )/2  \in \dot{ 2 }  
$
then  $ \xi_{ m - 3 } $ follows from
(\ref{equation: Xi_n -- n odd -- m = 3*p+2})
whence
%
%
\begin{equation}
\label{equation: n odd  -- DXi-- m=3p+2 -- mu=0 }
         \xi_{ m + 1 }  -  \xi_{ m - 3 }  =  \frac{ m ( m + 4 ) }{6}  -  \frac{ ( m - 2 )^2 }{6}  = \frac{ 4 \, m - 2 }{ 3 }   =  4 \, p  + 2
\end{equation}
%
%
%
\item  $  \mu_{ p }  = 1      \Rightarrow     f_{ 12 } ( m - 3 ) = f_{ 12 } ( m ) - 2 =  m - 2  $     
%
%
and 
$ 
          (m - 2 ) / 2  \in \dot{ 1 }  
$
then  $ \xi_{ m - 2 } $ follows from
(\ref{equation: Xi_n -- n odd -- m = 3*p+1})
whence
%
%
\begin{equation}
\label{equation: n odd  -- DXi-- m=3p+2 -- mu=1 }
         \xi_{ m }  -  \xi_{ m - 2 }  =  \frac{ ( m + 1 )^2 }{6}  -  \frac{  ( m - 3 ) ( m + 1 ) }{6}  = \frac{ 2 \, m + 2 }{3}   =  2 \, p + 2    
\end{equation}
%
%
\end{itemize}
%
%
give the expansion rate for the capacity in the square core.
%
%
%
\begin{figure}
\centering
\includegraphics[width=8cm]{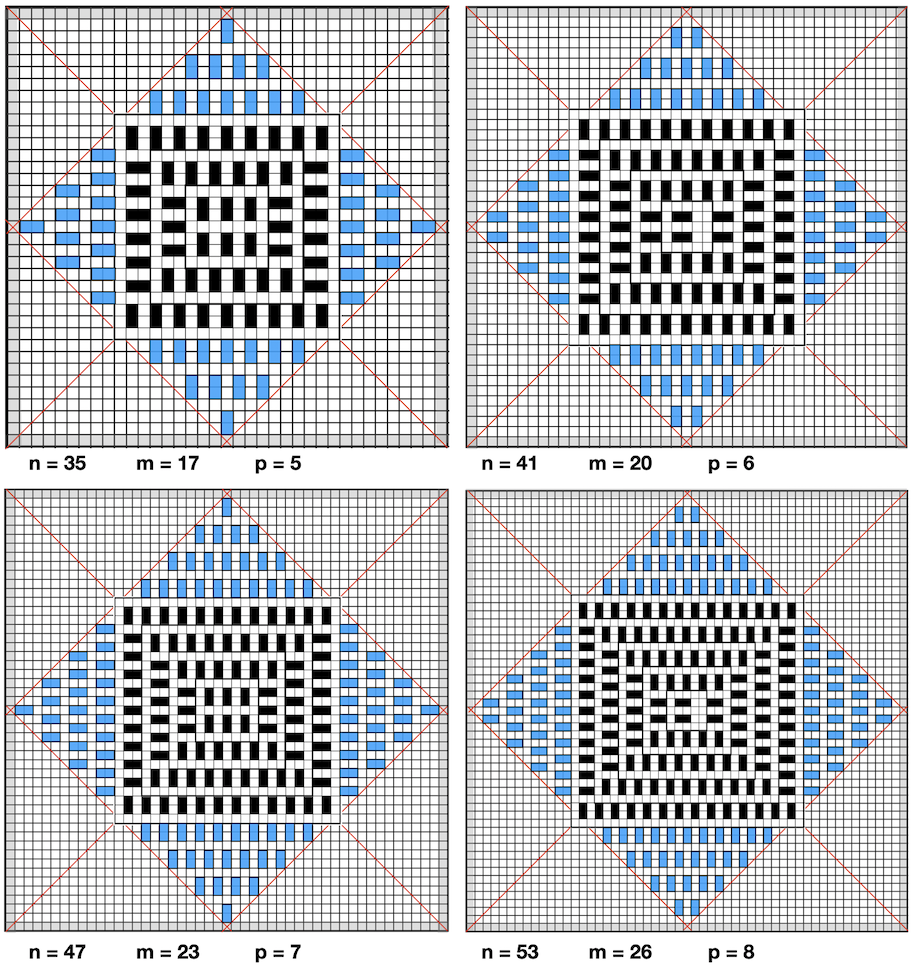}  
\caption{A 4--sequence of  $ \mathcal{D}_{ n } $ expansion in \texttt{Class} $ 12 $.
The class of the square core 
$
       \mathcal{S}_{ f_{12}( m )} 
$
alternates from $ ( \dot{ 12 } ) $ to $ ( \dot{ 11 } ) $ 
(refer to Fig.\,\ref{Figure: Dominos in the Square -- n odd -- n even-- })
whereas a new row is added in
$
         \mathcal{W}_{ 12 }
$
for $ p $ odd.
} 
\label{Figure: C12-N35N53 }
\end{figure}
%
%
%
%
\subsubsection{Capacity and expansion  of  $ \mathcal{D}_{ n } $ in \texttt{Class} $ 12 $ }                 
%
%
The capacity $ \psi_{ n } $ of  $ \mathcal{D}_{ n } $ in \texttt{Class} $ 12 $ is given by
%
%
\begin{equation}
\label{equation: Psi_n in Class 12 }
      \psi_{n}  =   \xi_{ f_{ 12 }( m ) }  + 4 \, W_{ 12 }( p ) 
\end{equation}
%
%
and with
$
          W_{ 12 } \equiv W_{ 10 }
$
we get
%
%
\begin{itemize}
%
%
\item  $  \mu_{ p }  = 0      \Rightarrow   $ from 
%
%
(\ref{equation: n odd  -- Xi-- m=3p+2 -- mu=0 })
and from
(\ref{equation: n odd  -- expression of W_10 -- m=3p })
\[        
             \psi_{n}  =  \frac{ m( m + 4 ) }{6}   +  \frac{  p  \, ( 3 \, p + 2 ) }{ 2 }  
\]
%
%
\item  $  \mu_{ p }  = 1      \Rightarrow    $ from   
%
%
(\ref{equation: n odd  -- Xi-- m=3p+2 -- mu=1 })
and from
(\ref{equation: n odd  -- expression of W_10 -- m=3p })
%
%
%
\[ 
         \psi_{n}  =   \frac{ ( m + 1 )^2 }{6}  +  \frac{   ( p + 1 )  \, ( 3 \, p + 1 ) }{ 2 }   
\]
%
%
\end{itemize}
%
%
and finally
%
%
\begin{equation}
\label{equation: n odd  -- PSi in C12 }
        \psi_{n}  =    ( 3 p + 2 ) ( p + 1 )  =   \frac{ ( n - 1 ) ( n + 1 ) }{ 12 }  
\end{equation}
%
%
whatever the parity of $ p $.

The expansion rate of  $ \psi_{ n } $ can also be obtained stepwise:
it follows from 
(\ref{equation: n odd  -- DXi-- m=3p+2 -- mu=0 }--\ref{equation: n odd  -- DXi-- m=3p+2 -- mu=1 })
and
from (\ref{equation: n odd  -- W_10 expansion -- m=3p })
that is, with $ p > 0 $  
\[
           \psi_{n} - \psi_{ n - 6 } =  ( \xi_{ f_{ 12 }( m )}  - \xi_{ f_{ 12 }( m - 3 )} )  +  4 \, ( W_{ 12 }( p ) - W_{ 12 }( p - 1 )  ) 
\]
\[
            =    \mu_{ p } ( 2 \, p + 2 ) \,  +  \, ( 1 -  \mu_{ p } ) ( 4 \, p + 2 )  \,  + 2 \, p \, ( 1 +  \mu_{ p } )  =  6 \, p  +  2       
\]
now 
$ 
           n = 2\,m + 1 
$ 
and 
$
          m = 3\,p + 2
$
then
%
%
\begin{equation}
\label{equation: n odd -- m == 2  -- Psi(n) - Psi( n - 6 ) }
          \psi_{ 5 }  = 2  \ ; \hspace{2mm}  \psi_{n} - \psi_{ n - 6 } =  n - 3  \, .
\end{equation}
%
%
%
The values of $ \psi_{n} $ in \texttt{Class} $ 12 $ are displayed in 
Table \ref{Table: Domino enumeration -- Class 12 }.
A 4--sequence of expansions is displayed in 
Fig.\,\ref{Figure: C12-N35N53 }.
%
%
%
%
%
%
%
%
%
%
\section{Dominos in the Diamond --- $n$ even}
\label{Section:Dominos in the Diamond -- n even}
\begin{figure}
\centering
\includegraphics[width=7cm]{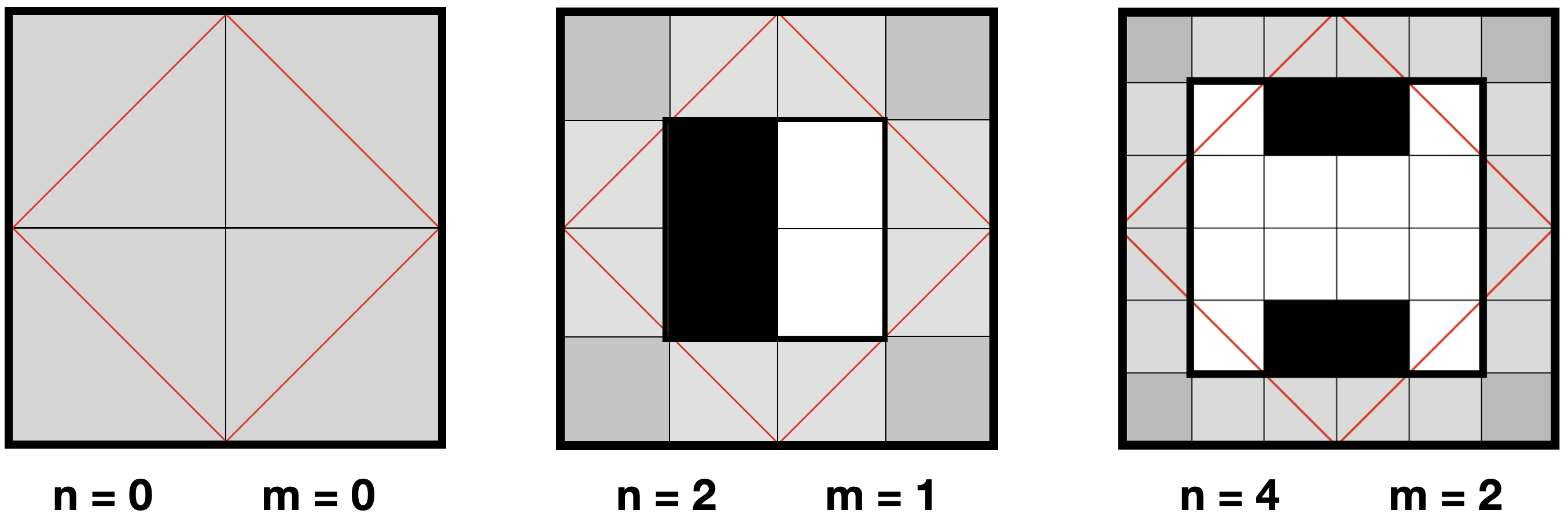}  
\caption{-- $n$ even. The first configurations with $ p = 0 $. From left to right:
Class $ { 00 } $ ($  \psi_{ 0 }  = 0 $),
Class $ { 01 } $ ($  \psi_{ 2 }  = 1 $),
Class $ { 02 } $ ($  \psi_{ 4 }  = 2 $). \newline
} 
\label{Figure: C0m--P0 } 
\end{figure}

%
%
Again and like for the {\em odd} case, the following constructions for 
$
        \overline{ 0\,m }
 $
 classes (or ``$ 0 \, m $'') will be valid only when the core--wedge structure actually appears from 
Fig.\,\ref{Figure: C0-n6-n28 -- n even }.
Therefore, the same conditions on $ p $ will be assumed everywhere.
The first configurations with $ p = 0 $ are shown in
Fig.\,\ref{Figure: C0m--P0 }.
%
%
%
%
%

Whereas the case --\,$ n $ odd\,--  argued for a ``vertical'' layout (dominos perpendicular to the adjacent core’s side), now the case --\,$ n $ even\,-- argues for a ``horizontal'' layout (dominos parallel to the adjacent core’s side).
Referring back to
Fig.\,\ref{Figure: Core-Wedges } ({\em right})
the exact size 
$
         f_{ 0 m } \approx  m 
$
of the core will be adjusted so that the configuration of the wedges is entirely fixed by the $p$ parameter.
Let
$
             f_{ 0 m }^{+} ( m ) =   f_{ 0 m } ( m )  + 2   
$
be the extended length of the side of 
$
       \mathcal{S}_{  f_{ 0 m }(m)} 
$
including borders, and let  $ h_{ 0 m } $ such that
\[
             f_{ 0 m }^{+} ( m )  = n -  2  \cdot  h_{ 0 m } ( p )  \hspace{10mm}  ( n = 2 m )
\]
where $ h_{ 0 m } $ is the height of the vertical median column of Region
$
       \mathcal{W}_{ 0 m } .
$
By observing that this wedge is an isosceles right triangle, we fix its height so that
%
%
\begin{equation}
\label{equation:  n even  --  2h_{ 0 m }( m )  }
               2  \cdot  h_{ 0 m } ( p ) =  m - 3  +  \left\{
                        \begin{array}{ll}
                                 r_{ p + 1 }		& 	 \mbox{ in  \texttt{Class} $ 00 $ }  \\
                                  r_{ p }		& 	 \mbox{ in  \texttt{Class} $ 01 $ }  \\
                                   r_{ p - 1 } 	& 	 \mbox{ in  \texttt{Class} $ 02 $ }  
                        \end{array}
                                                                                                \right.
\end{equation}
%
%
%
%
\begin{equation}
\label{equation:  n even  --  p = 4*q + r  -- C00 C01 C02 }
          \mbox{ with }   \hspace{2mm}  \left\{
                        \begin{array}{clll}
                                p + 1 	&	=   4 \, q_{ p + 1 } + r_{ p + 1 }	&  ( 0  \le  r_{ p + 1 } < 4)	&	 \mbox{ in  \texttt{Class} $ 00 $ }  \\
                                p     	&	=   4 \, q_{    p    } + r_{    p    }	&  ( 0  \le  r_{    p    } < 4)	&	 \mbox{ in  \texttt{Class} $ 01 $ }  \\
                                p - 1 	&      =   4 \, q_{ p - 1 } + r_{ p - 1 }	&  ( 0  \le  r_{ p - 1 } < 4)	&	 \mbox{ in  \texttt{Class} $ 02 $ }  \\
                        \end{array}
                                                                                                \right.
\end{equation}
%
%
whence finally the strict length of the side in the square core
%
%
\begin{equation}
\label{equation:  n even  --  f_{ 0 m }( m ) -- C00 C01 C02 }
            f_{ 0 m }( m ) =    f_{ 0 m }^{+} ( m ) - 2 = \left\{
                        \begin{array}{ll}
                                 m + 1 - r_{ p + 1 }	& 	 \mbox{ in  \texttt{Class} $ 00 $ }  \\
                                 m + 1 - r_{ p }		& 	 \mbox{ in  \texttt{Class} $ 01 $ }  \\
                                 m + 1 - r_{ p - 1 }	& 	 \mbox{ in  \texttt{Class} $ 02 $ }  \\
                        \end{array}
                                                                                                \right.
\end{equation}
%
%
expressed in terms of $ m $.
Or, by using a general unifying notation we get
%
%
\begin{equation}
\label{equation:  n even  --  f_{ 0 m }( m ) }
          f_{ 0 m }( m )   =   m + 1  -  r_{\widehat{ p }}	
\end{equation}
%
%
%
\begin{equation}
\label{equation:  n even  --  p = 4*q + r  }
          \widehat{ p }   =   4 \,   q_{\widehat{ p }}   +   r_{\widehat{ p }} \ :   \hspace{5mm}   \widehat{ p }  = p + 1 - m  \bmod 3  
           \hspace{5mm}    ( 0  \le  r_{\widehat{ p }}   <   4 )
\end{equation}
%
%
and where the 
$
               ( q_{\widehat{ p }}, r_{\widehat{ p }}  )  
$
are quotient and remainder of dividing by 4 the left--hand term in $ p $.
Moreover, with
$
            r_{\widehat{ p } - 1 }   =  ( r_{\widehat{ p }} - 1 )  \bmod  4 	    
$
we get
%
%
\begin{equation}
\label{equation:  n even  --  Square Core Elongation }
           f_{ 0 m }( m )  -  f_{ 0 m }( m - 3 ) =   3  -   r_{\widehat{ p }}  +  (  r_{\widehat{ p }} - 1 )  \bmod  4   =   \left\{
                        \begin{array}{ll}
                                 + 6	& 	 \mbox{ for }   r_{\widehat{ p }}  =  0	\\
                                 + 2	& 	 \mbox{ otherwise }  		 
                        \end{array}
                                                                                                \right.
\end{equation}
%
%
giving the elongation of the side of the core within a
$
          p - 1  \rightarrow  p  
$
expansion.
The evolution will follow a 4--cycle sequence, namely a ``$ (q_{\widehat{ p }}) $--cycle'' at constant 
$ 
          q_{\widehat{ p }}. 
$
%
%
\subsection{ \texttt{Class} $ 00 $ --- $ n $ even \, \&  \, $ m \equiv  0  \pmod 3 $}
\label{Subsection: n  even --- m = 3p -- }
%
%
%
The basic parameters of \texttt{Class} $ 00 $ verify: 
\[
		p \, \in  \, \mathbb{N} \,  ;  \hspace{5mm} m  \, = 3 \, p \, ;  \hspace{5mm} n \, =  2 \, m \, ;
\]
%
%
\subsubsection{Capacity and expansion rate of the square core $  \mathcal{S}_{  f_{ 00 }  } $ }                 
%
%
From
(\ref{equation:  n even  --  f_{ 0 m }( m ) -- C00 C01 C02 }--\ref{equation:  n even  --  p = 4*q + r  })
$
            f_{00} ( m ) = m  + 1 - r_{ p + 1 } 	
$
and we get the following four cases
%
%
\begin{itemize}
%
%
\item  $   r_{ p + 1 }  = 0      \Rightarrow    f_{00} ( m ) =   m + 1  $     
%
%
and
$ 
          ( m + 1 ) / 2  \in \dot{2}   \hspace{1mm} \mbox{in} \hspace{1mm}     \mathbb{Z} / 3 \, \mathbb{Z}   
$
then  $ \xi_{ m + 1 } $ follows from
(\ref{equation: Xi_n -- n even -- m = 3*p+2})
whence
%
%
\begin{equation}
\label{equation: n even  -- Xi-- m=3p -- r = 0 }
         \xi_{ m + 1 }    =  \frac{ ( m + 1 ) ( m + 3 ) }{ 6 }
\end{equation}
%
%
%
\item  $   r_{ p + 1 }  = 1      \Rightarrow    f_{00} ( m ) =   m  $     
%
%
and 
$ 
         m / 2  \in \dot{ 0 }   \hspace{1mm}   
$
then  $ \xi_{ m } $ follows from
(\ref{equation: Xi_n -- n even -- m = 3*p})
whence
%
%
\begin{equation}
\label{equation: n even  -- Xi-- m=3p -- r = 1 }
         \xi_{ m }   =  \frac{ m ( m + 2 ) }{6}    
\end{equation}
%
%
%
\item  $   r_{ p + 1 }  = 2      \Rightarrow    f_{00} ( m ) =   m - 1  $     
%
%
and 
$ 
         ( m - 1 ) / 2  \in \dot{ 1 }   \hspace{1mm}   
$
then  $ \xi_{ m - 1 } $ follows from
(\ref{equation: Xi_n -- n even -- m = 3*p+1})
whence
%
%
\begin{equation}
\label{equation: n even  -- Xi-- m=3p -- r = 2 }
         \xi_{ m - 1 } =  \frac{ ( m - 1 ) ( m + 1 ) - 2 }{6}     
\end{equation}
%
%
%
\item  $   r_{ p + 1 }  = 3      \Rightarrow    f_{00} ( m ) =   m - 2  $     
%
%
and 
$ 
         ( m - 2 ) / 2  \in \dot{ 2 }   \hspace{1mm}   
$
then  $ \xi_{ m - 2 } $ follows from
(\ref{equation: Xi_n -- n even -- m = 3*p+2})
whence
%
%
\begin{equation}
\label{equation: n even  -- Xi-- m=3p -- r = 3 }
         \xi_{ m - 2 }  =  \frac{ m ( m - 2 ) }{6}  
\end{equation}
%
%
\end{itemize}
%
%
giving the capacity of the square core. In the same way, from
(\ref{equation:  n even  --  Square Core Elongation })
we get the following four cases
%
%
\begin{itemize}
%
%
\item  $   r_{ p + 1 }  = 0      \Rightarrow    f_{00} ( m - 3 ) = f_{00} ( m ) - 6   $     
%
%
whence from
(\ref{equation: Induction Xi -- n even }) 
%
%
\begin{equation}
\label{equation: n even  -- DXi-- m=3p -- r = 0 }
         \xi_{ m + 1 }  -  \xi_{ ( m + 1 ) - 6 }  =  2 ( m - 1 )=  6 \, p - 2
\end{equation}
%
%
%
\item  $   r_{ p + 1 }  = 1      \Rightarrow    f_{00} ( m - 3 ) =   f_{00} ( m )  - 2  =  m - 2  $     
%
%
and 
$ 
          ( m - 2  ) / 2  \in \dot{ 2 }  
$
then  $ \xi_{ m - 2 } $ follows from
(\ref{equation: Xi_n -- n even -- m = 3*p+2})
whence
%
%
\begin{equation}
\label{equation: n even  -- DXi-- m=3p -- r = 1 }
         \xi_{ m }  -  \xi_{ m - 2 }  =  \frac{ m ( m + 2 ) }{6}  -   \frac{ m ( m - 2 ) }{6}  =  2 \, p    
\end{equation}
%
%
%
\item  $   r_{ p + 1 }  =  2      \Rightarrow    f_{00} ( m - 3 ) =   f_{00} ( m ) - 2   =  m - 3  $    
%
%
and 
$ 
          ( m - 3 ) / 2  \in \dot{ 0 }  
$
then  $ \xi_{ m - 3 } $ follows from
(\ref{equation: Xi_n -- n even -- m = 3*p})
whence
%
%
\begin{equation}
\label{equation: n even  -- DXi-- m=3p -- r = 2 }
         \xi_{ m - 1 }  -  \xi_{ m - 3 }  =  \frac{ ( m - 1 ) ( m + 1 ) - 2 }{6}  -   \frac{ ( m - 3 )( m - 1 ) }{6}  =  2 \, p  - 1  
\end{equation}
%
%
%
\item  $   r_{ p + 1 }  = 3      \Rightarrow    f_{00} ( m - 3 ) =   f_{00} ( m ) - 2   =   m - 4  $ 
%
%
and 
$ 
          ( m - 4 ) / 2  \in \dot{ 1 }  
$
then  $ \xi_{ m - 4 } $ follows from
(\ref{equation: Xi_n -- n even -- m = 3*p+1})
whence
%
%
\begin{equation}
\label{equation: n even  -- DXi-- m=3p -- r = 3 }
         \xi_{ m - 2 }  -  \xi_{ m - 4 }  =  \frac{ m ( m - 2 ) }{6}  -   \frac{ ( m - 4 )( m - 2 ) - 2 }{6}  =  2 \, p  - 1  
\end{equation}
%
%
\end{itemize}
%
%
giving the expansion rate for the capacity in the square core. 
%
%
\subsubsection{Capacity and expansion rate of the wedge $  \mathcal{W}_{ 00 } $ }                 
%
%
Let $ W_ { 00 }( p ) $ be the capacity of the wedge for any given $ p $. From
(\ref{equation:  n even  --  2h_{ 0 m }( m )  })
comes
\[
               2  \, (  h_{ 00 } ( p )  +  1  )  =  m - 1  +  \left\{
                        \begin{array}{ll}
                                 r_{ p + 1 }		& 	 \mbox{ in  \texttt{Class} $ 00 $ }  \\
                                  r_{ p }		& 	 \mbox{ in  \texttt{Class} $ 01 $ }  \\
                                   r_{ p - 1 } 	& 	 \mbox{ in  \texttt{Class} $ 02 $ }  
                        \end{array}
                                                                                                \right.
\]
and thus we get the capacity of the vertical median of Region
$
       \mathcal{W}_{ 00 } 
$
\[
              \overline{h}_{00} ( p ) =  (  h_{00} ( p ) + 1 ) / 2 =  \frac{1}{4} \, ( 3 \, p  - 1  +  r_{ p + 1 }  ) 
\]
in number of horizontal dominos --\,namely, in number of horizontal rows of dominos.
From
(\ref{equation:  n even  --  p = 4*q + r  -- C00 C01 C02 })
let us bring out the entity $ p - q_{ p + 1 }  $ as
%
%
\begin{equation}
\label{equation: n even  --  Capacity of vertical median of  W_{00} -- m = 3p -- }
               p - q_{ p + 1 }   = 3 \, q_{ p + 1 } + ( r_{ p + 1 }  - 1 )  \hspace{2mm} \mbox{or} \hspace{2mm} 
               4 \, ( p - q_{ p + 1 }  )  =  3 \, p + ( r_{ p + 1 }  - 1 )  
\end{equation}
%
%
by eliminating the $ q_{ p + 1 } $ term from the second member. 
It follows that
%
%
\begin{equation}
\label{equation: n even  --  median capacity h_(00) -- m = 3p -- }
	  \overline{h}_{00} ( p )  =    p - q_{ p + 1 }    \hspace{3mm}     \mbox{whence}   \hspace{3mm} 
           \overline{h}_{00} ( p )  -   \overline{h}_{00} ( p - 1 )  =   \left\{
                        \begin{array}{ll}
                                  0	& 	\mbox{ if }		 r_{ p + 1 }  =  0 	\\
                                 + 1 	& 	\mbox{ otherwise }		 
                        \end{array}
                                                                                                \right.
\end{equation}
%
%
and as a result, the number of rows increases within each internal step and remains constant at the start of a new cycle.
The baserow of the  \texttt{N}--wedge has the length
\[
              b_{00} ( p ) =  2 \cdot h_{00} ( p )  =    3 \, p  +  r_{ p + 1 }  - 3  
\]
rewritten as the sum
$
              b_{00} ( p ) = b'_{00} ( p )  +  b''_{00} ( p ) 
$
with
%
%
\begin{equation}
\label{equation: n even  --  baserow length b_(00) -- m = 3p -- }
        b'_{00} ( p ) =   3 \, ( p - q_{ p + 1 } )  \hspace{2mm} \mbox{and} \hspace{2mm}  
        b''_{00} ( p )  =  3 \, ( q_{ p + 1 } - 1 ) + r_{ p + 1 }  \hspace{5mm}  ( p > 1 )
\end{equation}
%
%
that yields
$
          \overline{b}_{00} ( p ) =   \overline{b'}_{00} ( p )  +  \overline{b''}_{00} ( p ) 
$
where 
$ 
            \overline{b}_{00} 
$ 
is the capacity in terms of (horizontal) dominos of the baserow of Region
$
         \mathcal{W}_{00}
$
and where
$ 
          \overline{b'}_{00} \hspace{1mm}  \mbox{and}  \hspace{1mm} \overline{b''}_{00}
$ 
denote the capacity of the baserow of
$
         \mathcal{W}'_{00}  \hspace{1mm}  \mbox{and}  \hspace{1mm} \mathcal{W}''_{00}
$
respectively.
%
%
\paragraph{Capacity and expansion of the subwedge $  \mathcal{W}'_{ 00 } $.}                
%
%
Regarding 
$
         \mathcal{W}'_{00}
$
it follows that
%
%
\begin{equation}
\label{equation: n even  -- b'_{00} ==  h'_{00}  -- m = 3p -- }
               \overline{b'}_{00} ( p )  =   \frac{1}{3} \, \cdot  \,  b'_{00} ( p )  =   p - q_{ p + 1 } 
               =   \overline{h}_{00} ( p ) 
\end{equation}
%
%
and 
$
           \overline{b'}_{00} 
$
behaves like
$
           \overline{h'}_{00} 
$
in
(\ref{equation: n even  --  median capacity h_(00) -- m = 3p -- }).
 Now
 $
          W'_{00} (0) = 0
 $
 and from
(\ref{equation: n even  --  median capacity h_(00) -- m = 3p -- })
and
(\ref{equation: n even  -- b'_{00} ==  h'_{00}  -- m = 3p -- })
%
%
%
\begin{equation}
\label{equation: n even  -- W'_{00} ( p ) - W'_{00} ( p - 1 ) -- m = 3p -- }
          W'_{00} ( p )  -   W'_{00} ( p - 1 )  =  \left\{
                        \begin{array}{ll}
                                   0				& 	\mbox{ if }		 r_{ p + 1 }  =  0 	\\
                                   p  -  q_{ p + 1 } 	& 	\mbox{ otherwise }		 
                        \end{array}
                                                                                                \right.
\end{equation}
%
%
and with
%
%
\begin{equation}
\label{equation: n even  -- W'_{00} ( p )  -- m = 3p -- }
            W'_{00} ( p )  =  \sum_{ k = 1 }^{ p - q_{ p + 1 } } k  =  ( p - q_{ p + 1 } ) \, ( p - q_{ p + 1 } + 1 ) / 2
\end{equation}
%
%
we finally obtain the total capacity of Region
$
         \mathcal{W}'_{00}.
$
Or expressed in another way
%
%
\begin{equation}
\label{equation: n even  -- W'_{0m} ( p ) - W'_{0m} ( p - 1 ) -- }
          W'_{0m} ( p )  -   W'_{0m} ( p - 1 )  =  \left\{
                        \begin{array}{ll}
                                   0				& 	\mbox{ if }		 r_{\widehat{ p }}  =  0 	\\
                                      \overline{b'}_{0m} ( p )  =  p  -  q_{ \widehat{ p } } 	& 	\mbox{ otherwise }		 
                        \end{array}
                                                                                                \right.
\end{equation}
%
%
%
%
\begin{equation}
\label{equation: n even  -- W'_{0m} ( p )  -- }
            W'_{0m} ( p )  =  ( p  -  q_{ \widehat{ p } } ) \, (  p  -  q_{ \widehat{ p } } + 1 ) / 2      
\end{equation}
%
%
%
%
\begin{equation}
\label{equation: n even  -- W'_{0m} ( p )  -- }
             4 \cdot  W'_{0m} ( p )  =  2 \, ( p  -  q_{ \widehat{ p } } ) \, (  p  -  q_{ \widehat{ p } } + 1 )     
\end{equation}
%
%
now without $ m $--class--aware statement.
In \texttt{Class} $ 00 $, from 
(\ref{equation: n even  --  Capacity of vertical median of  W_{00} -- m = 3p -- })
and
(\ref{equation: n even  -- W'_{00} ( p )  -- m = 3p -- })
it comes now
\[
          4 \cdot  W'_{00} ( p )  =   \frac{   ( 3 \, p + r_{ p + 1 }  - 1 ) ( 3 \, p +  r_{ p + 1 }   +  3 ) }{ 8 }   
\]
whence
%
%
\begin{equation}
\label{equation: n even  -- 4 * W'_{00} ( p )  -- }
          4 \cdot  W'_{00} ( p )  =  \left\{
                        \begin{array}{ll}
                                   ( 3 \, p  - 1 ) ( 3 \, p  +  3 ) / 8 	& 	\mbox{ if }		 r_{ p + 1 }  = 0  	\\
                                   ( 3 \, p  ) ( 3 \, p  +  4 ) / 8 	& 	\mbox{ if }		 r_{ p + 1 }  = 1  	\\
                                   ( 3 \, p  + 1 ) ( 3 \, p  +  5 ) / 8 	& 	\mbox{ if }		 r_{ p + 1 }  = 2  	\\
                                  ( 3 \, p  + 2 ) ( 3 \, p  +  6 ) / 8 	& 	\mbox{ if }		 r_{ p + 1 }  = 3  	\\
                        \end{array}
                                                                                                \right.
\end{equation}
%
%
%
giving the overall capacity of the four  \texttt{N}--\texttt{S}--\texttt{E}--\texttt{W} subregions $ \mathcal{W}'_{00} $ in terms of $ p $.
%
%
\paragraph{Capacity and expansion of the subwedge $  \mathcal{W}''_{ 00 } $.}                
%
%
From
(\ref{equation: n even  --  baserow length b_(00) -- m = 3p -- })
we get
%
%
\begin{equation}
\label{equation: n even  --  capacity of b"_(00) -- m = 3p -- }
              \overline{b''}_{00} ( p )   =   ( q_{ p + 1 } - 1 )  +   \delta_{ r_{ p + 1 } }    \hspace{3mm} \mbox{with} \hspace{3mm} 
              \delta_{ r_{ p + 1 } }  =   \Bigl\lfloor  \,\frac { r_{ p + 1 }  \, }  { 3 }  \Bigr\rceil 
\end{equation}
%
%
where the second term on the right side denotes the integer rounding closest to the rational $  r_{ p + 1 } / 3 $.
Region 
$
         \mathcal{W}''_{00}
$
is empty for $ p = 1 $ then for any $ p > 1 $
%
%
\begin{equation}
\label{equation: n even  -- W''_{00} ( p ) - W''_{00} ( p - 1 ) -- m = 3p -- }
         	W''_{00} ( p )  -   W''_{00} ( p - 1 )  =  \left\{
                        \begin{array}{ll}
                               	0							 	& 	\mbox{ if }	    r_{ p + 1 } = 0 	\\
                                  q_{ p + 1 } - 1  +   \delta_{ r_{ p + 1 } }  	& 	\mbox{ otherwise. }		 
                         \end{array}
                                                                           \right.
\end{equation}
%
%
The periodic sequence of expansion within a $ q_{ p + 1 } $--cycle (at constant $ q_{ p + 1 } $) is as follows:
%
%
\begin{itemize}
%
%
\item  $  r_{ p + 1 } = 0 $, $ \delta_{ r_{ p + 1 } } = 0  \Rightarrow  W''_{00}( p )  =   W''_{00}( p - 1 ) $
and the baserow can accommodate exactly $ q_{ p + 1 } - 1 $ dominos.  
%
%
%
\item  $  r_{ p + 1 } = 1 $, $ \delta_{ r_{ p + 1 } } = 0  \Rightarrow  W''_{00}( p )  =   W''_{00}( p - 1 )  +  q_{ p + 1 } - 1  $
and a new lacunar void is inserted in the new baserow.
%
%
%
\item  $  r_{ p + 1 } = 2 $, $ \delta_{ r_{ p + 1 } } = 1  \Rightarrow  W''_{00}( p )  =   W''_{00}( p - 1 )  +  q_{ p + 1 }  $
and a new domino is inserted in the new baserow.
%
%
%
\item  $  r_{ p + 1 } = 3 $, $ \delta_{ r_{ p + 1 } } = 1  \Rightarrow  W''_{00}( p )  =   W''_{00}( p - 1 )  +  q_{ p + 1 }  $
and a new baserow of the same capacity is added.
%
%
\end{itemize}
%
%
By adding up all items, the expansion in
$
         \mathcal{W}''_{00}
$
at the end of the $ q_{ p + 1 } $--cycle becomes 
$
         3 \, q_{ p + 1 } - 1.
$
Now, at the end of the first cycle $ W''_{00}( 2 )  = 0 $ and with $  r_{ p + 1 } = 3 $ we can get
%
%
\begin{equation}
\label{equation: n even  -- W''_{00} ( p )  -- m = 3p -- Sigma -- }
            W''_{ 00 }( p )  =   \sum_{ k = 1 }^{ q_{ p + 1 } } ( 3 \, k - 1 )      
            =    \frac{ q_{ p + 1 }  \, ( 3 \, q_{ p + 1 }  + 1 ) }{ 2 } 
\end{equation}
%
%
and subtracting step by step back
%
%
\begin{equation}
\label{equation: n even  -- W''_{00} ( p )  -- m = 3p -- Stepback -- }
           W''_{00} ( p )   =  \left\{
                        \begin{array}{ll}
                               	   q_{ p + 1 }  \, ( 3 \, q_{ p + 1 }  + 1 ) / 2		& 	\mbox{ if }	 r_{ p + 1 }  =  3 		\\
                               	  q_{ p + 1 } \, (  3 \, q_{ p + 1 }  -  1  ) / 2		& 	\mbox{ if }	 r_{ p + 1 }  =  2 		\\
                               	  (  q_{ p + 1 } - 1 ) \, 3 \, q_{ p + 1 } / 2 		& 	\mbox{ if }	 r_{ p + 1 }  =  1 		\\
                               	  (  q_{ p + 1 } - 1 ) \, (  3 \, q_{ p + 1 } - 2 ) 	/ 2 	& 	\mbox{ if }	 r_{ p + 1 }  =  0 		
                        \end{array}
                                                                           \right.
\end{equation}
%
%
or even in the form
%
%
\begin{equation}
\label{equation: n even  -- W''_{00} ( p )  -- m = 3p -- }
                W''_{00}( p )  =  \overline{b''}_{00} ( p )  \cdot
                \frac{  ( p - q_{ p + 1 }  -  3 )  +  ( p - q_{ p + 1 }  -  3  \cdot  \overline{b''}_{00} ( p ) ) }  { 2 } 
\end{equation}
%
%
in order to get rid of the dependence of $ W''_{00} $ on $  r_{ p + 1 } $ and also to emphasize that  $ W''_{00} $ grows with the area of Triangle
$
         \mathcal{W}''_{00}.
$
Or in another way, from
(\ref{equation: n even  --  capacity of b"_(00) -- m = 3p -- })
comes
%
%
\begin{equation}
\label{equation: n even  --  capacity of b"_(0m) -- }
              \overline{b''}_{0m} ( p )   =   (   q_{ \widehat{ p } }  - 1 )  +    \delta_{  r_{ \widehat{ p } } }    
              \hspace{3mm} \mbox{with} \hspace{3mm} 
               \delta_{  r_{ \widehat{ p } } } =   \Bigl\lfloor  \,\frac {r_{ \widehat{ p } }  + m  \bmod 3 \, }  { 3 }  \Bigr\rceil 
\end{equation}
%
%
then from
(\ref{equation: n even  -- W''_{00} ( p ) - W''_{00} ( p - 1 ) -- m = 3p -- })
and
(\ref{equation: n even  -- W''_{00} ( p )  -- m = 3p  -- })
we get
%
%
\begin{equation}
\label{equation: n even  -- W''_{0m} ( p ) - W''_{0m} ( p - 1 ) -- }
         	W''_{0m} ( p )  -   W''_{0m} ( p - 1 )  =  \left\{
                        \begin{array}{ll}
                    0							 & 	\mbox{ if }	     r_{ \widehat{ p } }  = 0 	\\
                    \overline{b''}_{0m} ( p )  =  q_{ \widehat{ p } }  - 1  +   \delta_{  r_{ \widehat{ p } } }  	& 	\mbox{ otherwise }		 
                         \end{array}
                                                                           \right.
\end{equation}
%
%
%
%
\begin{equation}
\label{equation: n even  -- W''_{0m} ( p )  -- }
                W''_{0m}( p )  =  \overline{b''}_{0m} ( p )  \cdot
                \frac{  ( p - q_{ \widehat{ p } }  -  3 )  +  ( p - q_{ \widehat{ p } } -  3  \cdot  \overline{b''}_{0m} ( p ) ) }  { 2 } 
\end{equation}
%
%
%
%
\begin{equation}
\label{equation: n even  -- W''_{0m} ( p )  -- }
               4  \cdot W''_{0m}( p )  =  2 \cdot  \overline{b''}_{0m} ( p )  \cdot
             ( ( p - q_{ \widehat{ p } }  -  3 )  +  ( p - q_{ \widehat{ p } } -  3  \cdot  \overline{b''}_{0m} ( p ) ) )  
\end{equation}
%
%
now expressed according to a $ m $--class--unaware statement.
In \texttt{Class} $ 00 $, substituting   
$
             q_{ p + 1 }  =  ( p + 1 - r_{ p + 1 } ) / 4
$
from
(\ref{equation:  n even  --  p = 4*q + r  -- C00 C01 C02 })
in
(\ref{equation: n even  -- W''_{00} ( p )  -- m = 3p -- Stepback -- })
it comes now
%
%
\begin{equation}
\label{equation: n even  -- 4 * W''_{00} ( p )  -- }
          4 \cdot  W''_{00} ( p )  =  \left\{
                        \begin{array}{ll}
                                   ( p - 3 ) ( 3 p - 5 ) / 8		& 	\mbox{ if }		 r_{ p + 1 }  = 0  	\\
                                    ( p - 4 ) ( 3 p ) / 8		& 	\mbox{ if }		 r_{ p + 1 }  = 1  	\\
                                    ( p - 1 ) ( 3 p - 7 ) / 8 	& 	\mbox{ if }		 r_{ p + 1 }  = 2  	\\
                                    ( p - 2 ) ( 3 p - 2 ) / 8	& 	\mbox{ if }		 r_{ p + 1 }  = 3  	\\
                        \end{array}
                                                                                                \right.
\end{equation}
%
%
giving the overall capacity of the four  \texttt{N}--\texttt{S}--\texttt{E}--\texttt{W} subregions $ \mathcal{W}''_{00} $ in terms of $ p $.
%
%
\subsubsection{Capacity and expansion  of  $ \mathcal{D}_{ n } $ in \texttt{Class} $ 00 $ }                 
%
%
The capacity $ \psi_{ n } $ of  $ \mathcal{D}_{ n } $ in \texttt{Class} $ 00 $ is simply given by
%
%
\begin{equation}
\label{equation: Psi_n in Class 00 }
      \psi_{n}  =   \xi_{ f_{00}( m ) }  + 4 \, ( W'_{00}( p )  +  W''_{00}( p ) )
\end{equation}
%
%
then from
(\ref{equation: n even  -- Xi-- m=3p -- r = 0 })--(\ref{equation: n even  -- Xi-- m=3p -- r = 3 })  
for
$
           \xi_{ f_{00} ( m ) }
$
and from
(\ref{equation: n even  -- 4 * W'_{00} ( p )  -- })
and 
(\ref{equation: n even  -- 4 * W''_{00} ( p )  -- })
for
$
         W'_{00} ( p )
$
and
$
         W''_{00} ( p )
$
it would easily come
%
%
\begin{equation}
\label{equation: n even  -- Psi_n ( n )  in  C00 }
            \psi_{n}   =    \left\{
                        \begin{array}{ll}
                                   ( n^2 + 2 n + 24 ) / 12	& 	\mbox{ if }		 r_{ p + 1 }  = 0  	\\
                                   ( n^2 + 2 n        ) / 12	& 	\mbox{ if }		 r_{ p + 1 }  = 1  	\\
                                    ( n^2 + 2 n + 12 ) / 12 	& 	\mbox{ if }		 r_{ p + 1 }  = 2  	\\
                                    ( n^2 + 2 n + 24 ) / 12	& 	\mbox{ if }		 r_{ p + 1 }  = 3  	\\
                        \end{array}
                                                                                                \right.
\end{equation}
%
%
with $ m = 3 p $ and $ n = 2 m $.
With
$
             \psi_{ 0 }  = 0
$
the expansion rate of  $ \psi_{ n } $ can also be obtained stepwise. In
Tab.\,\ref{Table: Expansion Psi -- m = 3p },
the sum of the three terms resulting from relations
(\ref{equation: n even  -- DXi-- m=3p -- r = 0 }--\ref{equation: n even  -- DXi-- m=3p -- r = 3 })
for 
$
\Delta \xi_{ f_{00} ( m ) }, 
$
from
(\ref{equation: n even  -- W'_{00} ( p ) - W'_{00} ( p - 1 ) -- m = 3p -- })
for 
$
        \Delta  W'_{00} ( p )
$
and from
(\ref{equation: n even  -- W''_{00} ( p ) - W''_{00} ( p - 1 ) -- m = 3p -- })
for 
$
        \Delta  W''_{00} ( p ) 
$
gives the result.
Moreover, by adding backwards the four terms in the last column of 
Tab.\,\ref{Table: Expansion Psi -- m = 3p },
it follows that 
\[
              \psi_n -  \psi_{ n - 24 }  =  ( n - 2 )  + ( ( n - 6 ) -1 )  + ( ( n - 12 ) - 1 )  +  ( ( n - 18 ) - 4 ) 
\]
$
             =  4\, n - ( \, 2 + 7 + 13 + 22 \, )
$
whence the induction
%
%
\begin{equation}
\label{equation: n even -- m == 0  -- Psi(n) - Psi( n - 24 ) }
            \psi_{ 0 }  =   0 \, ;   \hspace{2mm}  \psi_n =  \psi_{ n - 24 }\ +   4\,( n - 11 )     \hspace{5mm}  ( p  \ge  4 )
\end{equation}
%
%
within any  
$
          p  \rightarrow  p + 4
$
cycle.
%
%
\begin{figure}
\centering
\includegraphics[width=12cm]{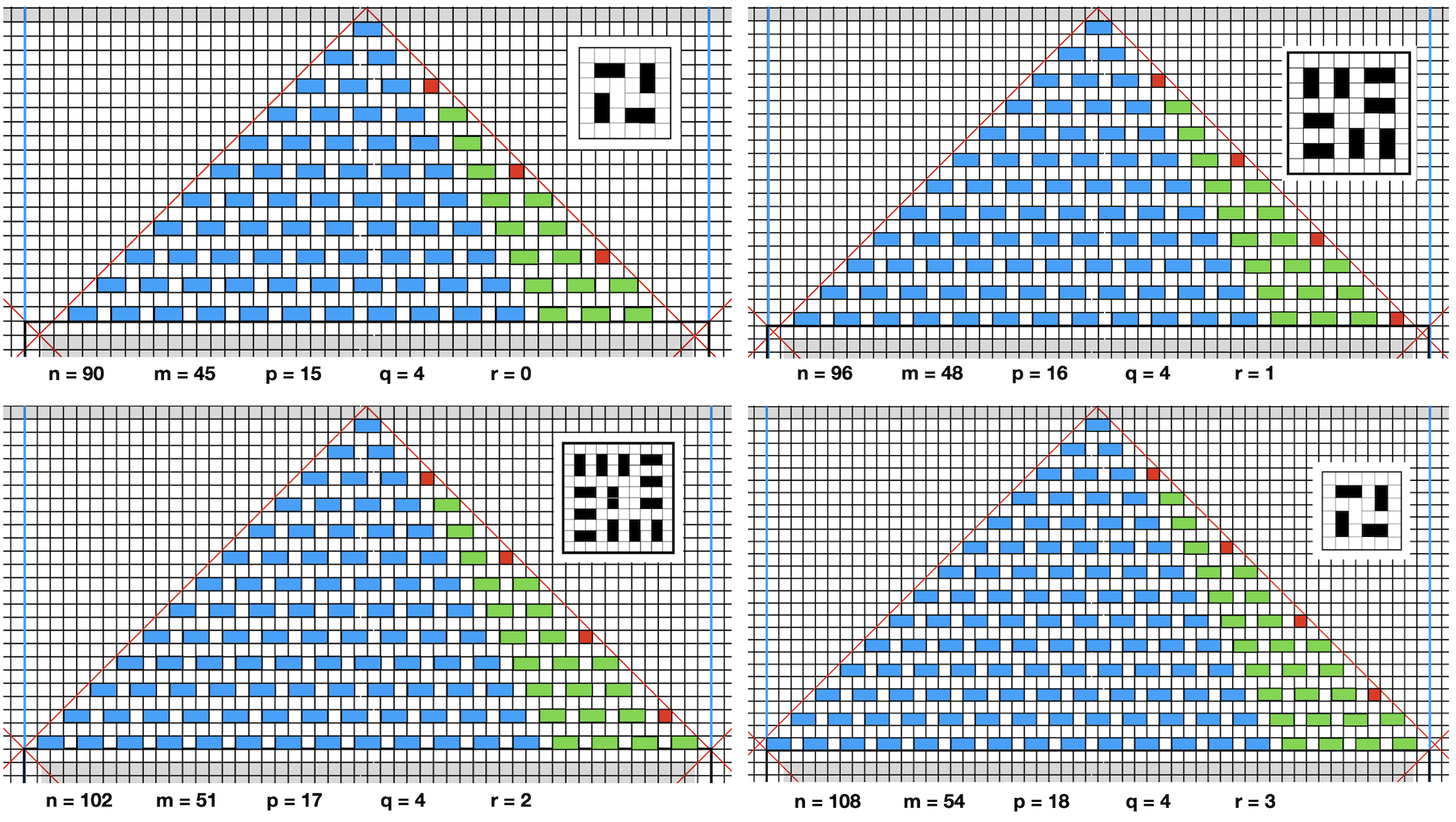}  
\caption{A 4--sequence of  $ \mathcal{D}_{ n } $ expansion in \texttt{Class} $ 00 $.
The class of the square core 
$
       \mathcal{S}_{ f_{00 }( m )} 
$
(refer to Fig.\,\ref{Figure: Dominos in the Square -- n odd -- n even-- })
evolves  as 
$ 
           ( \dot{ 02 } ) \rightarrow ( \dot{ 00 } ) \rightarrow ( \dot{ 01 })  \rightarrow ( \dot{ 02 } )
$ 
(highlighted in insets). A new sequence starts with $  p + 1 \equiv_4  0 $ whereas  
$
         \mathcal{W}_{00}
$
remains unchanged.
} 
\label{Figure: C00-P15P18 }
\end{figure}
%
%
%
\begin{table}[tb!]
\centering
\caption{Expansion in \texttt{Class} $ 00 $:  
$
           p \, \in  \, \mathbb{N} \,  ;  m  \, = 3 \, p    \, ;   n \, =  2 \, m  \, ;  \ - \  \hspace{1mm}  \psi_{ 0 }  = 0 \, ; \newline
           \mbox{ }  \hspace{13mm} 
            r_{ p + 1 } \equiv  p + 1  \pmod 4
$
}
\vspace{2mm} 
\label{Table: Expansion Psi -- m = 3p }
\begin{tabular}{|c||c|c|c||c|}
\hline
\rule{0pt}{10pt}
        $   r_{ p + 1 }    $	&    
        		$ \Delta \xi_{ f_{00} ( m ) } $	&   
							$ 4\, \Delta  W'_{00} ( p ) $	&   
														$ 4\, \Delta  W''_{00} ( p ) $	& 
																			$	\psi_{n}-\psi_{n-6}	$  \\  [0.5ex] 
\hline\hline	
	 0 	&	$   6\, p - 2 	$	&    $			0			$ 	&    $			0		$	&    	$	n - 2		$	\\   
 	 1 	&	$   2 \, p 		$	&    $  4 \, ( p - q_{ p + 1 } )	$  	&    $ 4 \, ( q_{ p + 1 } - 1 )	$ 	&     	$	n - 4		$	\\   	 	 
	 2  	&	$   2 \, p - 1	$	&    $  4 \, ( p - q_{ p + 1 } )	$  	&    $ 	4 \, q_{ p + 1 }	$	&	$	n - 1		$	\\
	 3  	&	$   2 \, p - 1	$	&    $  4 \, ( p - q_{ p + 1 } )	$  	&    $ 	4 \, q_{ p + 1 }	$ 	&	$	n - 1		$	\\  \hline	   	 
\end{tabular}
\end{table}
%
%
A 4--sequence of expansion is illustrated in  
Fig.\,\ref{Figure: C00-P15P18 }.
The values of $ \psi_{n} $ in \texttt{Class} $ 00 $ are displayed in 
Table \ref{Table: Domino enumeration -- Class 00 }.
\newpage 
%
%
%
\subsection{ \texttt{Class} $ 01 $ --- $ n $ even \, \&  \, $ m \equiv  1  \pmod 3 $}
\label{Subsection: n  even --- m = 3p+1 -- }
%
%
%
The basic parameters of \texttt{Class} $ 01 $ verify: 
\[
		p \, \in  \, \mathbb{N} \,  ;  \hspace{5mm} m  \, = 3 \, p \, + \, 1 \, ;  \hspace{5mm} n \, =  2 \, m \, ;
\]
%
%
\subsubsection{Capacity and expansion rate of the square core $  \mathcal{S}_{  f_{ 01 }  } $ }                 
%
%
From
(\ref{equation:  n even  --  f_{ 0 m }( m ) -- C00 C01 C02 }--\ref{equation:  n even  --  p = 4*q + r  })
$
            f_{01} ( m ) = m  + 1 - r_{ p } 	
$
and we get the following four cases
%
%
%
\begin{itemize}
%
%
\item  $   r_{ p }  = 0      \Rightarrow    f_{01} ( m ) =   m + 1  $     
%
%
and
$ 
          ( m + 1 ) / 2  \in \dot{ 1 }   \hspace{1mm} \mbox{in} \hspace{1mm}    \mathbb{Z} / 3 \, \mathbb{Z}
$
then  $ \xi_{ m + 1 } $ follows from
(\ref{equation: Xi_n -- n even -- m = 3*p+1})
whence
%
%
\begin{equation}
\label{equation: n even  -- Xi-- m=3p+1 -- r = 0 }
         \xi_{ m + 1 }    =  \frac{ ( m + 1 ) ( m + 3 )  - 2 }{ 6 }
\end{equation}
%
%
%
%
\item  $   r_{ p }  = 1      \Rightarrow    f_{01} ( m ) =   m  $     
%
%
and 
$ 
         m / 2  \in \dot{ 2 }   \hspace{1mm}   
$
then  $ \xi_{ m } $ follows from
(\ref{equation: Xi_n -- n even -- m = 3*p+2})
whence
%
%
\begin{equation}
\label{equation: n even  -- Xi-- m=3p+1 -- r = 1 }
         \xi_{ m }  =  \frac{ m ( m + 2 ) }{6}  
\end{equation}
%
%
%
%
\item  $   r_{ p }  = 2      \Rightarrow    f_{01} ( m ) =   m - 1  $     
%
%
and 
$ 
         ( m - 1 ) / 2  \in \dot{ 0 }   \hspace{1mm}   
$
then  $ \xi_{ m - 1 } $ follows from
(\ref{equation: Xi_n -- n even -- m = 3*p})
whence
%
%
\begin{equation}
\label{equation: n even  -- Xi-- m=3p+1 -- r = 2 }
         \xi_{ m - 1 }  =  \frac{ ( m - 1 ) ( m + 1 ) }{6}     
\end{equation}
%
%
%
%
\item  $   r_{ p }  = 3      \Rightarrow    f_{01} ( m ) =   m - 2  $     
%
%
and 
$ 
         ( m - 2 ) / 2  \in \dot{ 1 }   \hspace{1mm}   
$
then  $ \xi_{ m - 2 } $ follows from
(\ref{equation: Xi_n -- n even -- m = 3*p+1})
whence
%
%
\begin{equation}
\label{equation: n even  -- Xi-- m=3p+1 -- r = 3 }
         \xi_{ m - 2 }   =  \frac{ m ( m - 2 ) - 2 }{6}    
\end{equation}
%
%
%
\end{itemize}
%
%
giving the capacity of the square core. In the same way, from
(\ref{equation:  n even  --  Square Core Elongation })
we get the following four cases
%
%
\begin{itemize}
%
%
\item  $   r_{ p }  = 0      \Rightarrow    f_{01} ( m - 3 ) =    f_{01} ( m ) -  6  $    
%
%
whence from
(\ref{equation: Induction Xi -- n even }) 
%
%
\begin{equation}
\label{equation: n even  -- DXi-- m=3p+1 -- r = 0 }
         \xi_{ m + 1 }  -  \xi_{ ( m + 1 ) - 6 }  =  2 ( m - 1 )=  6 \, p 
\end{equation}
%
%
%
\item  $   r_{ p }  = 1      \Rightarrow    f_{01} ( m - 3 ) =    f_{01} ( m ) -  2  =  m - 2  $     
%
%
and 
$ 
          ( m - 2  ) / 2  \in \dot{ 1 }  
$
then  $ \xi_{ m - 2 } $ follows from
(\ref{equation: Xi_n -- n even -- m = 3*p+1})
whence
%
%
\begin{equation}
\label{equation: n even  -- DXi-- m=3p+1 -- r = 1 }
         \xi_{ m }  -  \xi_{ m - 2 }  =  \frac{ m ( m + 2 ) }{6}  -   \frac{ m ( m - 2 ) - 2 }{6}  =  2 \, p  + 1
\end{equation}
%
%
%
%
\item  $   r_{ p }  = 2      \Rightarrow    f_{01} ( m - 3 ) =    f_{01} ( m ) -  2  =  m - 3  $     
%
%
and 
$ 
          ( m - 3 ) / 2  \in \dot{ 2 }  
$
then  $ \xi_{ m - 3 } $ follows from
(\ref{equation: Xi_n -- n even -- m = 3*p+2})
whence
%
%
\begin{equation}
\label{equation: n even  -- DXi-- m=3p+1 -- r = 2 }
         \xi_{ m - 1 }  -  \xi_{ m - 3 }  =  \frac{ ( m - 1 ) ( m + 1 ) }{6}  -   \frac{ ( m - 3 )( m - 1 ) }{6}  =  2 \, p    
\end{equation}
%
%
%
%
\item  $   r_{ p }  = 3      \Rightarrow    f_{01} ( m - 3 ) =    f_{01} ( m ) -  2  =  m - 4  $     
%
%
and 
$ 
          ( m - 4 ) / 2  \in \dot{ 0 }  
$
then  $ \xi_{ m - 4 } $ follows from
(\ref{equation: Xi_n -- n even -- m = 3*p})
whence
%
%
\begin{equation}
\label{equation: n even  -- DXi-- m=3p+1 -- r = 3 }
         \xi_{ m - 2 }  -  \xi_{ m - 4 }  =  \frac{ m ( m - 2 ) - 2 }{6}  -   \frac{ ( m - 4 )( m - 2 ) }{6}  =  2 \, p  - 1  
\end{equation}
%
%
%
\end{itemize}
%
%
giving the expansion rate for the capacity in the square core.
%
%
\subsubsection{Capacity and expansion  of  $ \mathcal{D}_{ n } $ in \texttt{Class} $ 01 $ }                 
%
%
The capacity $ \psi_{ n } $ of  $ \mathcal{D}_{ n } $ in \texttt{Class} $ 01 $ is simply given by
%
%
\begin{equation}
\label{equation: Psi_n in Class 01 }
      \psi_{n}  =   \xi_{ f_{01}( m ) }  + 4 \, ( W'_{01}( p )  +  W''_{01}( p ) )
\end{equation}
%
%
%
and
$
           \xi_{ f_{01} ( m ) }
$
results from
(\ref{equation: n even  -- Xi-- m=3p+1 -- r = 0 })--(\ref{equation: n even  -- Xi-- m=3p+1 -- r = 3 }).
In \texttt{Class} $ 01 $ from
(\ref{equation:  n even  --  p = 4*q + r  -- C00 C01 C02 })
it comes
%
%
\begin{equation}
\label{equation: n even  --  Capacity of vertical median of  W_{01} -- m = 3p+1 -- }
            p - q_{ p }  =  3 \, q_{ p } + r_{ p }   \hspace{2mm} \mbox{or} \hspace{2mm} 
            4 \, ( p - q_{ p }  )  =  3 \, p +  r_{ p }  
\end{equation}
%
%
and from
(\ref{equation: n even  -- W'_{0m} ( p )  -- })
%
%
%
\begin{equation}
\label{equation: n even  -- W'_{01} ( p )  -- m = 3p+1 -- }
            W'_{01} ( p )  =  ( p - q_{ p } ) \, ( p - q_{ p } + 1 ) / 2
\end{equation}
%
%
rewritten from
(\ref{equation: n even  --  Capacity of vertical median of  W_{01} -- m = 3p+1 --  })
as
\[
          4 \cdot  W'_{01} ( p )  =   \frac{   ( 3 \, p + r_{ p } ) ( 3 \, p +  r_{ p }  +  4 ) }{ 8 }   
\]
whence
%
%
\begin{equation}
\label{equation: n even  -- 4 * W'_{01} ( p )  -- }
          4 \cdot  W'_{01} ( p )  =  \left\{
                        \begin{array}{ll}
                                   ( 3 \, p ) ( 3 \, p  +  4 ) / 8	 		& 	\mbox{ if }		 r_{ p }  = 0  	\\
                                   ( 3 \, p  + 1 ) ( 3 \, p  +  5 ) / 8		& 	\mbox{ if }		 r_{ p }  = 1  	\\
                                   ( 3 \, p  +  2 ) ( 3 \, p  +  6 ) / 8	 	& 	\mbox{ if }		 r_{ p }  = 2  	\\
                                  ( 3 \, p  +  3 ) ( 3 \, p  +  7 ) / 8	 	& 	\mbox{ if }		 r_{ p }  = 3  	\\
                        \end{array}
                                                                                                \right.
\end{equation}
%
%
%
giving the overall capacity of the four  \texttt{N}--\texttt{S}--\texttt{E}--\texttt{W} subregions $ \mathcal{W}'_{01} $ in terms of $ p $.
%
%

For subregion $ \mathcal{W}''_{01} $ we get
%
%
\begin{equation}
\label{equation: n even  --  capacity of b"_(01) -- m = 3p+1 -- }
              \overline{b''}_{01} ( p )   =   ( q_{ p } - 1 )  +   \delta_{ r_{ p } }    \hspace{3mm} \mbox{with} \hspace{3mm} 
              \delta_{ r_{ p } }  =   \Bigl\lfloor  \,\frac { r_{ p }  + 1 \, }  { 3 }  \Bigr\rceil 
\end{equation}
%
%
%
%
\begin{equation}
\label{equation: n even  -- W''_{01} ( p )  -- m = 3p+1 -- }
                W''_{01}( p )  =  \overline{b''}_{01} ( p )  \cdot
                \frac{  ( p - q_{ p }  -  3 )  +  ( p - q_{ p }  -  3  \cdot  \overline{b''}_{01} ( p ) ) }  { 2 } 
\end{equation}
%
%
%
from
(\ref{equation: n even  --  capacity of b"_(0m) -- })
and
(\ref{equation: n even  -- W''_{0m} ( p )  -- })
respectively. It follows that
\[           \overline{b''}_{01} ( p )   =   \left\{
                        \begin{array}{ll}
                                   q_{ p } - 1	& 	\mbox{ if }		 r_{ p }  =  0  	\\
                                   q_{ p } 		& 	\mbox{ otherwise }			\\
                        \end{array}
                                                                                                \right.
\]
and then
\[         W''_{01}( p )   =   \left\{
                        \begin{array}{ll}
                                  (  q_{ p } - 1 ) \cdot  ( 2\, p -  5 \,  q_{ p }  ) / 2		& 	\mbox{ if }		 r_{ p }  =  0  	\\
                                   q_{ p } \cdot	( 2 \, p -  5 \, q_{ p }  -  3 ) / 2		& 	\mbox{ otherwise. }			\\
                        \end{array}
                                                                                                \right.
\]
Replacing $ q_{ p } $ by $ ( p - r_{ p } ) / 4 $ from 
(\ref{equation:  n even  --  p = 4*q + r  -- C00 C01 C02 })
it comes
\[         4 \cdot  W''_{01}( p )   =   \left\{
                        \begin{array}{ll}
                                  (  p - r_{ p }  -  4 ) \cdot  ( 3 \, p +  5 \, r_{ p } ) / 8		& 	\mbox{ if }		 r_{ p }  =  0  	\\
                                    ( p - r_{ p } )    \cdot	  (  3 \, p +  5 \, r_{ p }  - 12 ) / 8		& 	\mbox{ otherwise }			\\
                        \end{array}
                                                                                                \right.
\]
whence
%
%
\begin{equation}
\label{equation: n even  -- 4 * W''_{01} ( p )  -- }
          	4 \cdot  W''_{01} ( p )   =  \left\{
                        \begin{array}{ll}
                               	   3 p \, ( p - 4 ) / 8		& 	\mbox{ if }	 r_{ p }  =  0 		\\
                               	  ( p - 1 ) ( 3 p - 7 ) / 8	& 	\mbox{ if }	 r_{ p }  =  1 		\\
                               	  ( p - 2 ) ( 3 p - 2 ) / 8	& 	\mbox{ if }	 r_{ p }  =  2 		\\
                               	  ( p - 3 ) ( 3 p + 3 ) / 8	& 	\mbox{ if }	 r_{ p }  =  3 		
                        \end{array}
                                                                           \right.
\end{equation}
%
%
giving the overall capacity of the four  \texttt{N}--\texttt{S}--\texttt{E}--\texttt{W} subregions $ \mathcal{W}''_{01} $ in terms of $ p $.

\begin{figure}
\centering
\includegraphics[width=12cm]{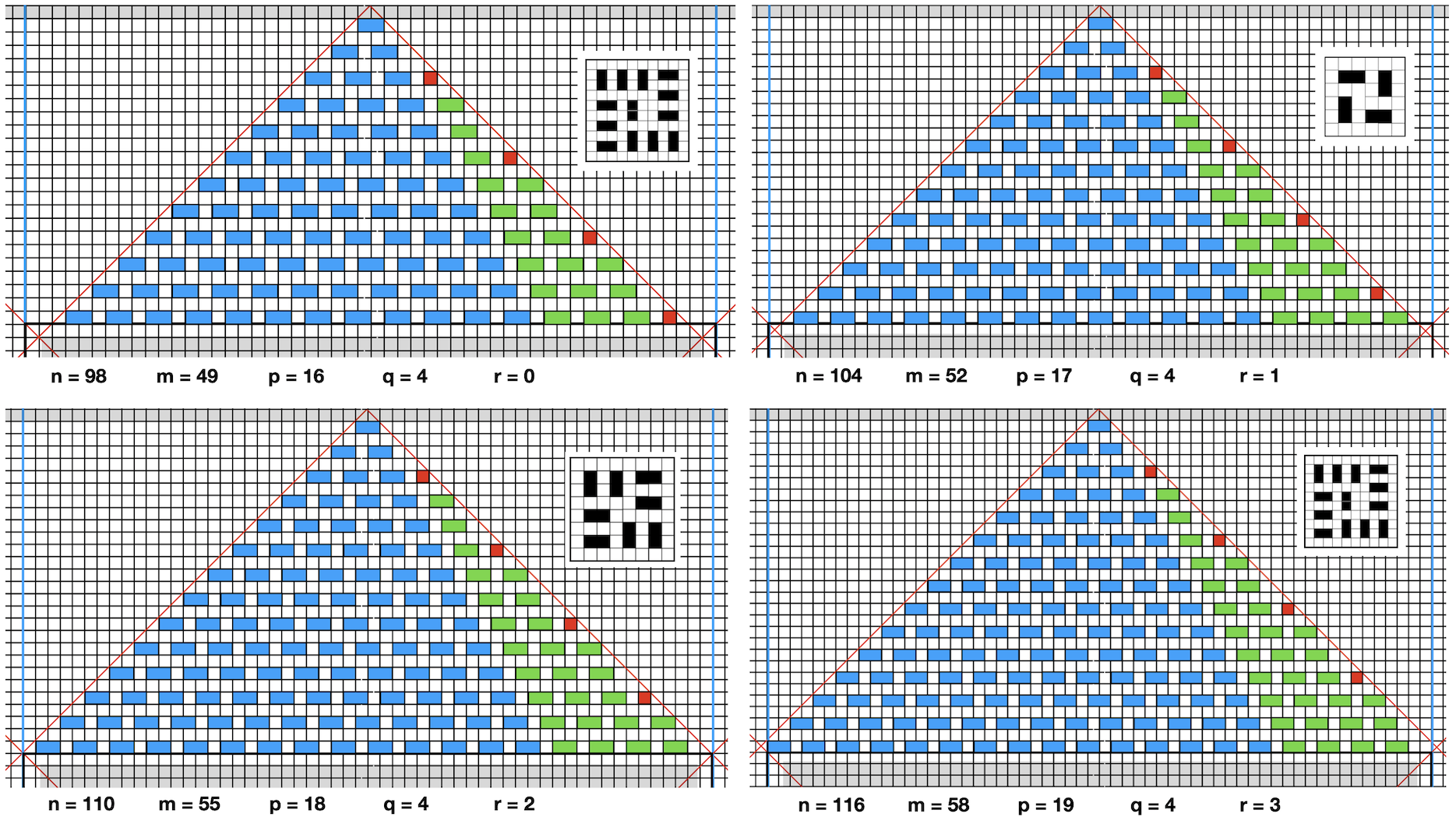}  
\caption{A 4--sequence of  $ \mathcal{D}_{ n } $ expansion in \texttt{Class} $ 01 $.
The class of the square core 
$
       \mathcal{S}_{ f_{ 01 }( m )} 
$
(refer to Fig.\,\ref{Figure: Dominos in the Square -- n odd -- n even-- })
evolves  as 
$ 
           ( \dot{ 01 } ) \rightarrow ( \dot{ 02 } ) \rightarrow ( \dot{ 00 })  \rightarrow ( \dot{ 01 } )
$ 
(highlighted in insets). A new sequence starts with $  p  \equiv_4  0 $ whereas  
$
         \mathcal{W}_{01}
$
remains unchanged.
}  
\label{Figure: C01-P16P19 }
\end{figure}

Finally, from
(\ref{equation: n even  -- Xi-- m=3p+1 -- r = 0 })--(\ref{equation: n even  -- Xi-- m=3p+1 -- r = 3 })  
for
$
           \xi_{ f_{01} ( m ) }
$
and from
(\ref{equation: n even  -- 4 * W'_{01} ( p )  -- })
and 
(\ref{equation: n even  -- 4 * W''_{01} ( p )  -- })
for
$
         W'_{01} ( p )
$
and
$
         W''_{01} ( p )
$
it would easily come
%
%
\begin{equation}
\label{equation: n even  -- Psi_n ( n )  in  C01 }
            \psi_{n}   =    \left\{
                        \begin{array}{ll}
                                   ( n^2 + 2 n + 4 ) / 12	& 	\mbox{ if }		 r_{ p }  =  0  	\\
                                   ( n^2 + 2 n + 16 ) / 12	& 	\mbox{ if }		 r_{ p }  =  1  	\\
                                   ( n^2 + 2 n + 16 ) / 12	& 	\mbox{ if }		 r_{ p }  =  2  	\\
                                   ( n^2 + 2 n + 4 ) / 12	& 	\mbox{ if }		 r_{ p }  =  3  	
                        \end{array}
                                                                                                \right.
\end{equation}
%
%
with $ m = 3 p + 1 $ and $ n = 2 m $.
With
$
             \psi_{ 2 }  = 1
$
the expansion rate of  $ \psi_{ n } $ can also be obtained stepwise. 
From
(\ref{equation: n even  -- W'_{0m} ( p ) - W'_{0m} ( p - 1 ) -- })
and
(\ref{equation: n even  -- W''_{0m} ( p ) - W''_{0m} ( p - 1 ) -- })
it comes
%
%
\begin{equation}
\label{equation: n even  -- W'_{01} ( p ) - W'_{01} ( p - 1 ) -- m = 3p+1 -- }
          W'_{01} ( p )  -   W'_{01} ( p - 1 )  =  \left\{
                        \begin{array}{ll}
                                   0				& 	\mbox{ if }		 r_{ p }  =  0 	\\
                                   p  -  q_{ p } 	& 	\mbox{ otherwise }		 
                        \end{array}
                                                                                                \right.
\end{equation}
%
%
%
%
\begin{equation}
\label{equation: n even  -- W''_{01} ( p ) - W''_{01} ( p - 1 ) -- m = 3p+1 -- }
         	W''_{01} ( p )  -   W''_{01} ( p - 1 )  =  \left\{
                        \begin{array}{ll}
                               	0							 		& 	\mbox{ if }	    r_{ p } = 0 	\\
                                  q_{ p } - 1  +   \delta_{ r_{ p } }  	=   q_{ p }		& 	\mbox{ otherwise }		 
                         \end{array}
                                                                           \right.
\end{equation}
%
%
and in
Tab.\,\ref{Table: Expansion Psi -- m = 3p+1 },
the sum of the three terms resulting from relations
(\ref{equation: n even  -- DXi-- m=3p+1 -- r = 0 }--\ref{equation: n even  -- DXi-- m=3p+1 -- r = 3 })
for 
$
\Delta \xi_{ f_{01} ( m ) }, 
$
from
(\ref{equation: n even  -- W'_{01} ( p ) - W'_{01} ( p - 1 ) -- m = 3p+1 -- })
for 
$
        \Delta  W'_{01} ( p )
$
and from
(\ref{equation: n even  -- W''_{01} ( p ) - W''_{01} ( p - 1 ) -- m = 3p+1 -- })
for 
$
        \Delta  W''_{01} ( p ) 
$
gives the result.
Moreover, by adding backwards the four terms in the last column of 
Tab.\,\ref{Table: Expansion Psi -- m = 3p+1 },
it follows that 
\[
              \psi_n -  \psi_{ n - 24 }  =  ( n - 2 )  + ( ( n - 6 ) - 3 )  + ( ( n - 12 ) - 2 )  +  ( ( n - 18 ) - 1 ) 
\]
$
             =  4\, n - ( \, 2 + 9 + 14 + 19 \, )
$
whence the induction
%
%
\begin{equation}
\label{equation: n even -- m == 1  -- Psi(n) - Psi( n - 24 ) }
            \psi_{ 2 }  =   1 \, ;   \hspace{2mm}  \psi_n =  \psi_{ n - 24 }\ +   4\,( n - 11 )     \hspace{5mm}  ( p  \ge  4 )
\end{equation}
%
%
within any  
$
          p  \rightarrow  p + 4
$
cycle.
A 4--sequence of expansion is illustrated in  
Fig.\,\ref{Figure: C01-P16P19 }.
%
%
%
%
The values of $ \psi_{n} $ in \texttt{Class} $ 01 $ are displayed in 
Table\,\ref{Table: Domino enumeration -- Class 01 }.
%
%
%
\begin{table}[tb!]
\centering
\caption{Expansion in \texttt{Class} $ 01 $:  
$
           p \, \in  \, \mathbb{N} \,  ;  m  \, = 3 \, p + 1    \, ;   n \, =  2 \, m  \, ;  \ - \  \hspace{1mm}  \psi_{ 2 }  = 1 \, ; \newline
           \mbox{ }  \hspace{13mm} 
           r_{ p } \equiv  p  \pmod 4
$
}
\vspace{2mm} 
\label{Table: Expansion Psi -- m = 3p+1 }
\begin{tabular}{|c||c|c|c||c|}
\hline
\rule{0pt}{10pt}
        $   r_{ p }    $	&    
        		$ \Delta \xi_{ f_{01} ( m ) } $	&   
							$ 4\, \Delta  W'_{01} ( p ) $	&   
														$ 4\, \Delta  W''_{01} ( p ) $	& 
																			$	\psi_{n}-\psi_{n-6}	$  \\   [0.5ex] 
\hline\hline	
	 0 	&	$      6\, p  	$	&    $			0			$ 	&    $			0		$	&    	$	n - 2		$	\\   
 	 1 	&	$   2 \, p + 1 	$	&    $  4 \, ( p - q_{ p } )		$  	&    $      4 \, q_{ p }  		$ 	&     	$	n - 1		$	\\   	 	 
	 2  	&	$       2 \, p	$	&    $  4 \, ( p - q_{ p } )		$  	&    $      4 \, q_{ p }		$	&	$	n - 2		$	\\
	 3  	&	$   2 \, p - 1	$	&    $  4 \, ( p - q_{ p } )		$  	&    $      4 \, q_{ p }		$ 	&	$	n - 3		$	\\  \hline	   	 
\end{tabular}
\end{table}
%
%
%
%
%
\subsection{ \texttt{Class} $ 02 $ --- $ n $ even \, \&  \, $ m \equiv  2  \pmod 3 $}
\label{Subsection: n  even --- m = 3p+2-- }
%
%
%
The basic parameters of \texttt{Class} $ 02 $ verify: 
\[
		p \, \in  \, \mathbb{N} \,  ;  \hspace{5mm} m  \, = 3 \, p \, + \, 2 \, ;  \hspace{5mm} n \, =  2 \, m \, ;
\]
%
%
\subsubsection{Capacity and expansion rate of the square core $  \mathcal{S}_{  f_{ 02 }  } $ }                 
%
%
%
From
(\ref{equation:  n even  --  f_{ 0 m }( m ) -- C00 C01 C02 }--\ref{equation:  n even  --  p = 4*q + r  })
$
            f_{02} ( m ) = m  + 1 - r_{ p - 1 } 	
$
and we get the following four cases
%
%
%
\begin{itemize}
%
%
\item  $   r_{ p - 1 }  = 0      \Rightarrow    f_{02} ( m ) =   m + 1  $     
%
%
and
$ 
          ( m + 1 ) / 2  \in \dot{ 0 }   \hspace{1mm} \mbox{in} \hspace{1mm}    \mathbb{Z} / 3 \, \mathbb{Z}
$
then  $ \xi_{ m + 1 } $ follows from
(\ref{equation: Xi_n -- n even -- m = 3*p})
whence
%
%
\begin{equation}
\label{equation: n even  -- Xi-- m=3p+2 -- r = 0 }
         \xi_{ m + 1 }  =   \frac{ ( m + 1 ) \, ( m + 3 )}{6} 
\end{equation}
%
%
%
\item  $   r_{ p - 1 }  = 1      \Rightarrow    f_{02} ( m ) =   m  $     
%
%
and 
$ 
         m / 2  \in \dot{ 1 }   \hspace{1mm}   
$
then  $ \xi_{ m } $ follows from
(\ref{equation: Xi_n -- n even -- m = 3*p+1})
whence
%
%
\begin{equation}
\label{equation: n even  -- Xi-- m=3p+2 -- r = 1 }
         \xi_{ m }   =  \frac{ m ( m + 2 ) - 2 }{6}   
\end{equation}
%
%
%
%
\item  $   r_{ p - 1 }  = 2      \Rightarrow    f_{02} ( m ) =   m - 1  $     
%
%
and 
$ 
         ( m - 1 ) / 2  \in \dot{ 2 }   \hspace{1mm}   
$
then  $ \xi_{ m - 1 } $ follows from
(\ref{equation: Xi_n -- n even -- m = 3*p+2})
whence
%
%
\begin{equation}
\label{equation: n even  -- Xi-- m=3p+2 -- r = 2 }
         \xi_{ m - 1 }   =  \frac{ ( m - 1 ) ( m + 1 ) }{6}   
\end{equation}
%
%
%
\item  $   r_{ p - 1 }  = 3      \Rightarrow    f_{02} ( m ) =   m - 2  $     
%
%
and 
$ 
         ( m - 2 ) / 2  \in \dot{ 0 }   \hspace{1mm}   
$
then  $ \xi_{ m - 2 } $ follows from
(\ref{equation: Xi_n -- n even -- m = 3*p})
whence
%
%
\begin{equation}
\label{equation: n even  -- Xi-- m=3p+2 -- r = 3 }
         \xi_{ m - 2 }  =  \frac{ m ( m - 2 ) }{6}      
\end{equation}
%
%
%
\end{itemize}
%
%
giving the capacity of the square core. In the same way, from
(\ref{equation:  n even  --  Square Core Elongation })
we get the following four cases
%
%
%
\begin{itemize}
%
%
\item  $   r_{ p + 1 }  = 0      \Rightarrow    f_{02} ( m - 3 ) = f_{02} ( m ) - 6   $     
%
%
whence from
(\ref{equation: Induction Xi -- n even }) 
%
%
\begin{equation}
\label{equation: n even  -- DXi-- m=3p+2 -- r = 0 }
         \xi_{ m + 1 }  -  \xi_{ ( m + 1 ) - 6 }  =  2 ( m - 1 )=  6 \, p  +  2
\end{equation}
%
%
%
\item  $   r_{ p - 1 }  = 1      \Rightarrow    f_{02} ( m - 3 ) =  f_{02} ( m ) - 2 = m - 2   $     
%
%
and 
$ 
          ( m - 2  ) / 2  \in \dot{ 0 }  
$
then  $ \xi_{ m - 2 } $ follows from
(\ref{equation: Xi_n -- n even -- m = 3*p})
whence
%
%
\begin{equation}
\label{equation: n even  -- DXi-- m=3p+2 -- r = 1 }
         \xi_{ m }  -  \xi_{ m - 2 }  =  \frac{ m ( m + 2 ) - 2 }{6}  -   \frac{ m ( m - 2 ) }{6}  =  2 \, p  + 1
\end{equation}
%
%
%
%
\item  $   r_{ p - 1 }  = 2      \Rightarrow     f_{02} ( m - 3 ) =  f_{02} ( m ) - 2 = m - 3   $     
%
%
and 
$ 
          ( m - 3 ) / 2  \in \dot{ 1 }  
$
then  $ \xi_{ m - 3 } $ follows from
(\ref{equation: Xi_n -- n even -- m = 3*p+1})
whence
%
%
\begin{equation}
\label{equation: n even  -- DXi-- m=3p+2 -- r = 2 }
         \xi_{ m - 1 }  -  \xi_{ m - 3 }  =  \frac{ ( m - 1 ) ( m + 1 ) }{6}  -   \frac{ ( m - 3 )( m - 1 ) - 2 }{6}  =  2 \, p  + 1  
\end{equation}
%
%
%
%
\item  $   r_{ p - 1 }  = 3      \Rightarrow    f_{02} ( m - 3 ) =  f_{02} ( m ) - 2 = m - 4  $     
%
%
and 
$ 
          ( m - 4 ) / 2  \in \dot{ 2 }  
$
then  $ \xi_{ m - 4 } $ follows from
(\ref{equation: Xi_n -- n even -- m = 3*p+2})
whence
%
%
\begin{equation}
\label{equation: n even  -- DXi-- m=3p+2 -- r = 3 }
         \xi_{ m - 2 }  -  \xi_{ m - 4 }  =  \frac{ m ( m - 2 ) }{6}  -   \frac{ ( m - 4 )( m - 2 ) }{6}  =  2 \, p    
\end{equation}
%
%
%
\end{itemize}
%
%
%
giving the expansion rate for the capacity in the square core.
%
%
%
\subsubsection{Capacity and expansion  of  $ \mathcal{D}_{ n } $ in \texttt{Class} $ 02 $ }                 
%
%
The capacity $ \psi_{ n } $ of  $ \mathcal{D}_{ n } $ in \texttt{Class} $ 02 $ is simply given by
%
%
\begin{equation}
\label{equation: Psi_n in Class 02 }
      \psi_{n}  =   \xi_{ f_{02}( m ) }  + 4 \, ( W'_{02}( p )  +  W''_{02}( p ) )
\end{equation}
%
%
%
and
$
           \xi_{ f_{02} ( m ) }
$
results from
(\ref{equation: n even  -- Xi-- m=3p+2 -- r = 0 })--(\ref{equation: n even  -- Xi-- m=3p+2 -- r = 3 }).
In \texttt{Class} $ 02 $ from
(\ref{equation:  n even  --  p = 4*q + r  -- C00 C01 C02 })
it comes
%
%
\begin{equation}
\label{equation: n even  --  Capacity of vertical median of  W_{02} -- m = 3p+2 -- }
            p - q_{ p - 1 }  =  3 \, q_{ p - 1 } + r_{ p - 1 } + 1   \hspace{2mm} \mbox{or} \hspace{2mm} 
            4 \, ( p - q_{ p - 1 }  )  =  3 \, p +  r_{ p - 1 }  + 1
\end{equation}
%
%
and from
(\ref{equation: n even  -- W'_{0m} ( p )  -- })
%
%
%
\begin{equation}
\label{equation: n even  -- W'_{02} ( p )  -- m = 3p+2 -- }
            W'_{02} ( p )  =  ( p - q_{ p - 1 } ) \, ( p - q_{ p - 1 } + 1 ) / 2
\end{equation}
%
%
rewritten from
(\ref{equation: n even  --  Capacity of vertical median of  W_{02} -- m = 3p+2 --  })
as
\[
          4 \cdot  W'_{02} ( p )  =   \frac{   ( 3 \, p + r_{ p - 1 } + 1 ) ( 3 \, p +  r_{ p - 1 }  +  5 ) }{ 8 }   
\]
whence
%
%
\begin{equation}
\label{equation: n even  -- 4 * W'_{02} ( p )  -- }
          4 \cdot  W'_{02} ( p )  =  \left\{
                        \begin{array}{ll}
                                   ( 3 \, p + 1 ) ( 3 \, p  +  5 ) / 8	 	& 	\mbox{ if }		 r_{ p - 1 }  = 0  	\\
                                   ( 3 \, p  + 2 ) ( 3 \, p  +  6 ) / 8		& 	\mbox{ if }		 r_{ p - 1 }  = 1  	\\
                                   ( 3 \, p  +  3 ) ( 3 \, p  +  7 ) / 8	 	& 	\mbox{ if }		 r_{ p - 1 }  = 2  	\\
                                   ( 3 \, p  +  4 ) ( 3 \, p  +  8 ) / 8	 	& 	\mbox{ if }		 r_{ p - 1 }  = 3  	\\
                        \end{array}
                                                                                                \right.
\end{equation}
%
%
%
giving the overall capacity of the four  \texttt{N}--\texttt{S}--\texttt{E}--\texttt{W} subregions $ \mathcal{W}'_{02} $ in terms of $ p $.

For subregion $ \mathcal{W}''_{02} $ we get
%
%
\begin{equation}
\label{equation: n even  --  capacity of b"_(02) -- m = 3p+2 -- }
              \overline{b''}_{02} ( p )   =   ( q_{ p - 1 } - 1 )  +   \delta_{ r_{ p - 1 } }    \hspace{3mm} \mbox{with} \hspace{3mm} 
              \delta_{ r_{ p - 1 } }  =   \Bigl\lfloor  \,\frac { r_{ p - 1 }  +  2 \, }  { 3 }  \Bigr\rceil 
\end{equation}
%
%
%
%
\begin{equation}
\label{equation: n even  -- W''_{02} ( p )  -- m = 3p+2 -- }
                W''_{02}( p )  =  \overline{b''}_{02} ( p )  \cdot
                \frac{  ( p - q_{ p - 1 }  -  3 )  +  ( p - q_{ p - 1 }  -  3  \cdot  \overline{b''}_{02} ( p ) ) }  { 2 } 
\end{equation}
%
%
%
from
(\ref{equation: n even  --  capacity of b"_(0m) -- })
and
(\ref{equation: n even  -- W''_{0m} ( p )  -- })
respectively. It follows that
\[           \overline{b''}_{02} ( p )   =   \left\{
                        \begin{array}{ll}
                                   q_{ p - 1 }  + 1		& 	\mbox{ if }		 r_{ p - 1 }  =  3  	\\
                                   q_{ p - 1 } 		& 	\mbox{ otherwise }			\\
                        \end{array}
                                                                                                \right.
\]
and then
\[         W''_{02}( p )   =   \left\{
                        \begin{array}{ll}
                                  (  q_{ p - 1 } + 1 ) \cdot  ( 2\, p -  5 \,  q_{ p - 1 } - 6 ) / 2		& 	\mbox{ if }		 r_{ p - 1 }  =  3  	\\
                                   q_{ p - 1 } \cdot	( 2 \, p -  5 \, q_{ p - 1 }  -  3 ) / 2		& 	\mbox{ otherwise. }			\\
                        \end{array}
                                                                                                \right.
\]
Replacing $ q_{ p - 1 } $ by $ ( p - 1 - r_{ p - 1 } ) / 4 $ from 
(\ref{equation:  n even  --  p = 4*q + r  -- C00 C01 C02 })
it comes
\[          4 \cdot  W''_{02}( p )   =   \left\{
                        \begin{array}{ll}
                                    p   \cdot  ( 3 \, p -  4  ) / 8										& 	\mbox{ if }		 r_{ p - 1}  =  3  	\\
                                   ( p - 1 - r_{ p - 1 } )    \cdot   ( 3 \, p +  5 \,  ( 1 +  r_{ p - 1 })  -  12 ) / 8	& 	\mbox{ otherwise }			\\
                        \end{array}
                                                                                                \right.
\]
whence
%
%
\begin{equation}
\label{equation: n even  -- 4 * W''_{02} ( p )  -- }
          	4 \cdot  W''_{02} ( p )   =  \left\{
                        \begin{array}{ll}
                               	  ( p - 1 ) ( 3 p - 7 ) / 8 	& 	\mbox{ if }	 r_{ p - 1 }  =  0 		\\
                               	  ( p - 2 ) ( 3 p - 2 ) / 8	& 	\mbox{ if }	 r_{ p - 1 }  =  1 		\\
                               	  ( p - 3 ) ( 3 p + 3 ) / 8	& 	\mbox{ if }	 r_{ p - 1 }  =  2 		\\
                                 	  p  \, ( 3 \, p -  4  ) / 8		& 	\mbox{ if }	 r_{ p - 1 }  =  3 		
                        \end{array}
                                                                           \right.
\end{equation}
%
%
giving the overall capacity of the four  \texttt{N}--\texttt{S}--\texttt{E}--\texttt{W} subregions $ \mathcal{W}''_{02} $ in terms of $ p $.

\begin{figure}
\centering
\includegraphics[width=12cm]{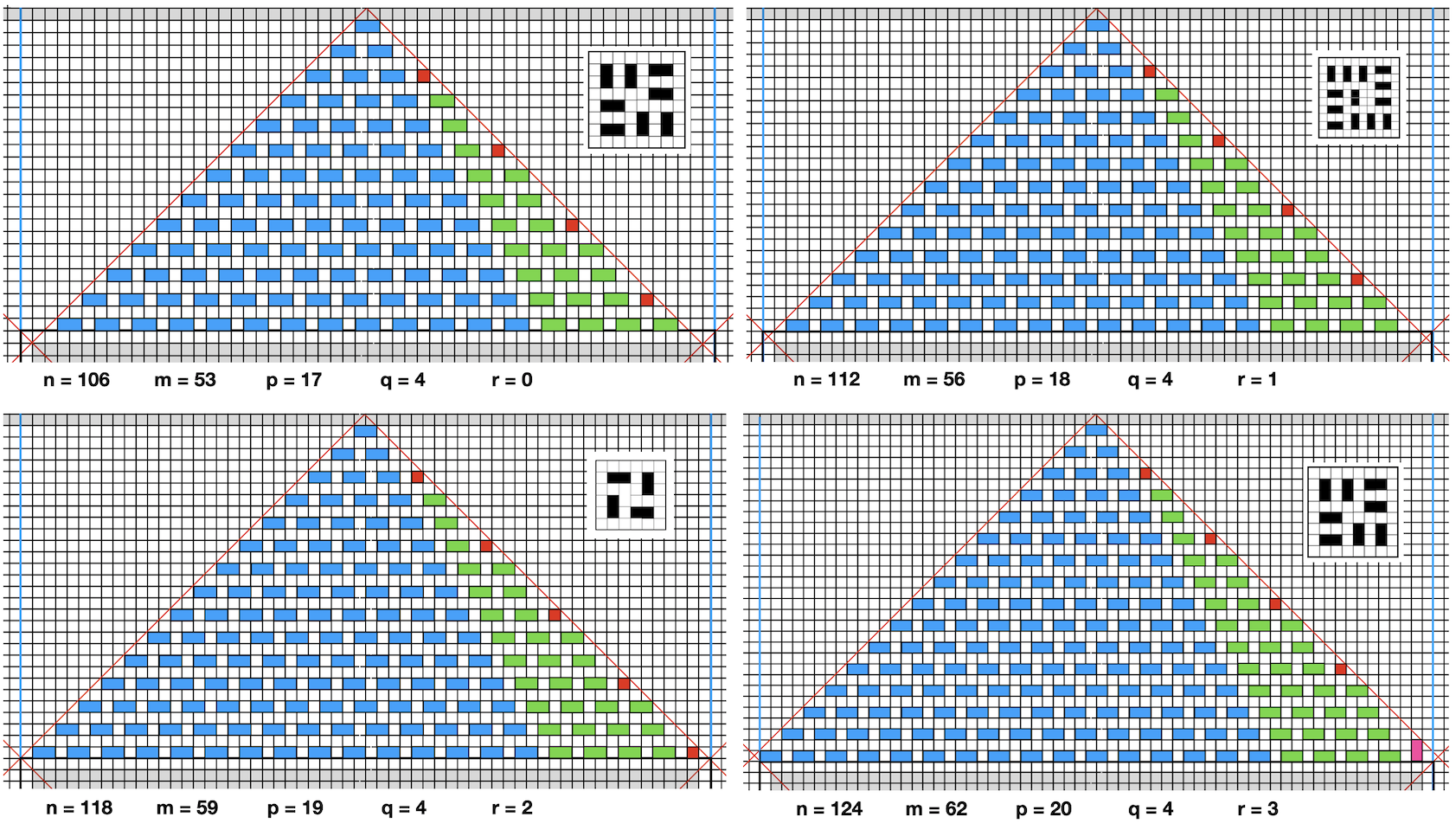}  
\caption{A 4--sequence of  $ \mathcal{D}_{ n } $ expansion in \texttt{Class} $ 02 $.
The class of the square core 
$
       \mathcal{S}_{ f_{ 02 }( m )} 
$
(refer to Fig.\,\ref{Figure: Dominos in the Square -- n odd -- n even-- })
evolves  as 
$ 
           ( \dot{ 00 } ) \rightarrow ( \dot{ 01 } ) \rightarrow ( \dot{ 02 })  \rightarrow ( \dot{ 00 } )
$ 
(highlighted in insets). A new sequence starts with $  p - 1  \equiv_4  0 $ whereas  
$
         \mathcal{W}_{02}
$
remains unchanged.
An additional, unstable vertical domino arises at $ r_{ p - 1 } = 3 $. It will return to its horizontal position at the start of the next sequence, restoring the emerging lacunar void at $ r_{ p - 1 } = 2 $.
}  
\label{Figure: C02-P17P20 }
\end{figure}

Finally, from
(\ref{equation: n even  -- Xi-- m=3p+2 -- r = 0 })--(\ref{equation: n even  -- Xi-- m=3p+2 -- r = 3 })  
for
$
           \xi_{ f_{02} ( m ) }
$
and from
(\ref{equation: n even  -- 4 * W'_{02} ( p )  -- })
and 
(\ref{equation: n even  -- 4 * W''_{02} ( p )  -- })
for
$
         W'_{02} ( p )
$
and
$
         W''_{02} ( p )
$
it would easily come
%
%
\begin{equation}
\label{equation: n even  -- Psi_n ( n )  in  C02 }
            \psi_{n}   =    \left\{
                        \begin{array}{ll}
                                   ( n^2 + 2 n + 24 ) / 12	& 	\mbox{ if }		 r_{ p - 1 }  =  0  	\\
                                   ( n^2 + 2 n + 12 ) / 12	& 	\mbox{ if }		 r_{ p - 1 }  =  1  	\\
                                   ( n^2 + 2 n  ) / 12		& 	\mbox{ if }		 r_{ p - 1 }  =  2  	\\
                                   ( n^2 + 2 n + 24 ) / 12	& 	\mbox{ if }		 r_{ p  - 1}  =  3  
                        \end{array}
                                                                                                \right.
\end{equation}
%
%
with $ m = 3 p + 2 $ and $ n = 2 m $.
With
$
             \psi_{ 4 }  =  2
$
the expansion rate of  $ \psi_{ n } $ can also be obtained stepwise. 
From
(\ref{equation: n even  -- W'_{0m} ( p ) - W'_{0m} ( p - 1 ) -- })
and
(\ref{equation: n even  -- W''_{0m} ( p ) - W''_{0m} ( p - 1 ) -- })
it comes
%
%
\begin{equation}
\label{equation: n even  -- W'_{02} ( p ) - W'_{02} ( p - 1 ) -- m = 3p+2 -- }
          W'_{02} ( p )  -   W'_{02} ( p - 1 )  =  \left\{
                        \begin{array}{ll}
                                   0				& 	\mbox{ if }		 r_{ p - 1 }  =  0 	\\
                                   p  -  q_{ p - 1 } 	& 	\mbox{ otherwise }		 
                        \end{array}
                                                                                                \right.
\end{equation}
%
%
%
%
\begin{equation}
\label{equation: n even  -- W''_{02} ( p ) - W''_{02} ( p - 1 ) -- m = 3p+2 -- }
         	W''_{02} ( p )  -   W''_{02} ( p - 1 )  =  \left\{
                        \begin{array}{ll}
                               	0							 	& 	\mbox{ if }	    r_{ p - 1 } = 0 	\\
                                  q_{ p - 1 } - 1  +   \delta_{ r_{ p - 1 } }  	& 	\mbox{ otherwise }		 
                         \end{array}
                                                                           \right.
\end{equation}
%
%
and in
Tab.\,\ref{Table: Expansion Psi -- m = 3p+2 },
the sum of the three terms resulting from relations
(\ref{equation: n even  -- DXi-- m=3p+2 -- r = 0 }--\ref{equation: n even  -- DXi-- m=3p+2 -- r = 3 })
for 
$
\Delta \xi_{ f_{02} ( m ) }, 
$
from
(\ref{equation: n even  -- W'_{02} ( p ) - W'_{02} ( p - 1 ) -- m = 3p+2 -- })
for 
$
        \Delta  W'_{02} ( p )
$
and from
(\ref{equation: n even  -- W''_{02} ( p ) - W''_{02} ( p - 1 ) -- m = 3p+2 -- })
for 
$
        \Delta  W''_{02} ( p ) 
$
gives the result.
Moreover, by adding backwards the four terms in the last column of 
Tab.\,\ref{Table: Expansion Psi -- m = 3p+2 },
it follows that 
\[
              \psi_n -  \psi_{ n - 24 }  =  ( n - 2 )  +  ( n - 6 )  + ( ( n - 12 ) - 3 )  +  ( ( n - 18 ) - 3 ) 
\]
$
             =  4\, n - ( \, 2 + 6 + 15 + 21 \, )
$
whence the induction
%
%
\begin{equation}
\label{equation: n even -- m == 2  -- Psi(n) - Psi( n - 24 ) }
            \psi_{ 4 }  =   2 \, ;   \hspace{2mm}  \psi_n =  \psi_{ n - 24 }\ +   4\,( n - 11 )     \hspace{5mm}  ( p  >  4 )
\end{equation}
%
%
within any  
$
          p  \rightarrow  p + 4
$
cycle.
A 4--sequence of expansion is illustrated in  
Fig.\,\ref{Figure: C02-P17P20 }.
The values of $ \psi_{n} $ in \texttt{Class} $ 02 $ are displayed in 
Table\,\ref{Table: Domino enumeration -- Class 02 }.
%
%
%
%
%
%
%
\begin{table}[tb!]
\centering
\caption{Expansion in \texttt{Class} $ 02 $:  
$
           p \, \in  \, \mathbb{N} \,  ;  m  \, = 3 \, p + 2    \, ;   n \, =  2 \, m  \, ;  \ - \  \hspace{1mm}  \psi_{ 4 }  = 2 \, ; \newline
           \mbox{ }  \hspace{13mm} 
            r_{ p - 1 } \equiv  p - 1  \pmod 4
$
}
\vspace{1mm} 
\label{Table: Expansion Psi -- m = 3p+2 }
\begin{tabular}{|c||c|c|c||c|}
\hline
\rule{0pt}{10pt}
        $   r_{ p - 1 }    $	&    
        		$ \Delta \xi_{ f_{02} ( m ) } $	&   
							$ 4\, \Delta  W'_{02} ( p ) $	&   
														$ 4\, \Delta  W''_{02} ( p ) $	& 
																			$	\psi_{n}-\psi_{n-6}	$  \\   [0.5ex] 
\hline\hline	
	 0 	&	$   6\, p + 2  	$	&    $			0			$ 	&    $			0		$	&    	$	n - 2		$	\\   
 	 1 	&	$   2 \, p + 1 	$	&    $  4 \, ( p - q_{ p - 1 } )		$  	&    $      4 \, q_{ p - 1 }  	$ 	&     	$	n - 3		$	\\   	 	 
	 2  	&	$   2 \, p + 1	$	&    $  4 \, ( p - q_{ p - 1 } )		$  	&    $      4 \, q_{ p - 1 }	$	&	$	n - 3		$	\\
	 3  	&	$      2 \, p 	$	&    $  4 \, ( p - q_{ p - 1 } )		$  	&    $ 4 \, ( q_{ p - 1 } + 1 )	$ 	&	$	   n 		$	\\  \hline	   	 
\end{tabular}
\end{table}
%
%
%
%
%
\section{Injection of Disorder}
\label{Section: Injection of Disorder}
%
%
Various scenarios of injection, depending on the class 
 $
 { \,  \overline{n} \,  \overline{m}  } 
 $
are now examined on a case--by--case basis. The resulting value
$
               \overline{\psi}_n  \ge   \psi_n  
$
will set the {\em upper} bound for $ \psi_n $ listed in the last column of each of the six 
$
          C{ \,  \overline{n} \,  \overline{m}  } 
$
tables.

The problem comes from the fact that our {\em axiomatic} construction is sub--optimal but not optimal.
In a disordered system {\em two} lacunar voids are likely to join together, that yields a potential of {\em one} additional domino. It could be helpful to refer back to
Fig.\,\ref{Figure: C1-n7-n29 -- n odd }
and 
Figs.\,\ref{Figure: C10-N31N49 }--\ref{Figure: C12-N35N53 }
for $ n$ odd
and to
Fig.\,\ref{Figure: C0-n6-n28 -- n even }
and 
Figs.\,\ref{Figure: C00-P15P18 }--\ref{Figure: C02-P17P20 }
for $ n $ even.
%
%
\subsection{Injection of Disorder -- $n$ odd}                     
%
%
A first observation of
Fig.\,\ref{Figure: C1-n7-n29 -- n odd }
shows that, for $ p $ even, a slight transformation of ``ridge--flattening'' can be carried out on the 2--fold tip of the wedge --\,by rotating both dominos\,-- thus releasing {\em four} \texttt{N}--\texttt{S}--\texttt{E}--\texttt{W} lacunar voids in all, giving a potential of {\em two} additional dominos. This case ($ \mu_p = 0 $) is illustrated thereafter.

The second observation is related to the central pattern of the square core: a vacant space, with or without lacunar void, is likely to cause a deficit in the whole.
Referring to
Fig.\,\ref{Figure: Dominos in the Square -- n odd -- n even-- }
such a situation exists either for 
$
                \mathcal{S}_{ 1 }
$
or for
$
                \mathcal{S}_{ 3 }
$
--\,condition noted briefly 
$
               ( \mbox{s}_1 \lor \mbox{s}_3 )
$
\hspace{ -1mm }-- while this is not the case for
$
                \mathcal{S}_{ 5 }.
$

Now, for $ \mu_p = 0 $ any deficiency of the core will be compensated by a ridge--flattening.
It follows that a non--optimality problem will be induced by the conjunction  
$
               ( \mbox{s}_1 \lor \mbox{s}_3 )  \land  ( \mu_p = 1 ).
$
This situation will occur for \texttt{Class} $ 10 $ and  \texttt{Class} $ 11 $.
It can be released either by {\em enlarging} or by {\em shrinking} the core size.
%
%
\subsubsection{ In \texttt{Class} $ 10 $ }                 
%
%
%
The following two cases are considered.
%
%
\begin{figure}
\centering
\includegraphics[width=7.2cm]{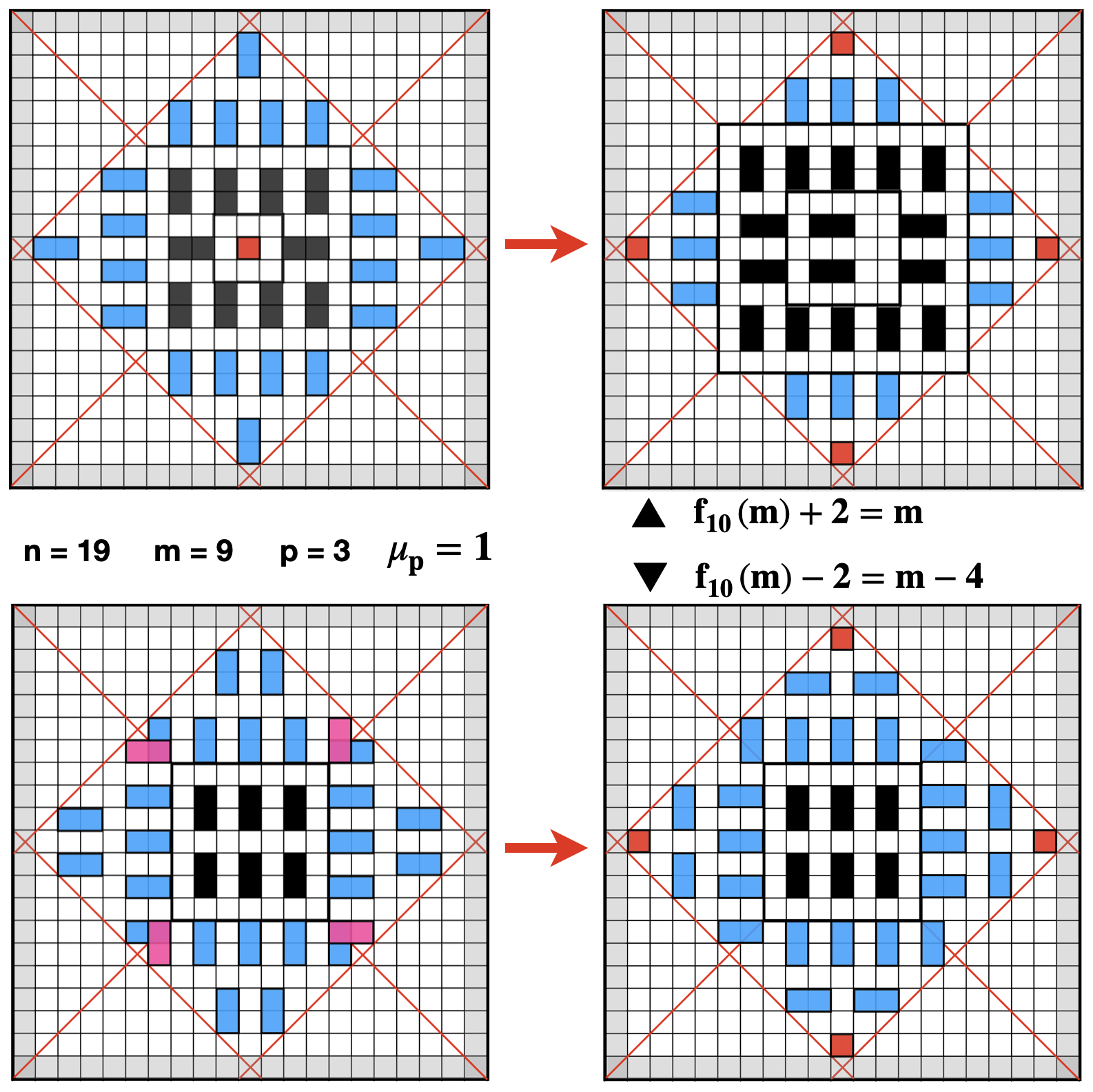}  
\caption{Injection of disorder in Class $ 10 $ for $ \mu_p = 1 $ -- 
 ($ \uparrow $)  Enlarging the square core from size  
 $
        f_{ 10 } ( m )  =  m  - 2   
 $
 to size
 $
       f_{ 10 } ( m ) + 2  =  m 
 $
with expansion  
 $         
         \xi_{ m }  -  \xi_{ m - 2 }  =   + \, 2 \, p   =  +  6
$
for this case;
overall deficit of $ 2 p + 2 = 8 $ for the wedges; overall loss of 2 dominos altogether; overall gain of 4 lacunar voids giving a potential of 2 additional dominos. 
The final additive potential is {\em zero}.\newline
 ($ \downarrow $)  Shrinking the square core from size  
 $
        f_{ 10 } ( m )     
 $
to size
 $
       f_{ 10 } ( m ) - 2  =  m  - 4
 $
with expansion
$         
         \xi_{ m - 4 }  -  \xi_{ m - 2 }  =   - \, 2 \, p  +  2   =  -  4
$
for this case;
disregarding the domino conflict on the diagonal, the new wedge is a copy of the wedge at $ p + 1 $
(Fig.\,\ref{Figure: C1-n7-n29 -- n odd });
overall gain of $ 2 p - 2 = + 4 $ for the wedges; the overall actual gain on $  \psi_n $ is {\em zero};
ridge--flattening releasing 4 lacunar voids.
Ultimate potential capacity of $ \psi_n  +  2 $.
 } 
\label{Figure: IoD-C10-P3 }
%
%
%
%
%
\centering
\includegraphics[width=7.2cm]{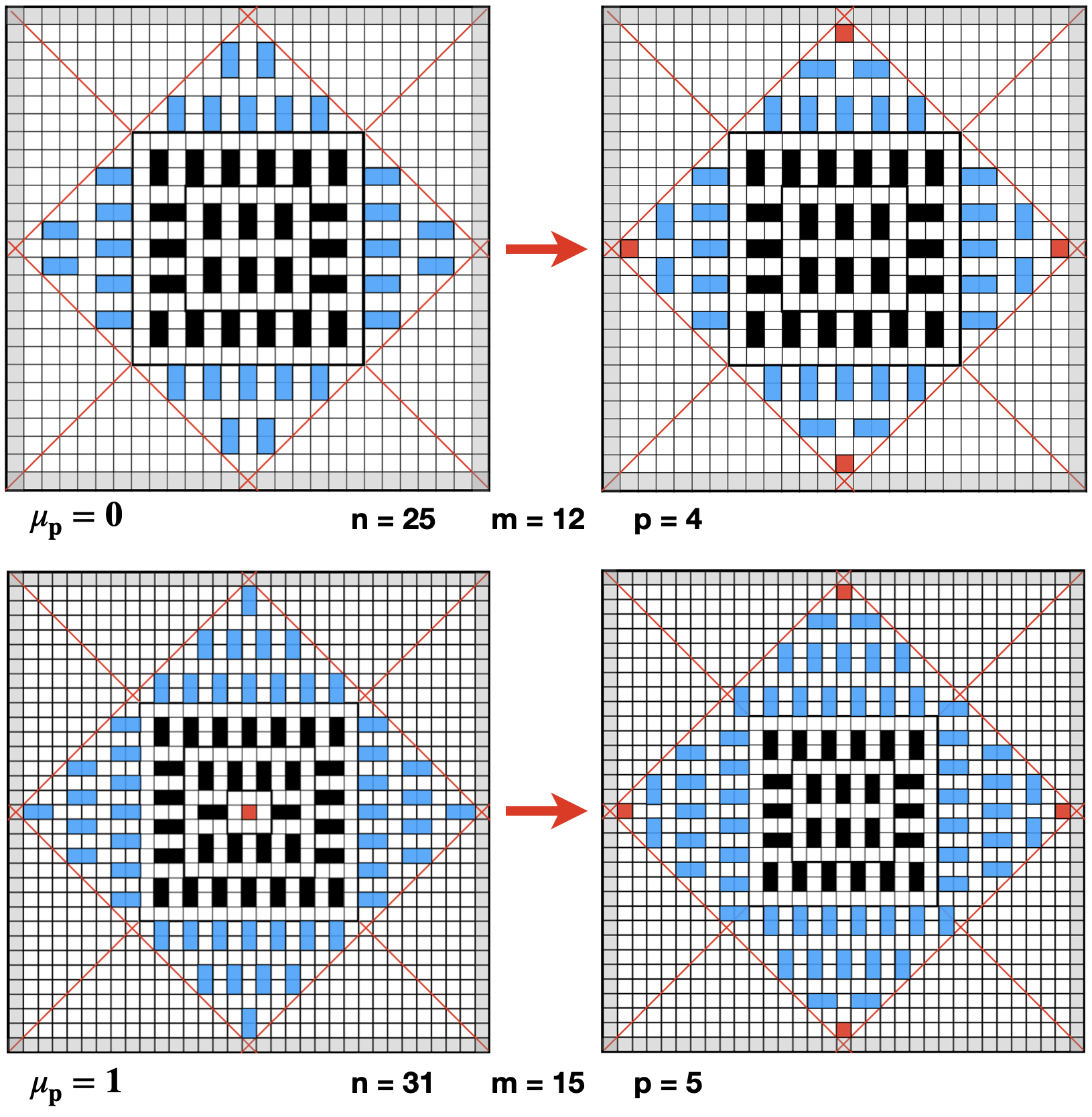}  
\caption{Injection of disorder in Class $ 10 $: this transformation yields a potential
$
               \overline{\psi}_n =  \psi_n  +  2
$
of {\em two} additional dominos whatever the parity of $ p $.
\newline
 ($ \uparrow $) $ \mu_p = 0 $ -- Ridge--flattening on the tip. 
 ($ \downarrow $) $ \mu_p = 1 $ -- Shrinking the core size.
} 
\label{Figure: IoD-C10-P4P5 }
\end{figure}
%
%
%
\begin{itemize}
%
%
\item  $ \mu_p = 1 $ --    
%
%
Fig.\,\ref{Figure: IoD-C10-P3 }($ \uparrow $)
displays a possible change by {\em enlarging} the square core from size  
 $
        f_{ 10 } ( m )   =  m - 2
 $
 to size
 $
       f_{ 10 } ( m ) + 2  = m 
 $
but this transformation is irrelevant and not suitable for optimality.
It would be easy to show from
(\ref{equation: Xi_n -- n odd -- m = 3*p+1})
and from
(\ref{equation: n odd  -- Xi-- m=3p -- mu=1 })
that
$
         \xi_{ m }  -  \xi_{ m - 2 }   =  2 p
$
but this results in the loss of one row for the wedge resulting from the reduction
$
          h_{ 10 } ( p )  \rightarrow   h_{ 10 } ( p ) - 1 
$
by turning the top domino into a lacunar void, as well as a deficit of one domino per row resulting from the reduction
$
          b_{ 10 } ( p )  \rightarrow   b_{ 10 } ( p ) - 2 
$
whence an overall deficit of $ ( p + 1 ) / 2 $ dominos per wedge, that is, a deficit of $ 2 \, p + 2 $ in all.  
We obtain finally an overall gain of
$
        - \, ( 2 \, p + 2 )  +  2 \, p   =  - 2 ,
$
i.e. an overall loss of {\em two} dominos, although compensated by the emergence of {\em four} lacunar voids.
Therefore, the overall additive potential is {\em zero} and the configuration is not better. 
%
%

%
Fig.\,\ref{Figure: IoD-C10-P3 }($ \downarrow $)
displays another change now by {\em shrinking} the square core to size  
$
        f_{ 10 } ( m ) - 2  =  m - 4 .
$
Then
$
         \xi_{ m - 4 }  -  \xi_{ m - 2 }   =  - 2 p + 2
$
from
(\ref{equation: n odd  -- DXi-- m=3p -- mu=1 })
but now
$
          h_{ 10 } ( p )  \rightarrow   h_{ 10 } ( p ) + 1  =    h_{ 10 } ( p + 1 ) 
$
and it follows that the new wedge is a copy of the wedge at $ p + 1 $. From
(\ref{equation: n odd  -- W_10 expansion -- m=3p })
this therefore gives a gain of
$
              ( p + 1 ) \, ( 1 + \mu_{ p + 1 } ) / 2   =    ( p + 1 ) / 2 
$
per wedge --first disregarding the domino conflict on the \texttt{N}--\texttt{W} junction.
By eliminating the redundant domino, the gain is then reduced to
$
            ( p + 1 ) / 2  - 1
$
that is,
$
            2 p - 2
$
for the four wedges. As a result, the overall actual gain on $ \psi_n $ is {\em zero}.
Now applying a ridge--flattening on the 2--fold tip of the wedge will release {\em four} lacunar voids in all, giving potential for {\em two} additional dominos, whence an ultimate potential capacity of $ \psi_n  +  2 $.
%
%
\item  $ \mu_p = 0 $ --    
%
%
As already mentioned in the general case for $ n $ {\em odd}, a slight transformation of ridge--flattening will yield {\em four} lacunar voids in all, giving a potential of {\em two} additional dominos, whence a potential capacity of $ \psi_n  +  2 $.
\end{itemize}
%
%
An illustration of this general transformation, whatever the parity of $ p $, is displayed in
Fig.\,\ref{Figure: IoD-C10-P4P5 }.
 This transformation holds for any $ p > 1 $ in  \texttt{Class} $ 10 $, whence
 %
%
\begin{equation}
\label{equation: n odd -- IoD -- C10 }
        \forall   m \equiv_3  0 ,  \  p > 1 , \  n  = 2 \, m + 1 :    \hspace{5mm}     
        \overline{\psi}_n =  \psi_n  +  2     
\end{equation}
%
%
where $  \overline{\psi}_n $ denotes the extended capacity resulting from this transformation.
%
%
%
\subsubsection{ In \texttt{Class} $ 11 $ }                 
%
%
%
The following two cases are considered.
%
%
%
\begin{figure}
\centering
\includegraphics[width=9cm]{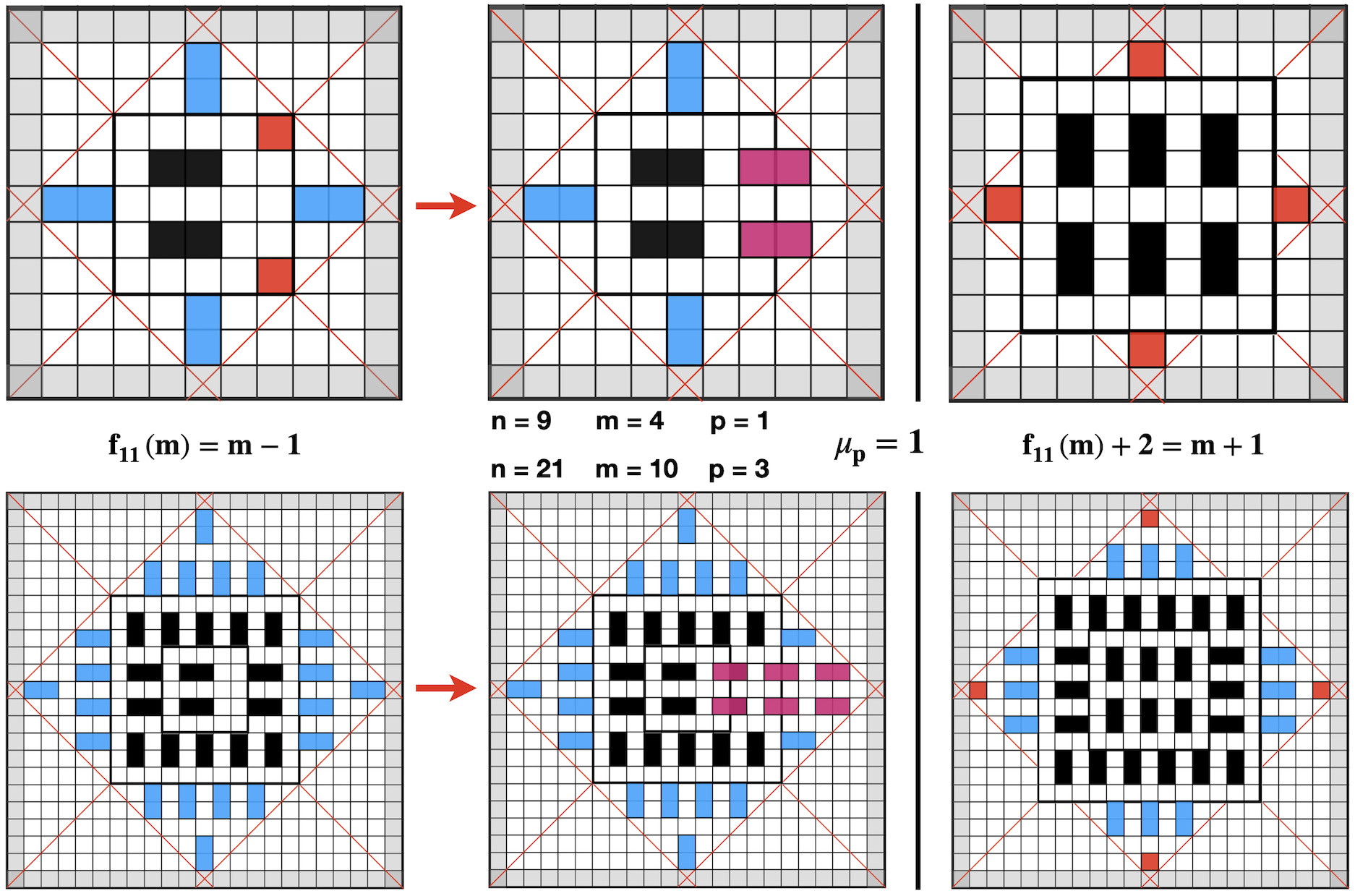}  
\caption{Injection of disorder in Class $ 11 $, case $ \mu_p = 1 $.
 ({\em Left}) Some vacant space, with or without lacunar void, can accommodate an additional domino but this anomaly cannot be propagated beyond small values of $ n $.
({\em Right}) Enlarging the square core from size  
 $
        f_{ 11 } ( m ) = m - 1   
 $
 to size
 $
       f_{ 11 } ( m ) + 2 = m + 1 .
 $
 The number of dominos remains unchanged but this transformation releases four lacunar voids on the tips.
 } 
\label{Figure: IoD-C11-P1P3 }
%
%
%
%
%
\vspace{5mm}
\centering
\includegraphics[width=7.5cm]{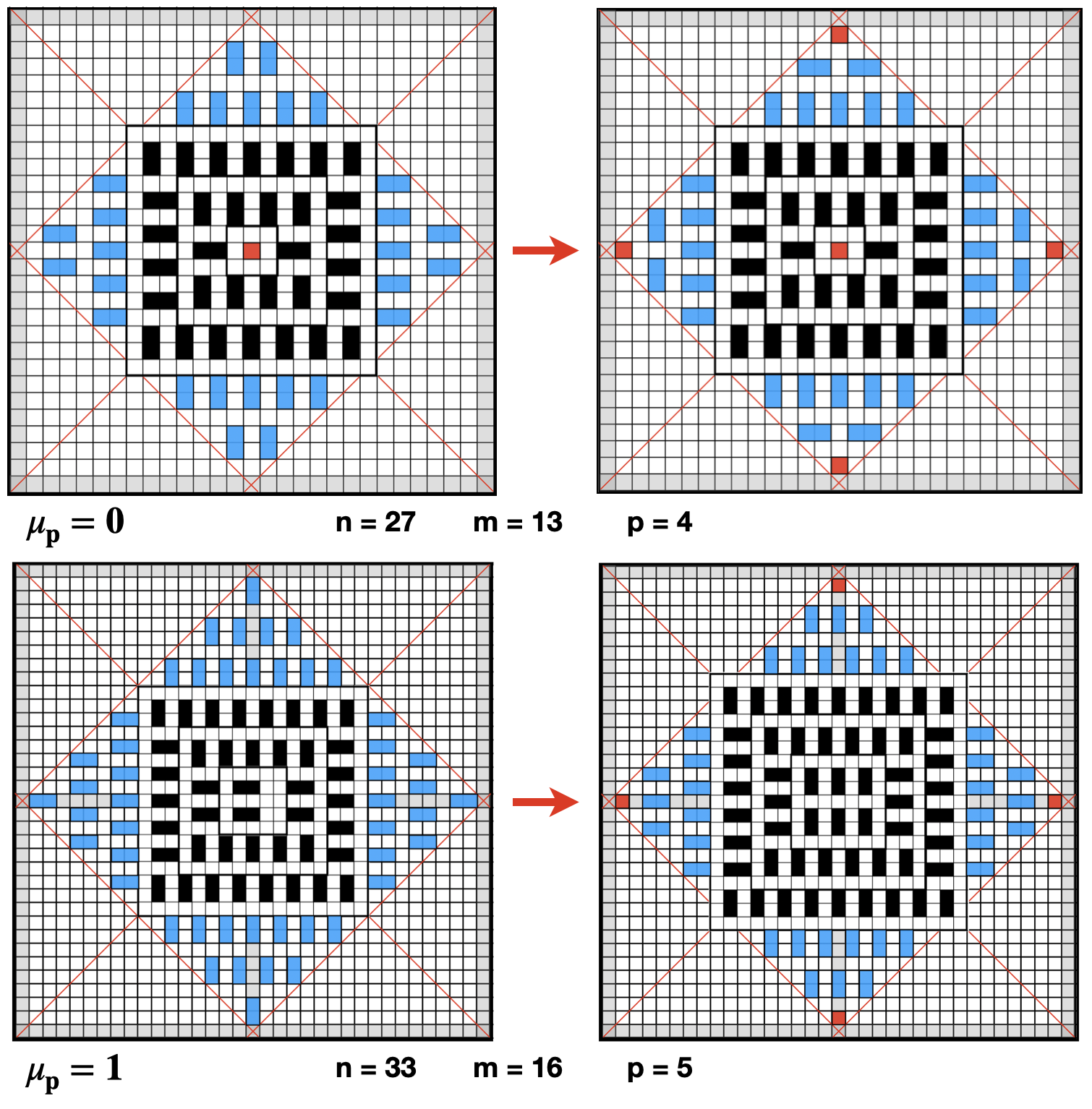}  
\caption{Injection of disorder in Class $ 11 $: this transformation yields a potential
$
               \overline{\psi}_n =  \psi_n  +  2
$
of {\em two} additional dominos whatever the parity of $ p $.
\newline
 ($ \uparrow $) $ \mu_p = 0 $ -- Ridge--flattening on the tip. 
 ($ \downarrow $) $ \mu_p = 1 $ -- Enlarging the core size.
} 
\label{Figure: IoD-C11-P4P5 }
\end{figure}
%
%
%
%
\begin{itemize}
%
%
\item  $ \mu_p = 1 $ --    
%
%
Fig.\,\ref{Figure: IoD-C11-P1P3 }
displays a possible change by breaking symmetry but this anomaly cannot be propagated beyond some low value of $ n $.
Instead, we choose another symmetric transformation by enlarging the square core from size  
 $
        f_{11} ( m )   =  m - 1
 $
 to size
 $
       f_{11} ( m ) + 2  = m + 1 
 $
 as follows.
From
(\ref{equation: Xi_n -- n odd -- m = 3*p+2})
and from
(\ref{equation: n odd  -- Xi-- m=3p+1 -- mu=1 })
the core gets the gain
$
          \xi_{ m + 1 }  -  \xi_{ m - 1 }  =   2 \, p + 2
$
but this results in the loss of one row for the wedge resulting from the reduction
$
          h_{11} ( m )  \rightarrow   h_{11} ( m ) - 1 
$
by turning the top domino into a lacunar void, as well as a deficit of one domino per row resulting from the reduction
$
          b_{11} ( m )  \rightarrow   b_{11} ( m ) - 2 
$
whence an overall deficit of $ ( p + 1 ) / 2 $ dominos per wedge, that is, a deficit of $ 2 \, p + 2 $ in total which compensates for the above core's gain.
The number of dominos thus remains unchanged but this transformation releases {\em four} additional lacunar voids, giving a potential of two additional dominos, whence an ultimate potential capacity of $ \psi_n  +  2 $.
%
%
\item  $ \mu_p = 0 $ --    
%
%
A ridge--flattening will yield {\em four} lacunar voids in all, giving a potential capacity of $ \psi_n  +  2 $.
\end{itemize}
%
%
An illustration of this general transformation, whatever the parity of $ p $, is displayed in
Fig.\,\ref{Figure: IoD-C11-P4P5 }.
This transformation holds for any $ p > 0 $ in  \texttt{Class} $ 11 $, whence
 %
%
\begin{equation}
\label{equation: n odd -- IoD -- C11 }
        \forall   m \equiv_3  1 ,  \  p > 0 , \  n  = 2 \, m + 1 :    \hspace{5mm}     
        \overline{\psi}_n =  \psi_n  +  2     
\end{equation}
%
%
where $  \overline{\psi}_n  $ denotes the extended capacity resulting from this transformation.
%
%
%
%
\begin{figure}
\centering
\includegraphics[width=11cm]{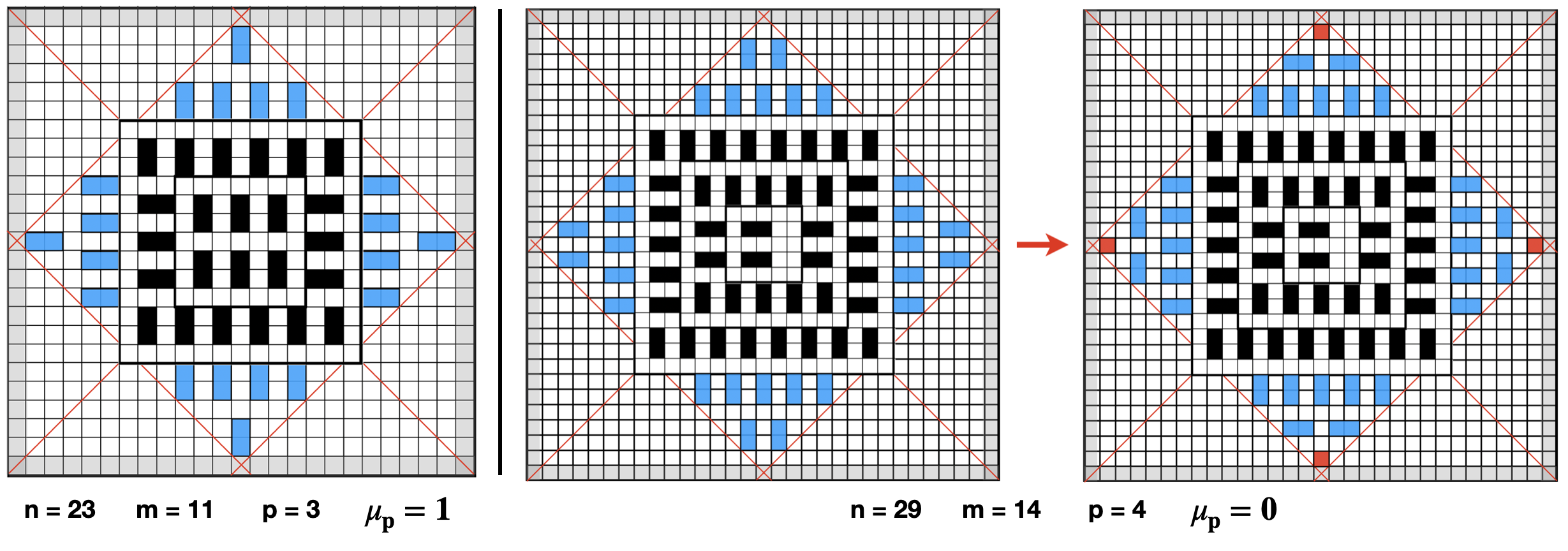}  
\caption{Injection of disorder in Class $ 12 $: this transformation yields a potential
$
               \overline{\psi}_n =  \psi_n  +  2 \,  ( 1 -  \mu_p ) .
$
\newline
 ({\em Left}) $ \mu_p = 1 $ -- The configuration is optimal.
 ({\em Right}) $ \mu_p = 0 $ -- Ridge--flattening. 
} 
\label{Figure: IoD-C12-P3P4 }
\end{figure}
%
%
%
\subsubsection{ In \texttt{Class} $ 12 $ }                 
%
%
%
According to
Fig.\,\ref{Figure: IoD-C12-P3P4 }
no change is made for $ \mu_p = 1 $ whereas for $ \mu_p = 0 $ a ridge--flattening will yield a potential 
$
               \overline{\psi}_n =  \psi_n  +  2
$
of two additional dominos.
 This transformation holds for any $ p > 0 $ in  \texttt{Class} $ 12 $, whence
 %
%
\begin{equation}
\label{equation: n odd -- IoD -- C12 }
        \forall   m \equiv_3  2 ,  \  p > 0 , \  n  = 2 \, m + 1 :    \hspace{5mm}    
         \overline{\psi}_n  =  \psi_n  +  2 \, ( 1 - \mu_{ p } )
\end{equation}
%
%
where $  \overline{\psi}_n  $ denotes the extended capacity resulting from this transformation.
For $ \mu_p = 1 $ the lower bound $  \psi_n  $ reaches the upper bound $  \overline{\psi}_n $ therefore, in this case, 
$  
          \psi_n  =     \overline{\psi}_n 
$ 
is the {\em exact} domino number for a maximal arrangement in
$
       \mathcal{D}_n .
$
%
%
\subsection{Injection of Disorder -- $n$ even}                     
%
%
A first observation of
Fig.\,\ref{Figure: C0-n6-n28 -- n even }
shows that some lacunar voids emerge as isolated cells from $ n > 20 $. 
It will even be stated that this emergence actually occurs from $ n = 20 $ in \texttt{Class} $ 01 $.
In other words, Region
$
         \mathcal{W}''_{0m}
$
is empty below this threshold, hereafter denoted 
$
       \mathbf{ n }_{ w'' }  .
$

The second observation is again related to the central pattern of the square core: a vacant space, with or without lacunar void, is likely to cause a deficit in the whole.
Referring to
Fig.\,\ref{Figure: Dominos in the Square -- n odd -- n even-- }
such a situation exists for 
$
                \mathcal{S}_{ 2 }
$
--\,condition noted briefly 
$
               ( \mbox{s}_2 )
$
\hspace{ -1mm }-- but neither for 
$
                \mathcal{S}_{ 0 }
$
nor for
$
                \mathcal{S}_{ 4 }.
$

Now, for 
$ 
          n   \ge   \mathbf{ n }_{ w'' } 
$ 
any deficiency of the core will be compensated by the potential induced by the emergence of lacunar voids.
It follows that a non--optimality problem will be induced by the conjunction  
$
               ( \mbox{s}_2 )  \land  ( n <  \mathbf{ n }_{ w'' } ).
$
It can again be released either by {\em enlarging} or by {\em shrinking} the core size, or by applying a local transformation.
One case per class is involved:
$ n = 6   $ for \texttt{Class} $ 00 $,  
$ n = 20 $ for \texttt{Class} $ 01 $,  
$ n = 16 $ for \texttt{Class} $ 02 $.
%
%
%
\subsubsection{ In \texttt{Class} $ 00 $ }                 
%
%
%
 By referring to the relation in
(\ref{equation:  n even  --  p = 4*q + r  -- C00 C01 C02 })
connecting $ p + 1 $ and $ r_{ p + 1 } $ and to the periodic sequence of expansion of $ W''_{00} $ in
(\ref{equation: n even  -- W''_{00} ( p ) - W''_{00} ( p - 1 ) -- m = 3p -- })
and beyond, we observe for any $ p > 0 $ the emergence of a new lacunar void per wedge at 
$
	r_{ p + 1 } = 1 ,
$
that is at 
$
	 p  \equiv_4  0  .
$
That gives
$
        4 \, \lfloor \, p / 4  \,  \rfloor   
$
voids in all and thus an additional potential of
$
        2 \, \lfloor \, p / 4  \,  \rfloor   
$
dominos, whence
%
%
\begin{equation}
\label{equation: n even -- IoD -- C00 }
        \forall   m \equiv_3  0  , \  n  = 2 \, m  :    \hspace{5mm}     
         \overline{\psi}_n  =  \psi_n  +  2 \, \lfloor \, p / 4  \,  \rfloor 
\end{equation}
%
%
where $   \overline{\psi}_n  $ denotes the resulting extended capacity.
The distribution of lacunar voids could already be observed in
Fig.\,\ref{Figure: C0-n6-n28 -- n even } and Fig.\,\ref{Figure: C00-P15P18 }
for  \texttt{Class} $ 00 $.

\begin{figure}
\centering
\includegraphics[width=8cm]{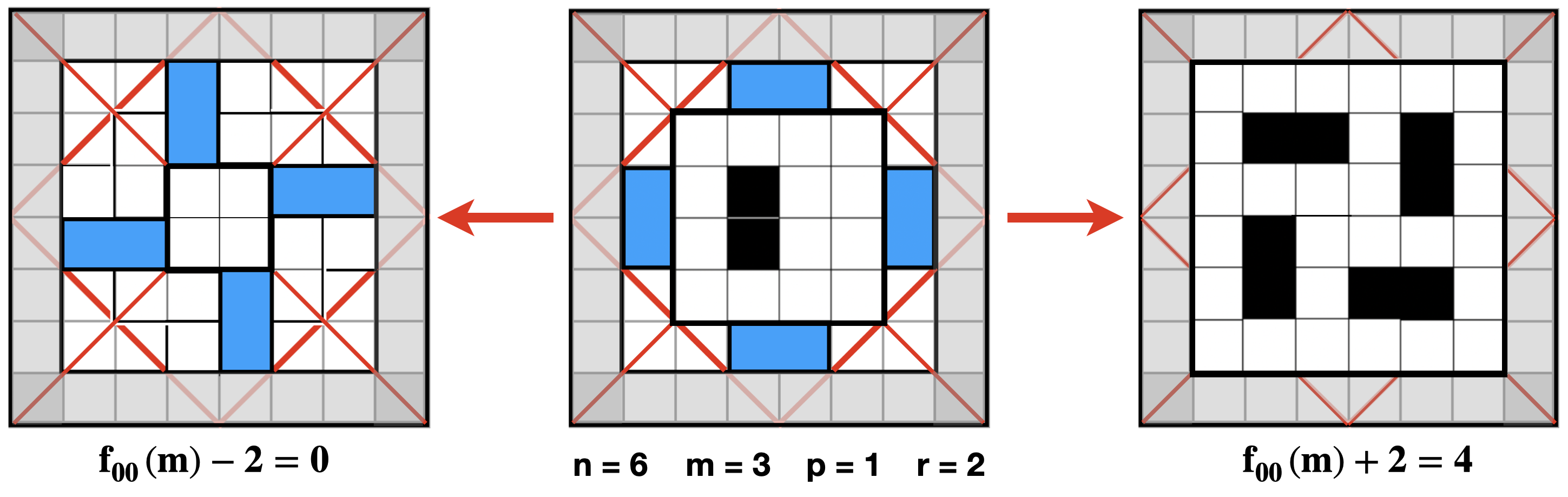}  
\caption{Injection of disorder in Class $ 00 $ for $ r_{ p + 1 } = 2 $.
Neither shrinking ($ \leftarrow $) the square core to size
 $
       f_{ 00 } ( m ) - 2   
 $
nor enlarging it ($ \rightarrow $) to size
 $
       f_{ 00 } ( m ) + 2   
 $
 gives better capacity than the basic configuration for $ n = 6 $ in 
Fig.\,\ref{Figure: C0-n6-n28 -- n even }.
} 
\label{Figure: IoD-C00-P1 }
\end{figure}

Moreover, it should be pointed out that the critical case mentioned above in this class ($ n = 6 $), while satisfying the condition 
$
               ( \mbox{s}_2 )  \land  ( n <  \mathbf{ n }_{ w'' } ),
$
necessary but not sufficient, remains optimal beyond the unsuccessful transformations of
Fig.\,\ref{Figure: IoD-C00-P1 }.
%
%
%
%
%
%
\subsubsection{ In \texttt{Class} $ 01 $ }                 
\begin{figure}
\centering
\includegraphics[width=12cm]{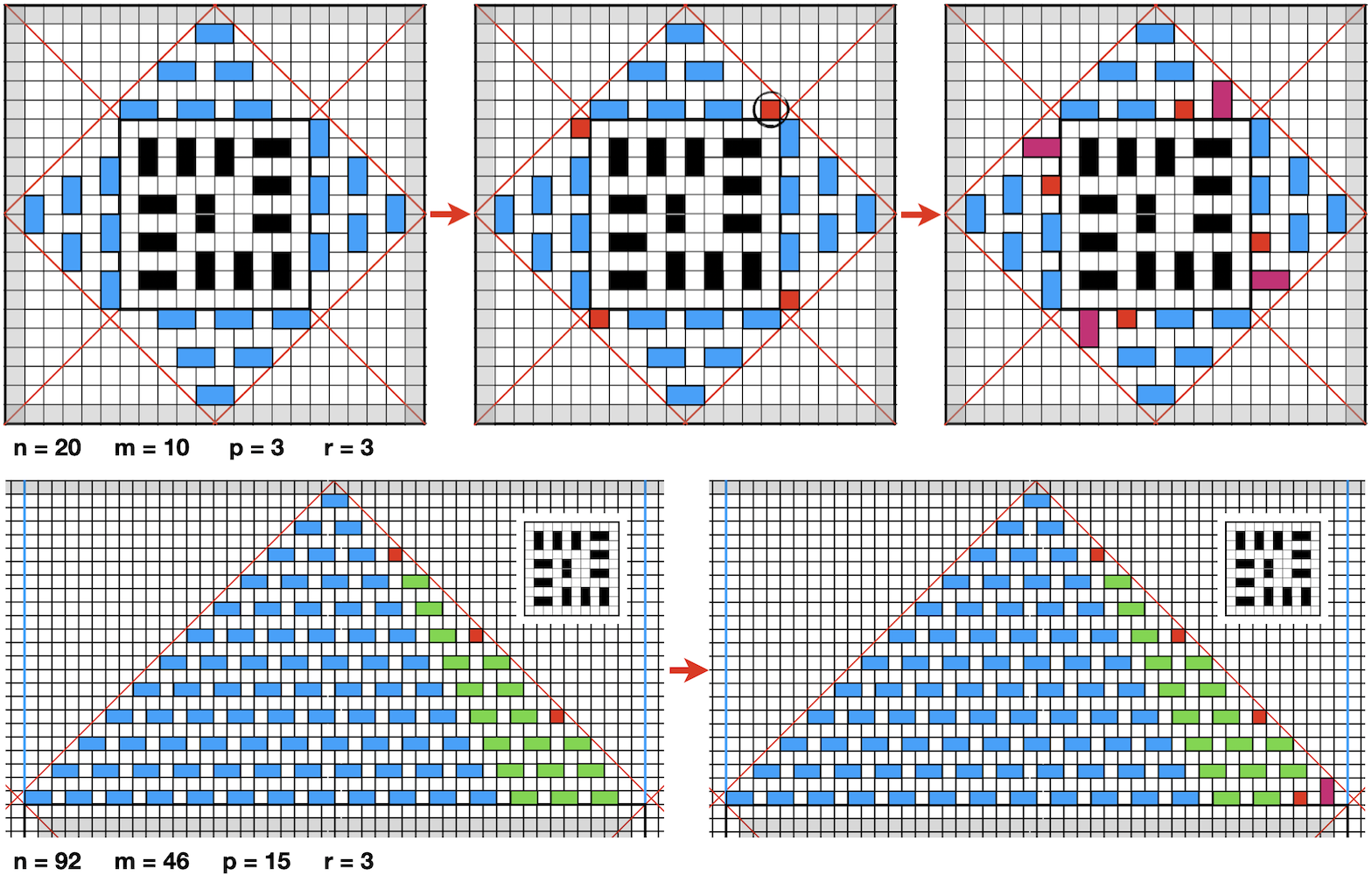}  
\caption{Injection of disorder in Class $ 01 $ -- 
Implicit void in the baserow of the wedge for $ r_{ p } = 3 $, taking an explicit form after a slight transformation.
} 
\label{Figure: IoD-C01P3P15 }
\end{figure}

%
%
%
By referring to the relation in
(\ref{equation:  n even  --  p = 4*q + r  -- C00 C01 C02 })
connecting $ p  $ and $ r_{ p } $ and to the periodic sequence of expansion of $ W''_{01} $ in
(\ref{equation: n even  -- W''_{01} ( p ) - W''_{01} ( p - 1 ) -- m = 3p+1 -- }), 
we observe for any $ p > 0 $ the emergence of a new lacunar void at 
$
	r_{ p } = 0 .
$
However, we also observe that there is a space at 
$
	r_{ p } = 3
$
to insert a new lacunar void but, by symmetry, this void remains closed by the baserow of the  \texttt{East} wedge.
In fact, everything happens as if this void nevertheless had the property of having to be taken into account.
This fact finally becomes true, hence at
$
	p + 1  \equiv_4  0  ,
$
by applying the slight transformation illustrated in 
Fig.\,\ref{Figure: IoD-C01P3P15 }.
That gives
$
        4 \, \lfloor \, ( p + 1 ) / 4  \,  \rfloor   
$
voids in all and thus an additional potential of
$
        2 \, \lfloor \, ( p + 1 ) / 4  \,  \rfloor   
$
dominos, whence
%
%
\begin{equation}
\label{equation: n even -- IoD -- C01 }
        \forall   m \equiv_3  1  , \  n  = 2 \, m  :    \hspace{5mm}    
         \overline{\psi}_n  =  \psi_n  +  2 \, \lfloor \, ( p + 1 ) / 4  \,  \rfloor 
\end{equation}
%
%
where $  \overline{\psi}_n $ denotes the extended capacity resulting from such a transformation.
The distribution of lacunar voids could already be observed in
Fig.\,\ref{Figure: C0-n6-n28 -- n even } and Fig.\,\ref{Figure: C01-P16P19 }
for  \texttt{Class} $ 01 $.

Note that the previously mentioned critical case in this class ($ n = 20 $) coincide with the threshold 
$
       \mathbf{ n }_{ w'' }  
$
of lacunar emergence.
%
%
%
%
%
%
\subsubsection{ In \texttt{Class} $ 02 $ }                 
%
%
%
%
\begin{figure}
\centering
\includegraphics[width=12cm]{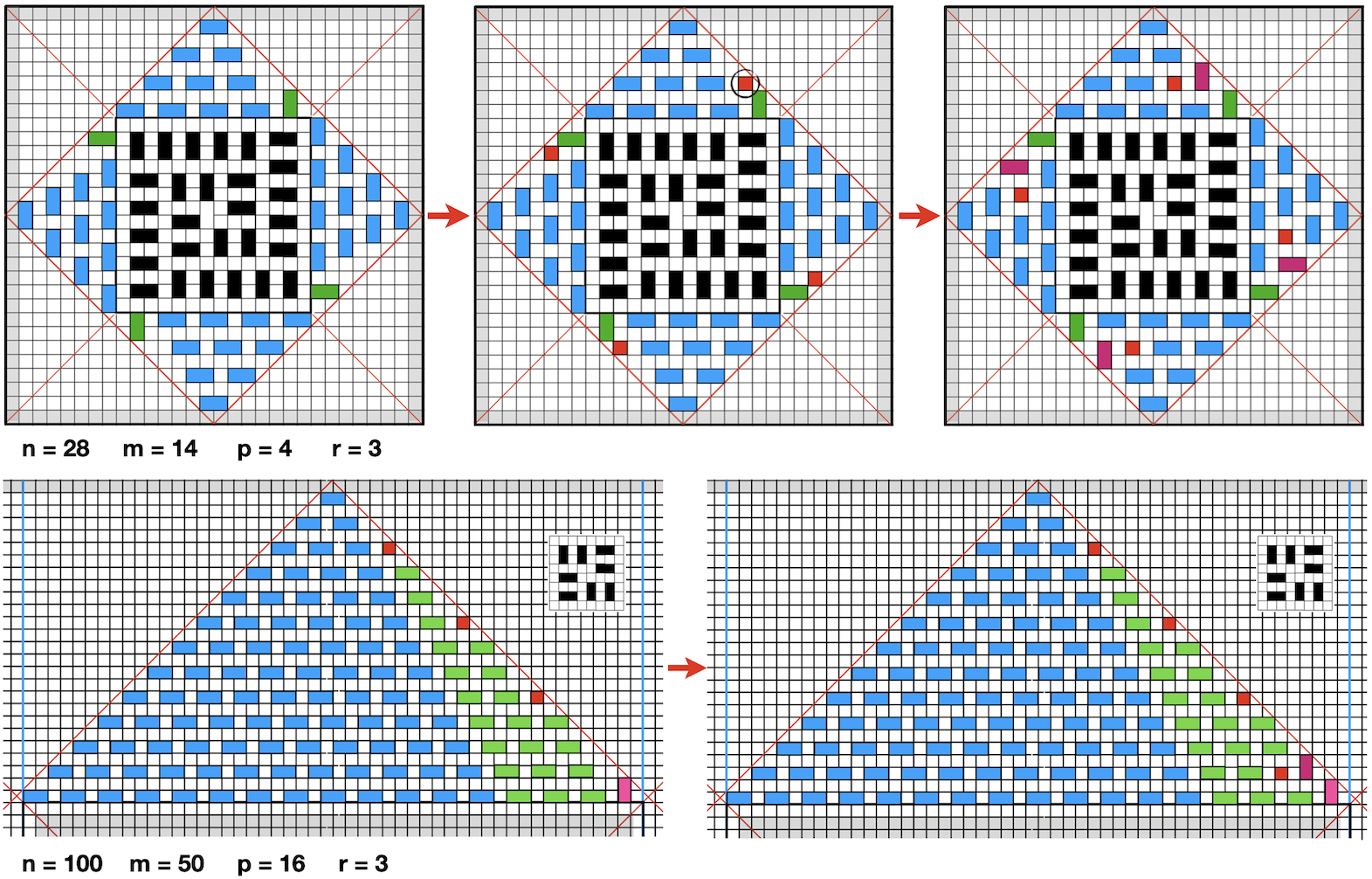}  
\caption{Injection of disorder in Class $ 02 $ for $ r_{ p - 1 } = 3 $ -- 
Implicit void in the penultimate row of the wedge, taking an explicit form after a slight transformation.
} 
\label{Figure: IoD-C02-P4P16--- r = 3 }
\end{figure}
%
%
This class is the most intricate and several cases should be examined. 
We refer again to the corresponding relation in
(\ref{equation:  n even  --  p = 4*q + r  -- C00 C01 C02 })
connecting $ p - 1 $ and $ r_{ p - 1 } $ and to the periodic sequence of expansion of $ W''_{02} $ in
(\ref{equation: n even  -- W''_{02} ( p ) - W''_{02} ( p - 1 ) -- m = 3p+2 -- }).
Without loss of generality, suppose first that 
%
%
\begin{equation}
\label{equation: n even -- IoD -- C02 -?- }
        \forall   m \equiv_3  2  , \  n  = 2 \, m  :    \hspace{5mm}     \psi_n^+  =  \psi_n  +  2 \, \lfloor \, ( p + 2 ) / 4  \,  \rfloor 
\end{equation}
%
%
where $  \psi_n^+  $ denotes the extended capacity that {\em would} result from the following counting and transformations.
%
%
%
\begin{itemize}
%
%
\item  $ r_{ p - 1 } = 3  $  \ --- 
%
%
The new lacunar void, emerging in the baserow of the previous step at $ r_{ p - 1 } = 2 $, disappears under the effect of the vertical rotation of the new domino inserted in the new baserow
(see Fig.\,\ref{Figure: C0-n6-n28 -- n even } and Fig.\,\ref{Figure: C02-P17P20 }).
In fact, everything happens as if this void nevertheless had the property of having to be taken into account.
This fact finally becomes true by applying the slight transformation illustrated in 
Fig.\,\ref{Figure: IoD-C02-P4P16--- r = 3 }
then the potential capacity $ \psi_n^+ $ assumed in 
(\ref{equation: n even -- IoD -- C02 -?- })
remains valid.
%
%
\begin{figure}
\centering
\includegraphics[width=8cm]{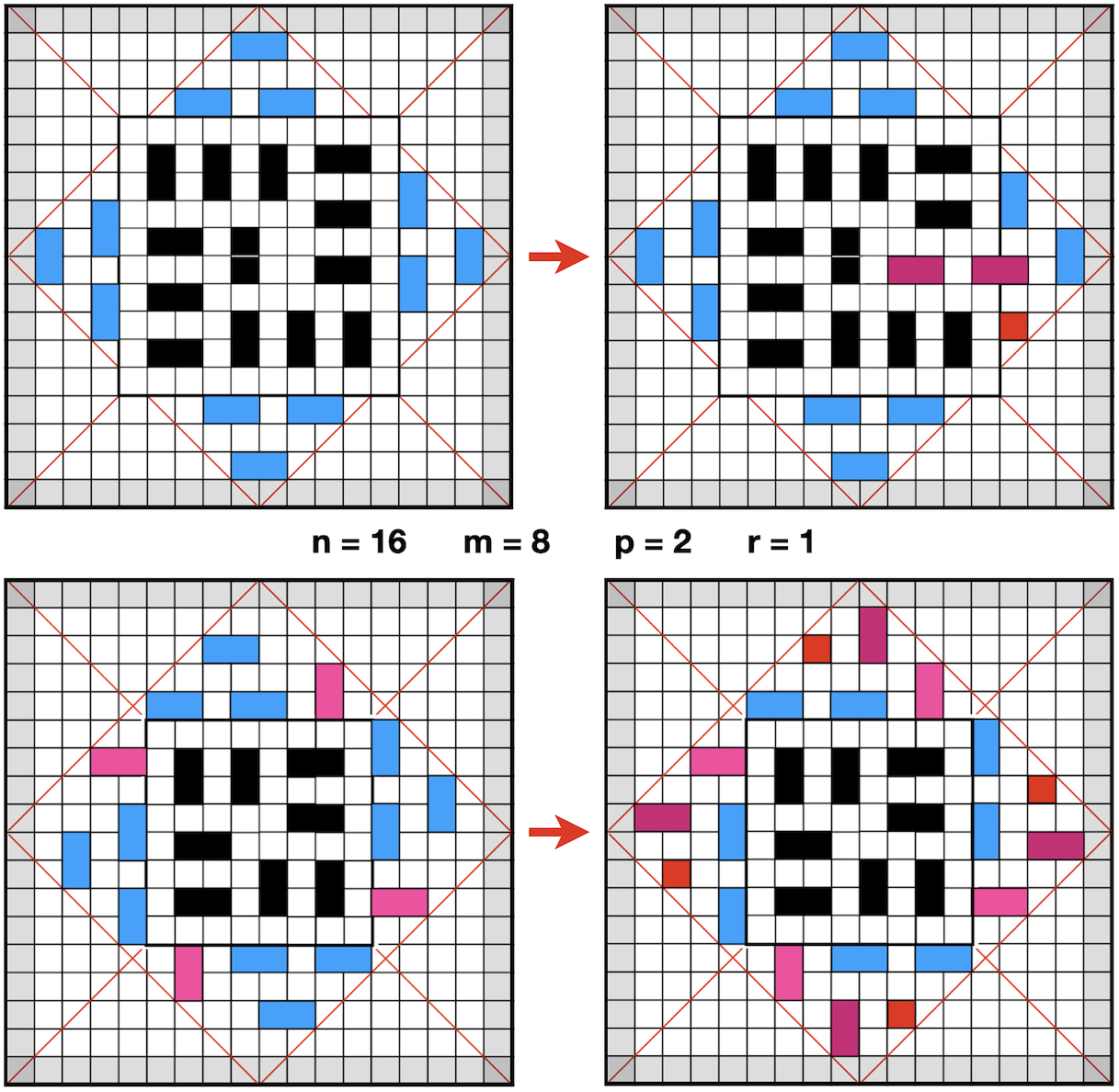}  
\caption{Injection of disorder in Class $ 02 $ for $ r_{ p - 1 } = 1 $ --
 ($ \uparrow $) Releasing a lacunar void from some vacancy in the core.
 ($ \downarrow $) Shrinking the square core from size  $ f_{02} ( m ) = m $ to size $ f_{02} ( m ) - 2 $ with reduction  
 $         
         \xi_{ m - 2 }  -  \xi_{ m }  =   - \, ( 2 \, p + 1 ) = - 5
$
for this case; overall gain of $ 2 p = 4 $ for the wedges; overall loss of {\em one} domino altogether; overall gain of {\em four} lacunar voids giving a potential of {\em two} additional dominos. 
Initial configuration with $ \psi_n = 25 $, final configuration with $ \psi_n - 1 = 24 $, potential capacity $  \overline{\psi}_n = 26 $.
} 
\label{Figure: IoD-C02-P2 }
\end{figure}
%
%
%
\begin{figure}
\centering
\includegraphics[width=12.4cm]{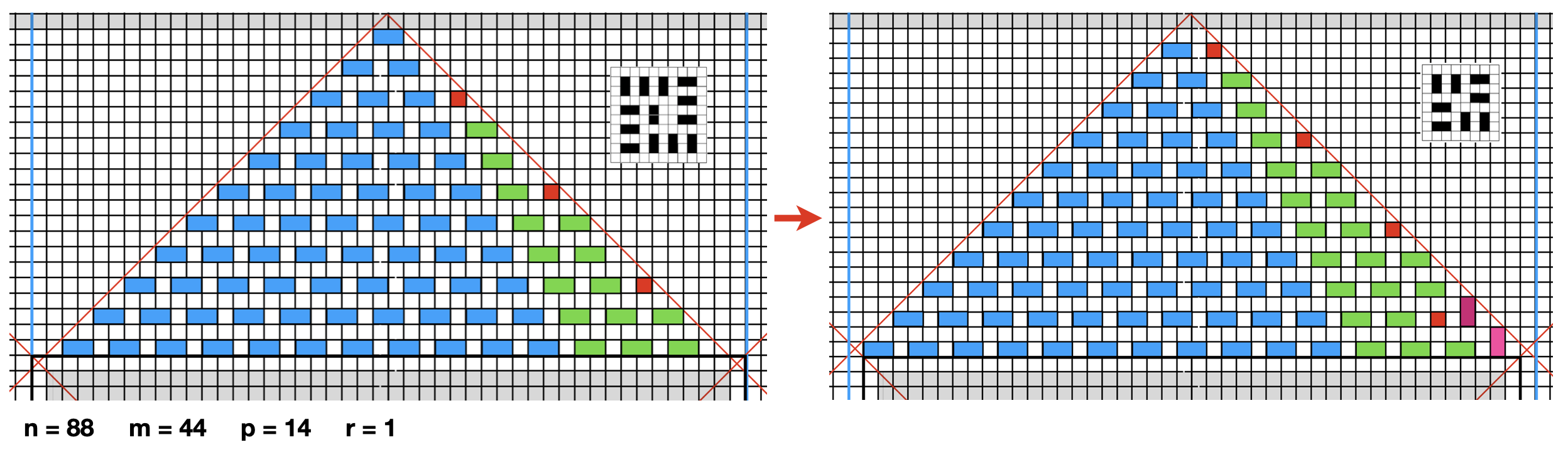}  
\caption{Injection of disorder in Class $ 02 $ for $ r_{ p - 1 } = 1 $ --
Shrinking the square core. After shrinking, the resulting configuration of Region 
$
         \mathcal{W}''_{02}
$
becomes a copy of itself after the same slight transformation for $ r_{ p - 1 } = 3 $ in 
Fig.\,\ref{Figure: IoD-C02-P4P16--- r = 3 }.
} 
\label{Figure: IoD-C02P14 -- r=1 }
\end{figure}
%
%
%
%
\item  $ r_{ p - 1 } = 1  $  \ --- 
%
%
In the center of the core, we can observe some vacancy which is likely to allow the emergence of voids as shown in
Fig.\,\ref{Figure: IoD-C02-P2 }.
We therefore choose a symmetric transformation by shrinking the square core from size  
 $
        f_{02} ( m )   =   m
 $
(from (\ref{equation:  n even  --  f_{ 0 m }( m ) -- C00 C01 C02 }))
 to size
 $
       f_{02} ( m ) - 2 .
 $
Then 
$         
         \xi_{ m - 2 }  -  \xi_{ m }  =   - \, ( 2 \, p + 1 )
$
from
(\ref{equation: n even  -- DXi-- m=3p+2 -- r = 1 })
but this results in the gain
$
          h_{02} ( m )  \rightarrow   h_{02} ( m ) + 1 
$
of one row for the wedge as well as the expansion
$
          b_{02} ( m )  \rightarrow   b_{02} ( m ) + 2 
$
of two cells for the new baserow.
As a consequence, the block of Region
$
         \mathcal{W}'_{02}
$
is shifted both downward and leftward while remaining adjacent to the shrunk core whereas capacity  $ W'_{02} $ remains unchanged.
 
 On the other hand, it follows that the configuration of Region
$
         \mathcal{W}''_{02}
$
turns into a copy of that at $ r_{ p - 1 } = 3  $ within the $ q_{ p - 1 } $--cycle (at constant $ q_{ p - 1 } $).
Thus by achieving the suitable assessments from 
(\ref{equation: n even  -- W''_{02} ( p ) - W''_{02} ( p - 1 ) -- m = 3p+2 -- })
we obtain the gain  
\[
         	W''_{02} ( p + 2 )  -   W''_{02} ( p )   =   
	        ( W''_{02} ( p + 2 )  -   W''_{02} ( p + 1 ) )  +  ( W''_{02} ( p + 1 )  -   W''_{02} ( p ) )		 
\]
namely
\[
	        (   q_{ p - 1 }  - 1  + \delta_{ 3 } )  +   (   q_{ p - 1 }  - 1  + \delta_{ 2 } )	 
\]
and since $ \delta_{ 3 } =  2 $ and $ \delta_{ 2 } =  1 $ from
(\ref{equation: n even  --  capacity of b"_(02) -- m = 3p+2 -- })
we obtain the gain
$        
                   2 \, q_{ p - 1 } + 1  =  p / 2
$
per wedge and therefore an overall gain of $ 2 \, p $ for the four wedges. 
Finally the overall gain becomes
$
        - \, ( 2 \, p + 1 )  +  2 \, p   =  - 1,
$
i.e. an overall loss of {\em one} domino.
\end{itemize}
%
%
Regarding Region
$
         \mathcal{W}''_{02}
$
it is worth comparing the last configuration in
Figs.\,\ref{Figure: IoD-C02-P2 } \& \ref{Figure: IoD-C02-P4P16--- r = 3 }
(resp. for ($ n = 16, r_{ p - 1 } = 1 $) \& ($ n = 28, r_{ p - 1 } = 3 $))
and the last configuration in 
Figs.\,\ref{Figure: IoD-C02P14 -- r=1 } \& \ref{Figure: IoD-C02-P4P16--- r = 3 }
(resp. for ($ n = 88, r_{ p - 1 } = 1 $) \& ($ n = 100, r_{ p - 1 } = 3 $)).

Finally, the extended capacity $  \psi_n^+ $ in
(\ref{equation: n even -- IoD -- C02 -?- })
holds almost everywhere, except for the loss of {\em one} domino for $ r_{ p - 1 } = 1 $ whence the final expression  
%
%
\begin{equation}
\label{equation: n even -- IoD -- C02 }
        \forall   m \equiv_3  2  , \  n  = 2 \, m  :    \hspace{5mm}     
        \overline{\psi}_n =  \psi_n  +  2 \, \lfloor \, ( p + 2 ) / 4  \,  \rfloor   -  \lambda_{ r_{ p - 1 } }
\end{equation}
%
%
where
$
          \lambda_{ r_{ p - 1 } }  =  1
$
for $ r_{ p - 1 } = 1 $ and 
$
          \lambda_{ r_{ p - 1 } }  =  0
$
otherwise
and where $  \overline{\psi}_n $ denotes the resulting extended capacity in 
\texttt{Class} $ 02 $.
%
%
%
\section{Discussion and Future Issues}
\label{Section: Discussion and Future Issues}
%
%
The subject herein was to find a maximal arrangement of {\em dominos} in the {\em diamond}
$
           \mathcal{D}_{ n } 
$
of size $ n $.
The words ``domino" and ``diamond" are understood according to their exact meaning defined in 
Sect.\,\ref{Section: General Statements}.
This construction was carried out in two stages.
First, a deterministic arrangement was proposed which led to a sub--optimal solution $ \psi_n $ as the {\em lower} bound.
Second, some disorder has been injected, leading to an {\em upper} bound $  \overline{\psi}_n  $ reachable or not.
Despite our inability to achieve the optimal $ \psi_{ n }^\ast $ for any $n$, we have nevertheless made a significant advance.
The main results are presented hereafter.
%
%
%
\paragraph{Linear evolution of the ratio $  \overline{\psi}_n  / n $.}                
%
%
%
\begin{figure}
\centering
\includegraphics[width=11cm]{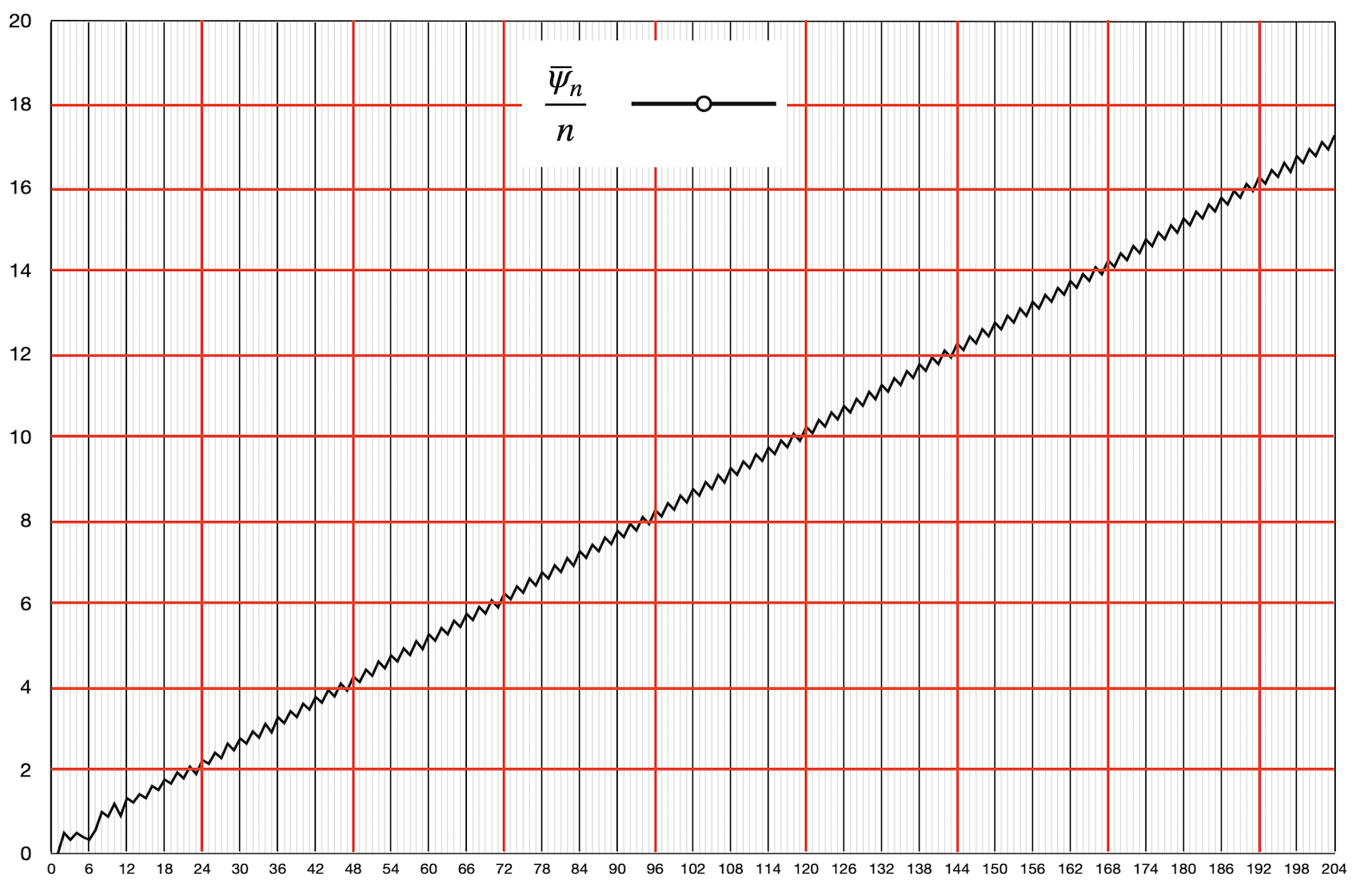}  
\caption{Linear evolution on average of the ratio $  \overline{\psi}_n  / n $  \textsf{vs} $ n $ after a chaotic behavior for small values of $ n $.
The graph is a discrete {\em dust} of points in 
$
         \mathbb{N} \times \mathbb{Q}
$
with jump discontinuities, represented as a sawtooth continuous line.
} 
\label{Figure: Conclusion -- PsiB_sur_n --  }
\end{figure}
%
%
It is not a breakthrough to observe in
Fig.\,\ref{Figure: Conclusion -- PsiB_sur_n --  }
that the slope of 
$
          { \overline{\psi}_n } /  { n } 
$
({\em vs. }$ n $) is linear on average after a chaotic behavior for small values of $ n $.
At this scale, we are allowed to disregard the negligible deviation
$
              \overline{\psi}_n   -    \psi_{ n } 
$
as it will be examined thereafter.
As a first approximation, the linearity is evidenced for the odd $ n $ in
(\ref{equation: n odd -- m == 0  -- Psi(n) - Psi( n - 6 ) },
\ref{equation: n odd -- m == 1  -- Psi(n) - Psi( n - 6 ) },
\ref{equation: n odd -- m == 2  -- Psi(n) - Psi( n - 6 ) })
by the linear rate of change
$
           \psi_{ n }  -  \psi_{ n - 6 }   
$
following a
$
          p  \rightarrow  p + 1
$
cycle within any  \texttt{Class}   $ 1 m $.
 
For the even $ n $, in spite of the fluctuations observed in 
Tabs.\,\ref{Table: Expansion Psi -- m = 3p }--\ref{Table: Expansion Psi -- m = 3p+2 }
a linear rate of change
$
           \psi_{ n }  -  \psi_{ n - 24 }  
$
is nevertheless observed in
(\ref{equation: n even -- m == 0  -- Psi(n) - Psi( n - 24 ) },
\ref{equation: n even -- m == 1  -- Psi(n) - Psi( n - 24 ) },
\ref{equation: n even -- m == 2  -- Psi(n) - Psi( n - 24 ) })
at a larger scale following a
$
          p  \rightarrow  p + 4
$
cycle within any  \texttt{Class}   $ 0 m $. 

\begin{figure}
\centering
\includegraphics[width=11cm]{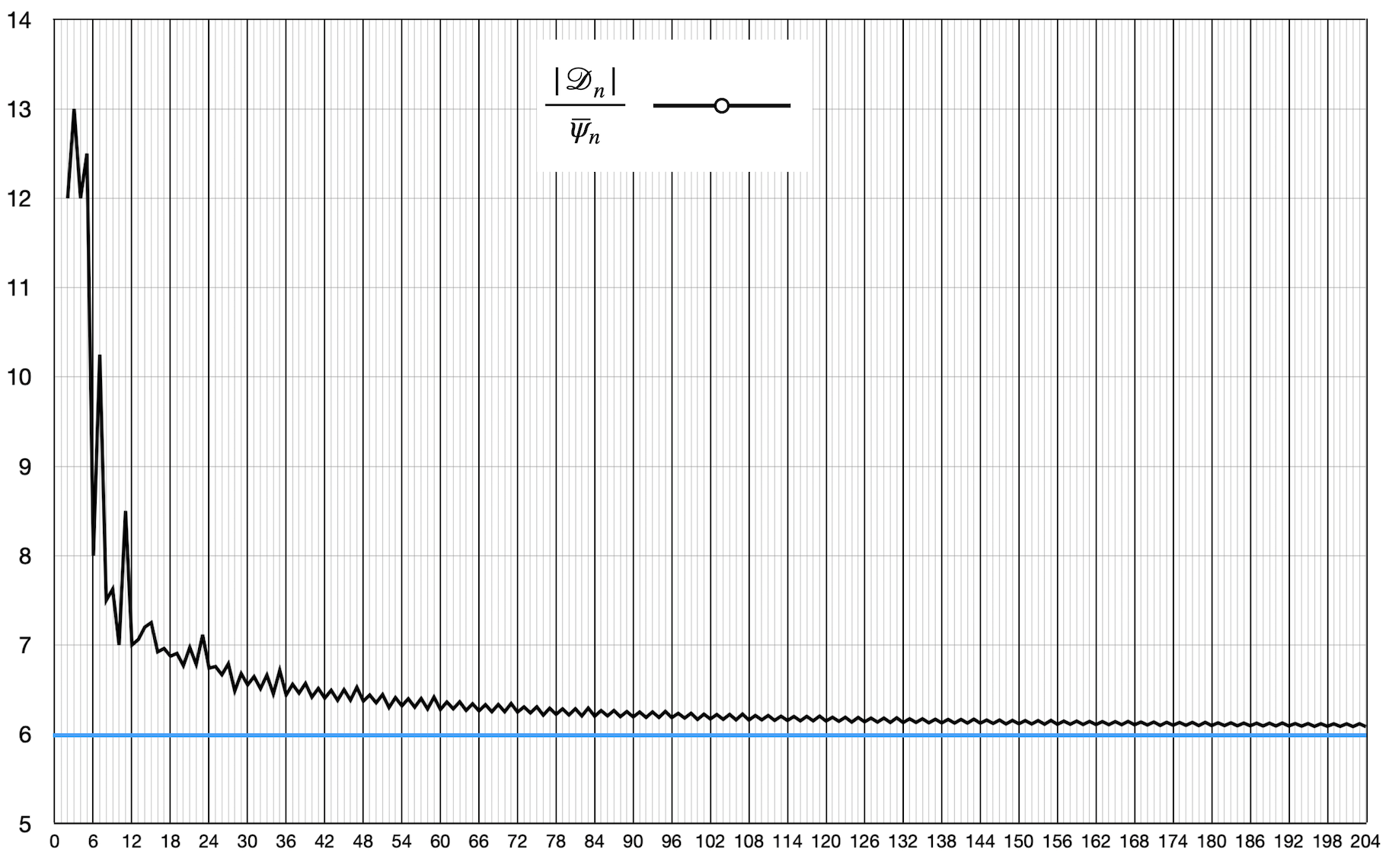}  
\caption{Beyond a chaotic behavior for small  $ n $, evolution of the global occupancy 
$
           | \mathcal{D}_{ n } | / { \overline{\psi}_n } 
$
 \textsf{vs} $ n $ towards the optimum --\,the blue line.
The graph is a sawtooth discrete {\em dust} of points in 
$
         \mathbb{N} \times \mathbb{Q}
$
with jump discontinuities.
} 
\label{Figure: Conclusion -- Dn-sur PsiB --  }
\end{figure}

The graph of
$
          { \overline{\psi}_n } /  { n } 
$
is a sawtooth curve (the ``curve'' should be seen here as a discrete {\em dust} of points in 
$
         \mathbb{N} \times \mathbb{Q}) 
$
with jump discontinuities.
This phenomenon is explained by the ``von Neumann $ \leftrightarrow $ Aztec'' expansion
$
            |  \mathcal{D}_{ n } |  -  | \mathcal{D}_{ n - 1 } | 
$
in (\ref{equation:  | D_{ n } |  --  | D_{ n - 1 } | -- })
where the diamond grows significantly only one time out of two.
%
%
%
\paragraph{Evolution of the global occupancy towards the optimum.}                
%
%
%
Fig.\,\ref{Figure: Conclusion -- Dn-sur PsiB --  }
shows the evolution with $ n $ of the global occupancy
$
           | \mathcal{D}_{ n } | / { \overline{\psi}_n } 
$
--\,as defined in
Subsect.\,\ref{Subsection: Density Measurement-- }\,--
which decreases asymptotically towards its optimal limit.
The graph is again a sawtooth curve (as discrete {\em dust} of points in 
$
         \mathbb{N} \times \mathbb{Q}) 
$
with jump discontinuities.
Again we are allowed to disregard the deviation
$
              \overline{\psi}_n   -    \psi_{ n } 
$
negligible at this scale.
We examine the following two cases
%
%
%
\begin{itemize}
%
%
\item  in  \texttt{Class}   $ 10 $ ---
From the cardinality of $ \mathcal{D}_{ n } $ in
(\ref{equation:  | D_{ n } |  -- n odd -- })
and the direct expression of $ \psi_{ n } $ in
(\ref{equation: n odd  -- PSi in C10 })
it comes
\[
             \lim\limits_{n \rightarrow +\infty}   \frac{\mid \mathcal{D}_{n}\mid} { \psi_{n}}  =  
              \lim\limits_{n \rightarrow +\infty}   \frac{ ( n^2 + 4 n + 5 ) / 2 } { ( n^2 - 1 ) / 12 }  = 
             \lim\limits_{n \rightarrow +\infty}  6 \cdot ( 1 +  \frac{ 4 } { n } \, ) 
\]
and the global occupancy tends to the optimum --\,the blue asymptote\,-- as $ n $ tends to infinity.
The same result is obtained in  \texttt{Class} $ 11 $ from 
(\ref{equation: n odd  -- PSi in C11 })
and in  \texttt{Class} $ 12 $ from 
(\ref{equation: n odd  -- PSi in C12 }).
%
%
%
\item  in  \texttt{Class}   $ 00 $ ---
For lack of a simple, direct expression of $ \psi_{ n } $ for this class, we assume $ p > 4 $ and refer to the linear rate of change
$
           \psi_{ n }  -  \psi_{ n - 24 }  
$
in
(\ref{equation: n even -- m == 0  -- Psi(n) - Psi( n - 24 ) })
following a
$
          p  \rightarrow  p + 4
$
cycle.
Now, from the cardinality of $ \mathcal{D}_{ n } $ in
(\ref{equation:  | D_{ n } |  -- n even -- })
it comes
\[
      \mid \mathcal{D}_{ n - 24 } \mid  =    ( ( n - 24 )^2  +  6 \,( n - 24 )  +  8 ) / 2   =
        \,   \mid \mathcal{D}_{ n }\mid   - \  24 \,  ( n - 9 )   
\] 
whereas 
$
           \psi_{ n - 24 }  =  \psi_{ n }  -  4 \, ( n - 11 )
$ 
from
(\ref{equation: n even -- m == 0  -- Psi(n) - Psi( n - 24 ) }).
Therefore, if the ratio
$
           | \mathcal{D}_{ n } | /  \psi_{ n }
$
admits a limit when $ n $ goes to infinity, then it must follow a Cauchy sequence. Now
\[
       \frac{\mid \mathcal{D}_{ n }\mid} { \psi_{ n }}     \sim      \frac{\mid \mathcal{D}_{ n - 24 }\mid} { \psi_{ n - 24 }}  =  
       \frac{\mid \mathcal{D}_{n}\mid - \, 24  \, ( n - 9 ) } { \psi_{n}  -  4 \,( n - 11 ) } 
\]
whence
\[
        \mid \mathcal{D}_{n}\mid  ( \psi_{ n }  -  4 \, ( n - 11 ) )     \sim    \psi_{n} \, ( \mid \mathcal{D}_{n}\mid   -  \, 24  \, ( n - 9 ) )
\]
\[
        \Rightarrow     4 \,( n - 11 )  \mid \mathcal{D}_{n}\mid    \  \sim    24  \, ( n - 9 )  \, \psi_{n}
\]
and finally
\[
             \lim\limits_{n \rightarrow +\infty}   \frac{\mid \mathcal{D}_{n}\mid} { \psi_{n}}  = 
             \lim\limits_{n \rightarrow +\infty}  6 \cdot  \frac{ n - 9 } { n - 11 }  =  
             \lim\limits_{n \rightarrow +\infty}  6 \cdot ( 1 +  \frac{ 2 } { n } \, ) 
\]
and the global occupancy tends to the optimum.
The same result is obtained in  \texttt{Class} $ 01 $ from 
(\ref{equation: n even -- m == 1  -- Psi(n) - Psi( n - 24 ) })
and in  \texttt{Class} $ 02 $ from 
(\ref{equation: n even -- m == 2  -- Psi(n) - Psi( n - 24 ) }).
\end{itemize}
%
%
%
      %
%
%
%
%
%
%
%
\paragraph{Absolute and relative deviations between {\em lower} and {\em upper} bounds.}                
%
%
%
We now focus on the deviation
$
             \Delta  \psi_{ n } =  \overline{\psi}_n  -  \psi_{ n }  \ge 0
$
between the suboptimal {\em lower} bound $  \psi_{ n } $ and the {\em upper} bound $  \overline{\psi}_n $ reachable or not, produced after injection of disorder.
The two following cases depend on the parity of $ n $.
%
%
\begin{itemize}
%
%
\item  $ n $ {\em odd} ---
By grouping all the cases resulting from 
(\ref{equation: n odd -- IoD -- C10 }--\ref{equation: n odd -- IoD -- C12 })
it comes
\[          \overline{\psi}_n - \psi_{ n }  =   \left\{
                        \begin{array}{ll}
                                  0	& 	\mbox{ in \texttt{Class} $ 12 $ \  if  $ p $  odd }	\\
                                  2 	& 	\mbox{ otherwise }		 
                        \end{array}
                                                                                                \right.
\]
namely
$
             \overline{\psi}_n - \psi_{ n }  = 2
$
in general, except in one case where the {\em exact} value is obtained.
Moreover, from the conditions on $ p $ in
(\ref{equation: n odd -- IoD -- C10 }--\ref{equation: n odd -- IoD -- C12 })
the exact value is also checked for small $ p < 2 $ ($ n < 12 $) except for $ n = 9 $.
The relative deviation
$
           \Delta  \psi_{ n } / \psi_{n} 
$
in the general case, is only relevant beyond the chaotic domain highlighted in 
Fig.\,\ref{Figure: Conclusion -- Dn-sur PsiB --  },
namely for $ n \ge 24 $ and can be expressed as follows.
For $ n $ large, we can consider from
(\ref{equation: n odd  -- PSi in C10 }, \ref{equation: n odd  -- PSi in C11 }, \ref{equation: n odd  -- PSi in C12 })
that
$
        \psi_{n}  \sim  n^2 / 12
$
and
$
           \Delta  \psi_{ n } / \psi_{n}  \sim  24 / n^2
$
thereof.
There is a near coincidence of the lower and upper curves.
%
%
\item  $ n $ {\em even} ---
Without going into details, we can estimate from
(\ref{equation: n even -- IoD -- C00 },\,\ref{equation: n even -- IoD -- C01 },\,\ref{equation: n even -- IoD -- C02 })
that the absolute deviation behaves like
\[
                  \overline{\psi}_n - \psi_{ n }  \approx   2 \, \Bigl\lfloor \,  \frac{ p + m  \bmod 3 }{ 4 }  \,  \Bigr\rfloor 
\]
approximately. As a result, we get a collection of exact values of $  \psi_{n} $ for almost small $ p < 4 $ ($ n < 24 $).
For the relative deviation with  large $ n $, we can approach $ \psi_{ n } $ as above with 
$
        \psi_{n}  \sim  n^2 / 12
$
and the absolute deviation with
$
        \Delta  \psi_{ n }   \sim  p / 2
$
to get 
$
           \Delta  \psi_{ n } / \psi_{n}  \sim  1 / n 
$
thereof.
For instance, for $ n = 200 $ we get 
$
           \Delta  \psi_{ n } / \psi_{n}  \approx 
$
5\,\textperthousand.
\end{itemize}
%
%
%
%
%
\begin{figure}
\centering
\includegraphics[width=12cm]{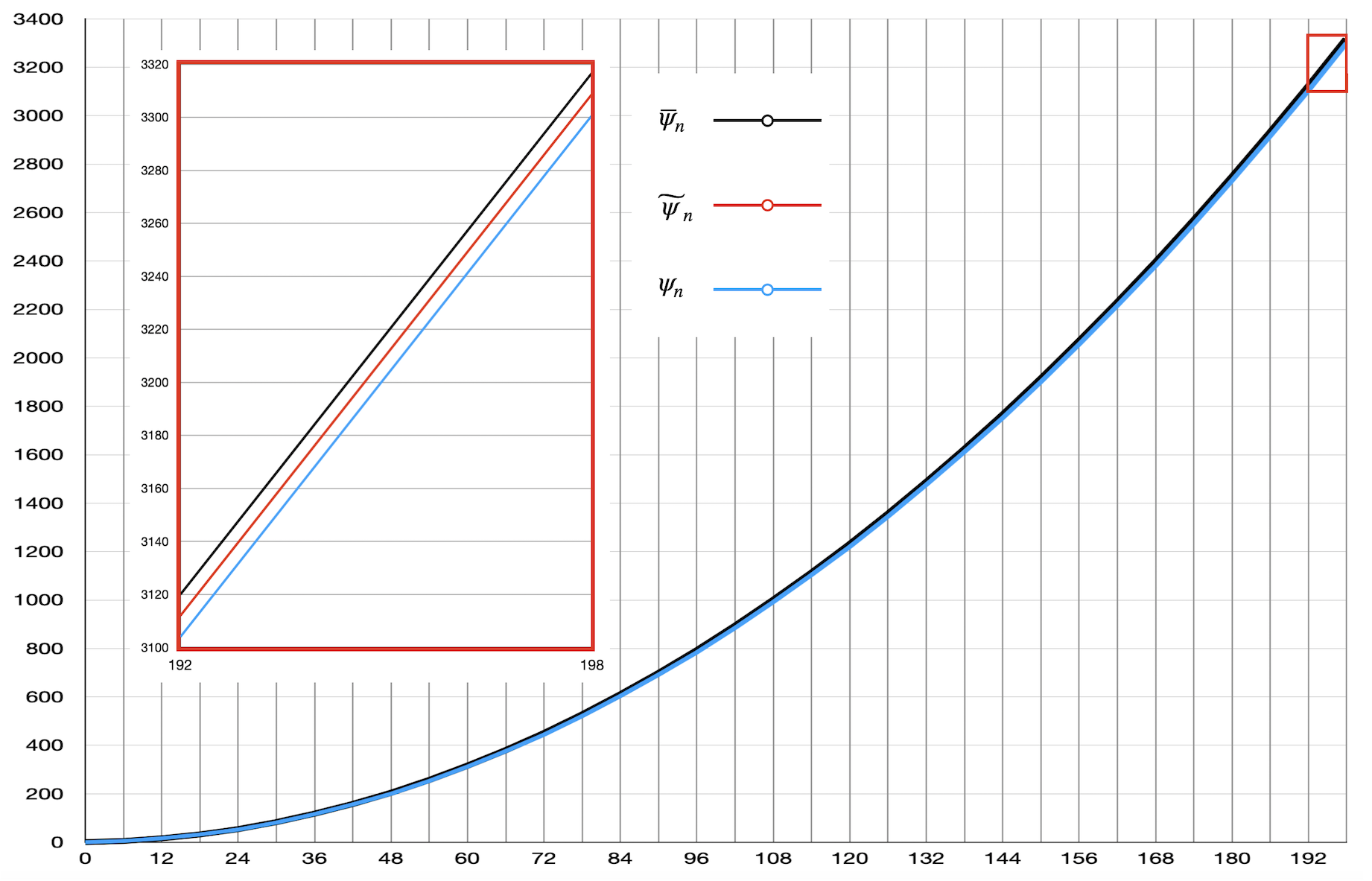}  
\caption{Evolution of the population of dominos in diamond $ \mathcal{D}_{ n } $ --  \texttt{Class} $ 00 $. 
Lower bound  $ \psi_{ n }  $, upper bound $  \overline{\psi}_n $ and median estimation 
$ 
          \widetilde{\psi}_n   =    ( \psi_n  +  \overline{\psi}_n ) / 2
$.
Each graph of $ ( \psi_{ n },  \overline{\psi}_n ) $  \textsf{vs} $ n $ is now a discrete dust of points in
$
         \mathbb{N} \times \mathbb{N}
$
located on a {\em smooth} parabola. 
Zoom on the small red window in the interval $ n \in [192, 198] $.
} 
\label{Figure: Conclusion -- Psi00-LB-M-UB --  }
\end{figure}
%
%
Fig.\,\ref{Figure: Conclusion -- Psi00-LB-M-UB --  }
shows the evolution of the population in diamond $ \mathcal{D}_{ n } $ --\,\texttt{Class} $ 00 $\,--
and the deviation between lower bound $ \psi_{ n }  $ and upper bound $  \overline{\psi}_n $.
Each graph is now a discrete dust of points in
$
         \mathbb{N} \times \mathbb{N}
$
located on a {\em smooth} parabola and this would also be the case for each of the six classes $ n m $. 
The divergence
$
           \overline{\psi}_n - \psi_{ n } 
$
within the small red window is highlighted.

A question arises whether the upper bound is excessive or not. 
The expression for $  \overline{\psi}_n  $ in 
Sect.\,\ref{Section: Injection of Disorder}
 was based on this harsh assumption that ``in a disordered system two lacunar voids are likely to join together, that yields a {\em potential} of one additional domino'' namely, {\em only} a potential.
 A median estimate
$ 
          \widetilde{\psi}_n   =    ( \psi_n  +  \overline{\psi}_n ) / 2
$ 
--\,as highlighted in
Fig.\,\ref{Figure: Conclusion -- Psi00-LB-M-UB --  }\,--
would perhaps be more appropriate by estimating a {\em fifty} percent chance of marrying two monominos (the lacunar voids) to give birth to a domino.
%
%
\paragraph{Conclusion.}                
%
%
In this study, we can deplore  the weakness of obtaining only an estimate $  \widetilde{\psi}_n $ in general and of obtaining the exact maximum only in a few cases.
This discrepancy is indeed an {\em open} problem.
It could be partially filled by better theory as well as by the use of simulation and, if possible, improved by high performance computing.
This technical study will serve as a companion paper to support further comparative study with the simulation model already tackled in
\cite{Hoffmann:Deserable-Seredynski-2021a}.
These new results will be examined elsewhere.
Obtaining new occurrences
$
           \psi_{ n }^{ ( s ) }  =   \overline{\psi}_{ n }
$
--\,where $   \psi_{ n }^{ ( s ) }  $ is the result of a simulation for a given $ n $\,-- would indeed be likely to expand the collection of exact values.
Be that as it may, this search for a maximum could at least remain the subject of a mathematical game, or of entertainment, or of a brain training challenge.
%
%
%
{ }
%
%
%
%
\section*{Appendix  -- Domino Enumeration}
%
%
%
This appendix brings together all the expressions and formulas that have been developed throughout this study and gives a numerical evaluation of them on a sample of 200 values of $ n $.

The first two tables \hspace{-1mm}
(Tabs.\,\ref{Table: Core, Wedge, Overall capacity in C1m ---  n  odd  }--\ref{Table: Core + Wedge  capacity in C0m ---  n  even  })
give an overview of the capacity of the regions in $ \mathcal{D}_{ n } $ according to the construction rule, and leading to the different expressions of $ \psi_n $.

The following eight tables
(Tabs.\,\ref{Table: Dominos in the Square --  n odd  }--\ref{Table: Domino enumeration -- Class 02 })
list the various counts of the maximum layouts of dominos in the {\em square} and the {\em diamond}.

The first two relate to the square $ \mathcal{S}_n $ of {\em odd} dimension
(Table \ref{Table: Dominos in the Square --  n odd  })
and then of {\em even} dimension
(Table \ref{Table: Dominos in the Square --  n even  }).

The next six relate to the diamond $ \mathcal{D}_n $ and are divided into two groups of three, first of {\em odd} dimension
(Tabs.\,\ref{Table: Domino enumeration -- Class 10 }--\ref{Table: Domino enumeration -- Class 12 })
then of {\em even} dimension
(Tabs.\,\ref{Table: Domino enumeration -- Class 00 }--\ref{Table: Domino enumeration -- Class 02 }).
The penultimate column lists the theoretical estimate $ \psi_{n} $ and the last column lists the upper bound $  \overline{\psi}_n $.

For the sake of clarity, each of the tables is preceded by the list of its various components.

\noindent\hrulefill
\newpage
%
%
%
%
\begin{table}[h]
\centering
\caption{Capacity 
$
        \psi_{ n } 
$ 
in 
$ 
          \mathcal{ D}_{n} 
$
according to the three classes  $ \texttt{C} 1m $  --- $ n = 2 m + 1 $. 
Detail of the different regions (core and wedges).
The 
$
        \psi_{ n } 
$ 
is expressed in various forms, either direct or following a recurrence relation.
}
\vspace{0mm}
\label{Table: Core, Wedge, Overall capacity in C1m ---  n  odd  }
\includegraphics[width=12cm]{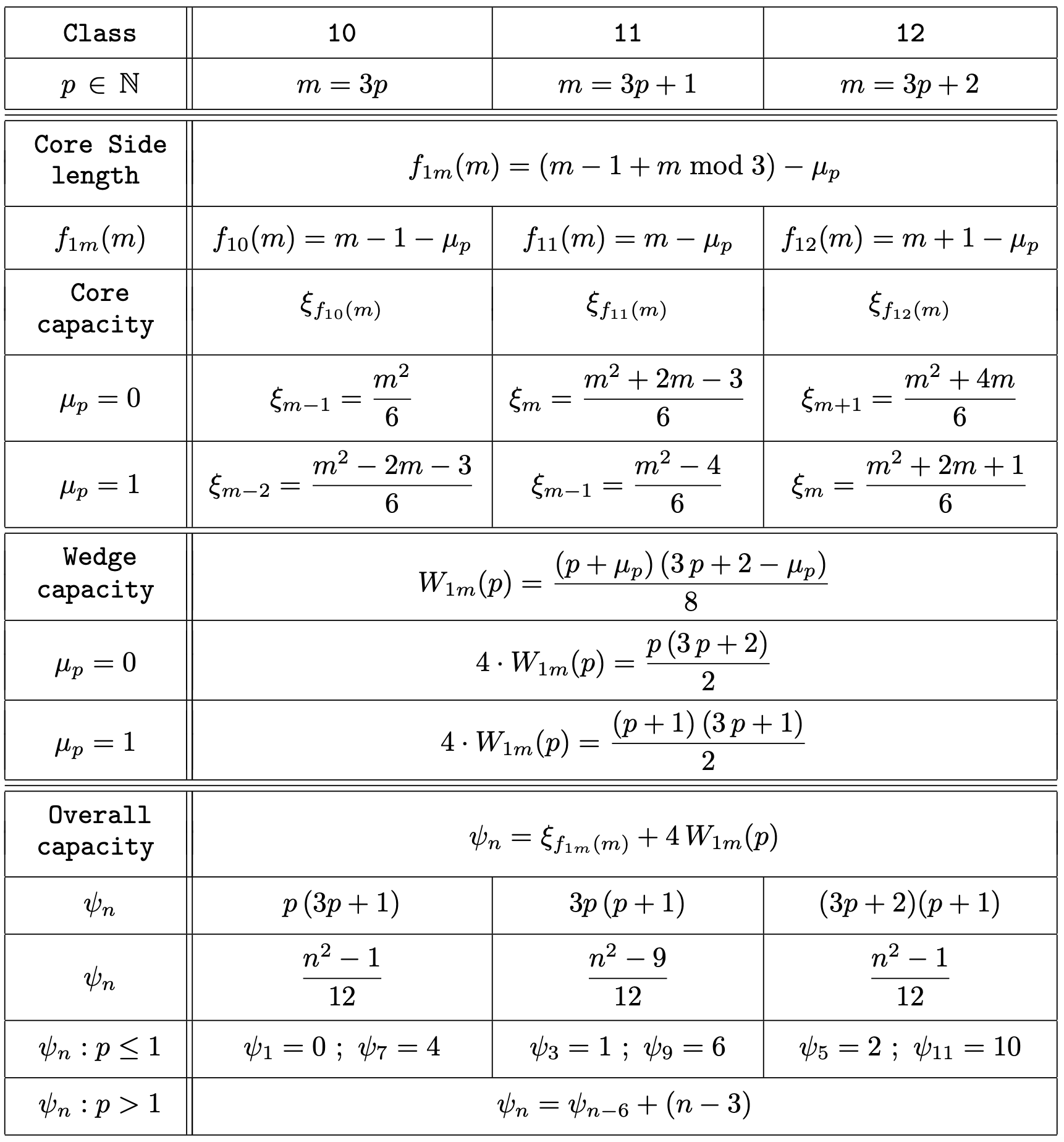}  
\end{table}
%
%
%
%
%
%
%
\begin{table}[h]
\centering
\caption{Part I --- Capacity 
$
        \psi_{ n } 
$ 
in 
$ 
          \mathcal{ D}_{n} 
$
according to the three classes  $ \texttt{C} 0m $  --- $ n = 2 m $. 
Detail of the different regions (core and wedges).
In the table, cells with a given color for the first and second subwedges denote a same wedge--pattern in the diamond (compare
Figs.\,\ref{Figure: C00-P15P18 }--\ref{Figure: C02-P17P20 }).
 }
\vspace{0mm}
\label{Table: Core + Wedge  capacity in C0m ---  n  even  }
\includegraphics[width=12cm]{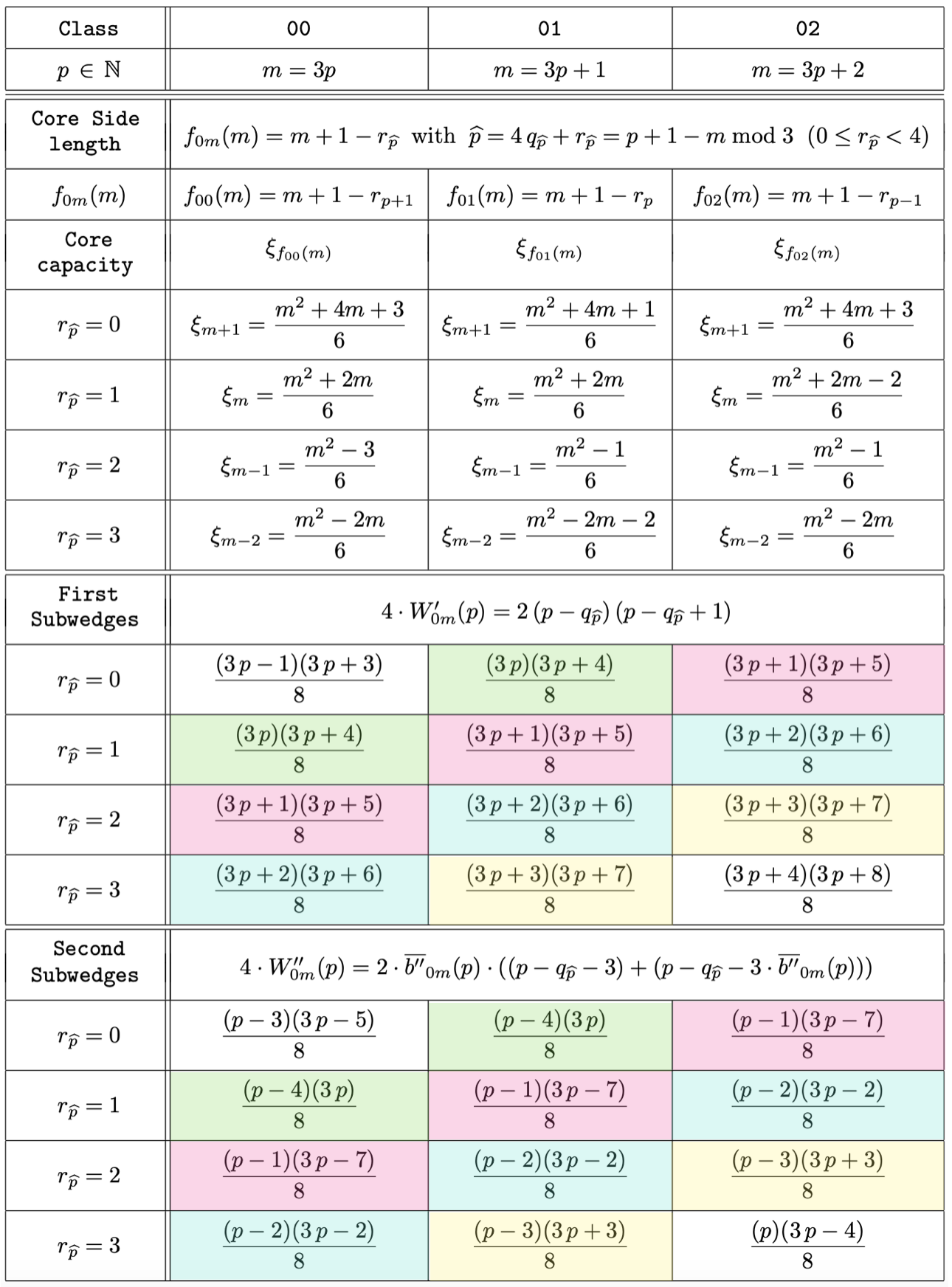}  
\end{table}
%
%
%
%
%
%
\begin{table}[h]
\centering
\addtocounter{table}{-1}
\caption{Part II --- Overall capacity 
$
        \psi_{ n } 
$ 
in 
$ 
          \mathcal{ D}_{n} 
$
according to the three classes  $ \texttt{C} 0m $  --- $ n = 2 m $. 
The 
$
        \psi_{ n } 
$ 
is expressed in various forms, either direct or following a recurrence relation.
A color highlights similar expressions for the direct form.
}
\vspace{0mm}
\label{Table: Overall capacity in C0m ---  n  even  }
\includegraphics[width=12cm]{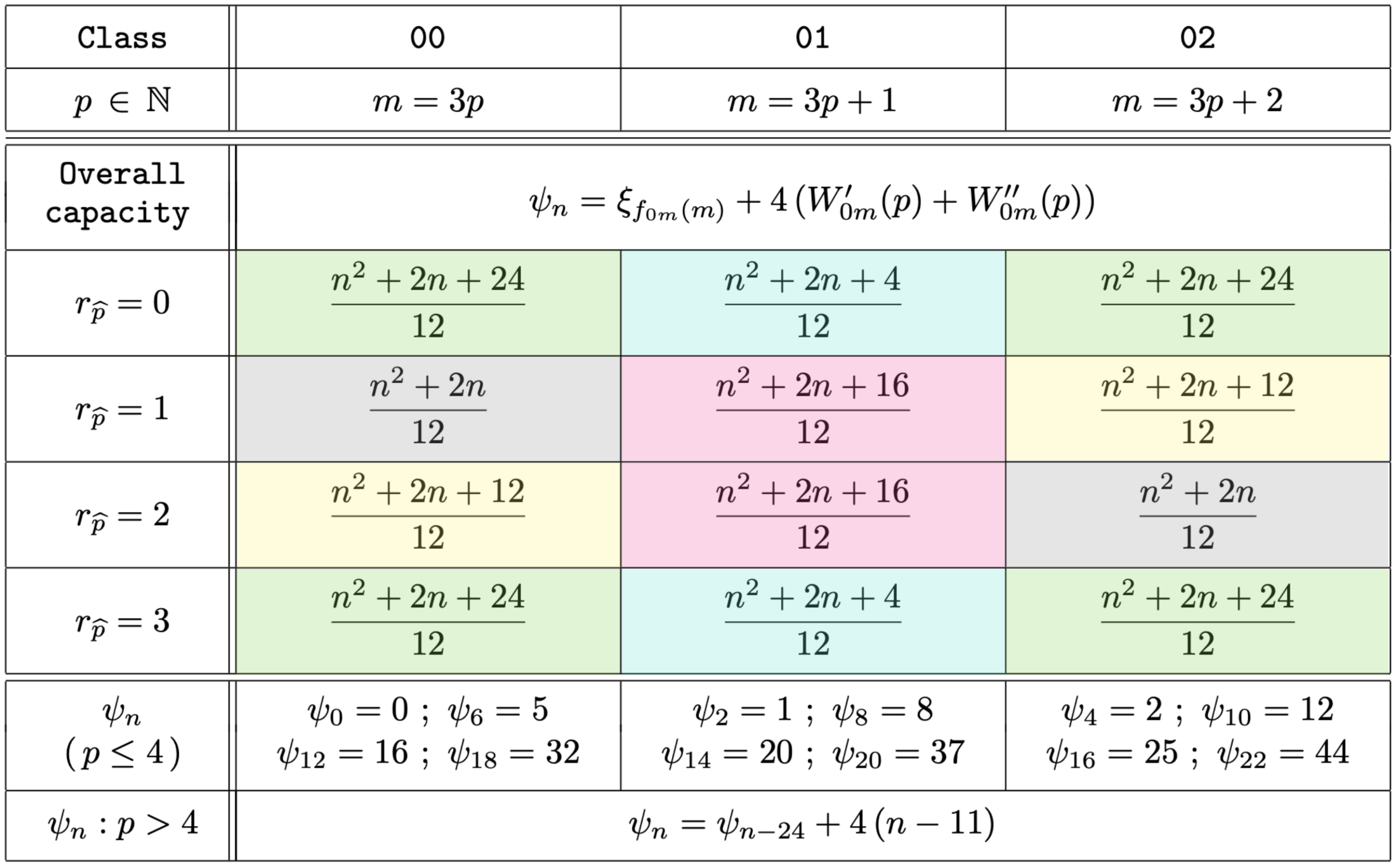}  
\end{table}
%
%
%
%
%
%
%
%
%
\begin{table}[h]
\centering
\caption{ Domino enumeration in the Square  --- 
$
           p \, \in  \, \mathbb{N} \,  ;  \hspace{2mm} n \, \mbox{ odd } \, ;
$ 
\newline\newline\normalsize{ 
$
        \xi_1 = 0,   \,\xi_3 = 2,  \, \xi_5 = 6   \hspace{3mm}   \mbox{and}   \hspace{3mm}     \xi_n  =  \xi_{n-6} + 2 ( n - 2 ) \, ;
$
\newline\newline
$
    \texttt{Class}  \ ( \dot{ 10 } )  :  \  m  \, = 3 \, p   \, ;  \hspace{2mm} n \, =  2 \, m  + 1 \, ; \hspace{3mm}   \xi_{ n }  = ( n - 1 ) ( n + 3 ) / \,6 \, ;   
$
\newline\newline
$
    \texttt{Class}  \ ( \dot{ 11 } )  : \  m  \, = 3 p + 1   \, ;  \hspace{2mm} n \, =  2 \, m  + 1 \, ; \hspace{3mm}   \xi_{ n }  = ( n - 1 ) ( n + 3 ) / \,6 \, ;   
$
\newline\newline
$
    \texttt{Class}  \ ( \dot{ 12 } )  : \  m  \, = 3 p + 2   \, ;  \hspace{2mm} n \, =  2 \, m  + 1 \, ; \hspace{3mm}   \xi_{ n }  = ( n + 1 )^2 / \,6 \, ;   
$
}
\newline
}
\label{Table: Dominos in the Square --  n odd }
\begin{tabular}{|c|||c|c|c|||c|c|c|||c|c|c|}
\hline
\rule{0pt}{11pt}
 $ \ \     \ \ $   	&   \multicolumn{ 3 }{c|||}{ Class   \ ($ \dot{ 10 } )$ } 
                         &   \multicolumn{ 3 }{c|||}{ Class   \ ($ \dot{ 11 } )$ }  
                         &   \multicolumn{ 3 }{c|}{   Class   \ ($ \dot{ 12 } )$  }  \\   [0.5ex]
\hline\hline	
\rule{0pt}{11pt}
 $ \ \ p  \ \ $   	&   $ \ \ m  \ \ $    &   $ \ \  n  \ \ $   &  $ \ \  \xi_n  $ \ \ 
 			&   $ \ \ m  \ \ $    &   $ \ \  n  \ \ $   &  $ \ \  \xi_n  $ \ \  
			&   $ \ \ m  \ \ $    &   $ \ \  n  \ \ $   &  $ \ \  \xi_n  $ \ \  \\   [0.5ex]
\hline\hline	
%
%
 0 	&   	0  	&    1		&    0		&   	1  	&    3		&    2		&   	2  	&    5		&    6		\\   
 1 	&   	3  	&    7		&    10	&   	4  	&    9		&    16	&   	5  	&    11	&    24	 \\   	 	 
 2  	&      6	&    13	&    32	&      7	&    15	&    42	&      8	&    17	&    54	 \\ 
 3  	&      9	&    19	&    66	&     10	&    21	&    80	&     11	&    23	&    96	 \\   
 4  	&    12	 &   25	&    112	&     13	 &   27	&    130	&     14	 &   29	&    150	\\  
 5  	&    15	 &   31	&    170	&     16	 &   33	&    192	&     17	 &   35	&    216	\\  \hline
6	&    18	 &    37	&    240	&     19	 &    39	&    266	 &    20	 &    41	&    294	 \\  
7	&    21	 &    43	&    322	&     22	 &    45	&    352	 &    23	 &    47	&    384	 \\  
8	&    24	 &    49	&    416	&     25	 &    51	&    450	&     26	 &    53	&    486	\\  
9	&    27	 &    55	&    522	&     28	 &    57	&    560	&     29	 &    59	&    600	\\  
10	&    30	 &    61	&    640	&     31	 &    63	&    682	&     32	 &    65	&    726	\\  
11 	&    33	 &    67	&    770	&    34	 &    69	&    816	&    35	 &    71	&    864	 \\  \hline
12 	&    36	 &    73	&    912	&     37	 &    75	&    962	&     38	 &    77	&     1014	\\ 
13 	&    39	 &    79	&    1066	&     40	 &    81	&    1120	&     41	 &    83	&    1176	\\ 
14 	&    42	 &    85	&    1232	&     43	 &    87	&    1290	&     44	 &    89	&    1350	\\ 
15 	&    45	 &    91	&    1410	&     46	 &    93	&    1472	&     47	 &    95	&    1536	\\ 
16 	&    48	 &    97	&    1600	&     49	 &    99	&    1666	&     50	 &    101	&    1734	\\ 
17 	&    51	 &    103	&    1802	&     52	 &   105	&    1872	&     53	 &   107	&    1944	\\  \hline
\end{tabular}
\end{table}
%
%
%
\newpage
%
%
%
%
\begin{table}[h]
\centering
\caption{ Domino enumeration in the Square  --- 
$
           p \, \in  \, \mathbb{N} \,  ;  \hspace{2mm} n \, \mbox{ even } \, ;
$ 
\newline\newline\normalsize{ 
$
        \xi_0 = 0,   \,\xi_2 = 1,  \, \xi_4 = 4   \hspace{3mm}   \mbox{and}   \hspace{3mm}     \xi_n  =  \xi_{n-6} + 2 ( n - 2 ) \, ;
$
\newline\newline
$
    \texttt{Class}  \ ( \dot{ 00 } ) :  \  m  \, = 3 \, p   \, ;  \hspace{2mm} n \, =  2 \, m  \, ; \hspace{3mm}   \xi_{ n }  = n ( n + 2 ) / \,6 \, ;    
$
\newline\newline
$
    \texttt{Class}  \ ( \dot{ 01 } ) : \  m  \, = 3 p + 1   \, ;  \hspace{2mm} n \, =  2 \, m  \, ; \hspace{3mm}   \xi_{ n }  =  ( n \, ( n + 2 ) - 2 ) / \,6 \, ;  
$
\newline\newline
$
    \texttt{Class}  \ ( \dot{ 02 } ) : \  m  \, = 3 p + 2   \, ;  \hspace{2mm} n \, =  2 \, m  \, ; \hspace{3mm}   \xi_{ n }  = n ( n + 2 ) / \,6 \, ;  
$
}
\newline
}
\label{Table: Dominos in the Square --  n even }
\begin{tabular}{|c|||c|c|c|||c|c|c|||c|c|c|}
\hline
\rule{0pt}{11pt}
  $ \ \     \ \ $   	&   \multicolumn{ 3 }{c|||}{ Class   \ ($ \dot{ 00 } )$ } 
                         &   \multicolumn{ 3 }{c|||}{ Class   \ ($ \dot{ 01 } )$ }  
                         &   \multicolumn{ 3 }{c|}{   Class   \ ($ \dot{ 02 } )$  }  \\   [0.5ex]
\hline\hline	
\rule{0pt}{11pt}
 $ \ \ p  \ \ $   	&   $ \ \ m  \ \ $    &   $ \ \  n  \ \ $   &  $ \ \  \xi_n  $ \ \ 
 			&   $ \ \ m  \ \ $    &   $ \ \  n  \ \ $   &  $ \ \  \xi_n  $ \ \  
			&   $ \ \ m  \ \ $    &   $ \ \  n  \ \ $   &  $ \ \  \xi_n  $ \ \  \\   [0.5ex]
\hline\hline	
%
%
 0 	&   	0  	&    0		&    0		&   	1  	&    2		&    1		&   	2  	&    4		&    4		\\   
 1 	&   	3  	&    6		&    8		&   	4  	&    8		&    13	&   	5  	&    10	&    20	 \\   	 	 
 2  	&      6	&    12	&    28	&      7	&    14	&    37	&      8	&    16	&    48	 \\ 
 3  	&      9	&    18	&    60	&     10	&    20	&    73	&     11	&    22	&    88	 \\   
 4  	&    12	 &   24	&    104	&     13	 &   26	&    121	&     14	 &   28	&    140	\\  
 5  	&    15	 &   30	&    160	&     16	 &   32	&    181	&     17	 &   34	&    204	\\  \hline
6	&    18	 &    36	&    228	&     19	 &    38	&    253	 &    20	 &    40	&    280	 \\  
7	&    21	 &    42	&    308	&     22	 &    44	&    337	 &    23	 &    46	&    368	 \\  
8	&    24	 &    48	&    400	&     25	 &    50	&    433	&     26	 &    52	&    468	\\  
9	&    27	 &    54	&    504	&     28	 &    56	&    541	&     29	 &    58	&    580	\\  
10	&    30	 &    60	&    620	&     31	 &    62	&    661	&     32	 &    64	&    704	\\  
11 	&    33	 &    66	&    748	&    34	 &    68	&    793	&    35	 &    70	&    840	 \\  \hline
12 	&    36	 &    72	&    888	&     37	 &    74	&    937	&     38	 &    76	&     988	\\ 
13 	&    39	 &    78	&    1040	&     40	 &    80	&    1093	&     41	 &    82	&    1148	\\ 
14 	&    42	 &    84	&    1204	&     43	 &    86	&    1261	&     44	 &    88	&    1320	\\ 
15 	&    45	 &    90	&    1380	&     46	 &    92	&    1441	&     47	 &    94	&    1504	\\ 
16 	&    48	 &    96	&    1568	&     49	 &    98	&    1633	&     50	 &    100	&    1700	\\ 
17 	&    51	 &    102	&    1768	&     52	 &   104	&    1837	&     53	 &    106	&    1908	\\  \hline
\end{tabular}
\end{table}
%
%
%
\newpage
%
%
%
%
\begin{table}[h]
\centering
\caption{ Domino enumeration in the Diamond --  \texttt{Class} $ 10 $.  \newline --- 
$
           p \, \in  \, \mathbb{N} \,  ;  \hspace{2mm} m  \, = 3 \, p   \, ;  \hspace{2mm} n \, =  2 \, m  + 1 \, ;
$ 
\newline\newline\normalsize{ 
$
           f_{10} ( m ) = m  - 1 - \mu_{ p } \, ;
          \hspace{3mm}   \psi_{n}  =   \xi_{ f_{10}( m ) }  + 4 \, W_{10}( p ) \, ;  
$
\newline\newline
$
        W_{10}( p )  =  ( p +  \mu_{ p } )  \,  ( 3 \, p + 2 -  \mu_{ p } ) / 8   \, ;
        \hspace{3mm}   \xi_{ f_{10}( m ) }  : 
$
refer to 
Tab.\,\ref{Table: Dominos in the Square --  n odd } ; 
\newline\newline
$
         \psi_{ 1 }  =   0 \, ;   \hspace{3mm}  \psi_n =  \psi_{ n - 6 }  +  ( n - 3 )  \hspace{4mm}  ( p > 0 )  \, ;
$
\newline\newline
$
         \psi_{ n }  =   p \, ( 3 p + 1 ) = ( n - 1 ) ( n + 1 )/ 12 \, ;   
$
\newline\newline
$
          \overline{\psi}_n  =  \psi_n  +  2     \hspace{4mm}  ( p > 1 ) 
$
}
\newline
}
\label{Table: Domino enumeration -- Class 10 }
\begin{tabular}{|c|c|c|c|c||c|c||c|c|}
\hline 
\rule{0pt}{11pt}
 $ \ \ n  \ \ $   &   $ \ \ m  \ \ $    &   $ \ \  p  \ \ $   &	$ \   \mu_{ p }  \   $   &   $   f_{10}( m )  $  &   
 $   \xi_{ f_{10}( m ) }  $   &  $  W_{10}( p ) $  &  $ \ \  \psi_{n}  $ \ \ &  $ \ \      \overline{\psi}_n  $ \ \   \\  [0.5ex]
\hline\hline	
%
%
 1 	&   	0  	&    0		&    0		&    --  	&    --  	&     0	&        0	&        0 	\\   
 7 	&   	3  	&    1		&    1		&    1 	&    0 	&     1	&        4	&        4	 \\   	 	 
 13  	&      6	&    2		&    0		&    5 	&    6		&     2 	&       14	&       16	 \\ 
 19  	&      9	&    3		&    1		&    7 	&    10  	&     5  	&       30	&       32 	 \\   
  25  	&    12	 &   4		&    0		&    11  	&    24	&     7	&       52	&       54	\\  
 31  	&    15	 &   5		&    1		&    13  	&    32  	&     12  	&       80	&       82	\\  \hline
37	&    18	 &    6	&    0		&    17  	&    54 	&      15 	&      114	&      116	 \\  
43	&    21	 &    7	&    1		&    19  	&    66 	&      22 	&      154	&      156	  \\  
49	&    24	 &    8	&    0		&    23  	&    96 	&      26 	&      200	&      202	 \\  
55	&    27	 &    9	&    1		&    25  	&    112 	&      35 	&      252	&      254	 \\  
61	&    30	 &    10	&    0		&    29  	&    150 	&      40 	&      310	&      312	  \\  
67 	&    33	 &    11	&    1		&    31    	&    170    	&      51  	&      374	&      376	  \\  \hline
73 	&    36	 &    12	&    0		&     35   	&     216   	&      57  	&      444	&      446	\\ 
79 	&    39	 &    13	&    1		&     37   	&     240   	&      70  	&      520	&      522	\\ 
85 	&    42	 &    14	&    0		&     41   	&     294   	&      77  	&      602 	&      604 	\\ 
91 	&    45	 &    15	&    1		&     43   	&     322   	&      92  	&      690 	&      692	\\ 
97 	&    48	 &    16	&    0		&     47   	&     384   	&      100  	&      784 	&      786	\\ 
103 	&    51	 &    17	&    1		&     49   	&     416   	&      117  	&      884	&      886	\\  \hline
109 	&    54	 &    18	&    0		&     53   	&     486   	&      126 	&      990	&      992 	\\ 	
115 	&    57	 &    19	&    1		&     55   	&     522   	&      145 	&      1102	&      1104	\\ 	
121 	&    60	 &    20	&    0		&     59   	&     600   	&      155 	&      1220	&      1222	\\ 	
127 	&    63	 &    21	&    1		&     61   	&     640   	&      176 	&      1344	&      1346	\\ 	
133 	&    66	 &    22	&    0		&     65   	&     726   	&      187 	&      1474 &      1476 \\ 	
139 	&    69	 &    23	&    1		&     67   	&     770   	&      210 	&      1610	&      1612	  \\  \hline	
 145 	&    72	 &    24	&    0		&     71   	&     864   	&      222 	&      1752	&      1754	 \\  	
 151 	&    75	 &    25	&    1		&     73   	&     912   	&      247 	&      1900	&      1902	 \\  	
 157 	&    78	 &    26	&    0		&     77   	&   1014   	&      260 	&      2054	&      2056 \\  	
 163 	&    81	 &    27	&    1		&     79   	&   1066   	&      287 	&      2214	&      2216 \\  	
 169 	&    84	 &    28	&    0		&     83   	&   1176   	&      301 	&      2380	&      2382 \\  	
 175 	&    87	 &    29	&    1		&     85   	&   1232   	&      330 	&      2552 &      2554  \\  \hline	
 181 	&    90	 &    30	&    0		&     89   	&   1350   	&      345 	&      2730	&      2732 \\  	
 187 	&    93	 &    31	&    1		&     91   	&   1410   	&      376 	&      2914	&      2916 \\  	
 193 	&    96	 &    32	&    0		&     95   	&   1536   	&      392 	&      3104	&      3106  \\  	
 199 	&    99	 &    33	&    1		&     97   	&   1600   	&      425 	&      3300	&      3302	 \\  \hline	
\end{tabular}
\end{table}
%
%
%
%
%
%
\begin{table}[h]
\centering
\caption{ Domino enumeration in the Diamond --  \texttt{Class} $ 11 $.  \newline --- 
$
           p \, \in  \, \mathbb{N} \,  ;  \hspace{2mm} m  \, = 3 \, p + 1  \, ;  \hspace{2mm} n \, =  2 \, m  + 1 \, ;
$
\newline\newline\normalsize{ 
$
           f_{ 11 } ( m ) = m  - \mu_{ p } \, ;
          \hspace{3mm}   \psi_{n}  =   \xi_{ f_{ 11 }( m ) }  + 4 \, W_{ 11 }( p ) \, ;  
$
\newline\newline
$
        W_{ 11 }( p )  =  ( p +  \mu_{ p } )  \,  ( 3 \, p + 2 -  \mu_{ p } ) / 8   \, ;
        \hspace{3mm}   \xi_{ f_{ 11 }( m ) }  : 
$
refer to 
Tab.\,\ref{Table: Dominos in the Square --  n odd } ; 
\newline\newline
$
         \psi_{ 3 }  =   1 \, ;   \hspace{3mm}  \psi_n =  \psi_{ n - 6 }  +  ( n - 3 )  \hspace{4mm}  ( p > 1 )  \, ;
$
\newline\newline
$
         \psi_{ n }  =   3\,p \, ( p + 1 ) = ( n - 3 ) ( n + 3 )/ 12 \, ;    \hspace{4mm}  ( p > 0 ) 
$
\newline\newline
$
         \overline{\psi}_n  =  \psi_n  +  2    \hspace{4mm}  ( p > 0 ) 
$
}
\newline
}
\label{Table: Domino enumeration -- Class 11 }
\begin{tabular}{|c|c|c|c|c||c|c||c|c|}
\hline
\rule{0pt}{11pt}
 $ \ \ n  \ \ $   &   $ \ \ m  \ \ $    &   $ \ \  p  \ \ $   &	$ \   \mu_{ p }  \   $   &   $   f_{11}( m )  $  &   
 $   \xi_{ f_{11}( m ) }  $   &  $  W_{11}( p ) $  &  $ \ \  \psi_{n}  $ \ \ &  $ \ \      \overline{\psi}_n $ \ \  \\  [0.5ex]
\hline\hline	
%
%
 3 	&   	1  	&    0		&    0		&    1 	&    0		&     0	&        1	&        1		 \\   
 9 	&   	4  	&    1		&    1		&    3 	&    2 	&     1	&        6	&        8		 \\   	 	 
 15  	&      7	&    2		&    0		&    7 	&    10	&     2 	&       18	&       20		 \\ 
 21  	&     10	&    3		&    1		&    9 	&    16  	&     5  	&       36	&       38 		 \\   
 27  	&    13	 &   4		&    0		&    13  	&    32	&     7	&       60	&       62		  \\  
 33  	&    16	 &   5		&    1		&    15  	&    42  	&     12  	&       90	&       92		  \\  \hline
39	&    19	 &    6	&    0		&    19  	&    66 	&      15 	&      126	&      128		  \\  
45	&    22	 &    7	&    1		&    21  	&    80 	&      22 	&      168	&      170		  \\  
51	&    25	 &    8	&    0		&    25  	&    112 	&      26 	&      216	&      218		  \\  
57	&    28	 &    9	&    1		&    27  	&    130 	&      35 	&      270	&      272		  \\  
63	&    31	 &    10	&    0		&    31  	&    170 	&      40 	&      330	&      332		  \\  
69 	&    34	 &    11	&    1		&    33    	&    192    	&      51  	&      396	&      398		\\  \hline
75 	&    37	 &    12	&    0		&     37   	&     240   	&      57  	&      468	&      470  		\\ 
81 	&    40	 &    13	&    1		&     39   	&     266   	&      70  	&      546	&      548  		\\ 
87 	&    43	 &    14	&    0		&     43   	&     322   	&      77  	&      630	&      632  		\\ 
93 	&    46	 &    15	&    1		&     45   	&     352   	&      92  	&      720 	&      722  		\\ 
99 	&    49	 &    16	&    0		&     49   	&     416   	&      100  	&      816 	&      818  		\\ 
105 	&    52	 &    17	&    1		&     51   	&     450   	&      117  	&      918	&      920  		\\  \hline
111 	&    55	 &    18	&    0		&     55   	&     522   	&      126 	&      1026	&      1028  	\\ 	
117 	&    58	 &    19	&    1		&     57   	&     560   	&      145 	&      1140	&      1142  	\\ 	
123 	&    61	 &    20	&    0		&     61   	&     640   	&      155 	&      1260	&      1262  	\\ 	
129 	&    64	 &    21	&    1		&     63   	&     682   	&      176 	&      1386	&      1388  	\\ 	
135 	&    67	 &    22	&    0		&     67   	&     770   	&      187 	&      1518	&      1520  	\\ 	
141 	&    70	 &    23	&    1		&     69   	&     816   	&      210 	&      1656	&      1658  	\\  \hline	
 147 	&    73	 &    24	&    0		&     73   	&     912   	&      222 	&      1800	&      1802  	\\  	
 153 	&    76	 &    25	&    1		&     75   	&     962   	&      247 	&      1950	&      1952  	\\  	
 159 	&    79	 &    26	&    0		&     79   	&   1066   	&      260 	&      2106	&      2108  	\\  	
 165 	&    82	 &    27	&    1		&     81   	&   1120   	&      287 	&      2268	&      2270  	\\  	
 171 	&    85	 &    28	&    0		&     85   	&   1232   	&      301 	&      2436	&      2438  	\\  	
 177 	&    88	 &    29	&    1		&     87   	&   1290   	&      330 	&      2610	&      2612  	\\  \hline	
 183 	&    91	 &    30	&    0		&     91   	&   1410   	&      345 	&      2790	&      2792  	\\  	
 189 	&    94	 &    31	&    1		&     93   	&   1472   	&      376 	&      2976	&      2978  	\\  	
 195 	&    97	 &    32	&    0		&     97   	&   1600   	&      392 	&      3168	&      3170  	\\  	
 201 	&    100	 &    33	&    1		&     99   	&   1666   	&      425 	&      3366	&      3368  	\\  \hline	
\end{tabular}
\end{table}
%
%
%
%
%
\begin{table}[h]
\centering
\caption{ Domino enumeration in the Diamond --  \texttt{Class} $ 12 $.  \newline --- 
$
           p \, \in  \, \mathbb{N} \,  ;  \hspace{2mm} m  \, = 3 \, p + 2  \, ;  \hspace{2mm} n \, =  2 \, m  + 1 \, ;
$
\newline\newline\normalsize{ 
$
           f_{ 12 } ( m ) = m + 1  - \mu_{ p } \, ;
          \hspace{3mm}   \psi_{n}  =   \xi_{ f_{ 12 }( m ) }  + 4 \, W_{ 12 }( p ) \, ;  
$
\newline\newline
$
        W_{ 12 }( p )  =  ( p +  \mu_{ p } )  \,  ( 3 \, p + 2 -  \mu_{ p } ) / 8   \, ;
        \hspace{3mm}   \xi_{ f_{ 12 }( m ) }  : 
$
refer to 
Tab.\,\ref{Table: Dominos in the Square --  n odd } ; 
\newline\newline
$
         \psi_{ 5 }  =   2 \, ;   \hspace{3mm}  \psi_n =  \psi_{ n - 6 }  +  ( n - 3 )  \hspace{4mm}  ( p > 0 )  \, ;
$
\newline\newline
$
         \psi_{ n }  =   ( 3\,p + 2 ) \, ( p + 1 ) = ( n - 1 ) ( n + 1 )/ 12 \, ;   
$
\newline\newline
$
         \overline{\psi}_n =  \psi_n  +  2 \, ( 1 - \mu_{ p } )   \hspace{4mm}  ( p > 0 )  
$
}
\newline
}
\label{Table: Domino enumeration -- Class 12 }
\begin{tabular}{|c|c|c|c|c||c|c||c|c|}
\hline
\rule{0pt}{11pt}
 $ \ \ n  \ \ $   &   $ \ \ m  \ \ $    &   $ \ \  p  \ \ $   &	$ \   \mu_{ p }  \   $   &   $   f_{12}( m )  $  &   
 $   \xi_{ f_{12}( m ) }  $   &  $  W_{12}( p ) $  &  $ \ \  \psi_{n}  $ \ \ &  $ \ \      \overline{\psi}_n  $ \ \  \\   [0.5ex]
\hline\hline	
%
%
 5 	&   	2  	&    0		&    0		&    3 	&    2		&     0	&       2	&       2		 \\   
 11 	&   	5  	&    1		&    1		&    5 	&    6 	&     1	&       10	&       10		 \\   	 	 
 17  	&      8	&    2		&    0		&    9 	&    16	&     2 	&       24	&       26		 \\ 
 23  	&     11	&    3		&    1		&    11 	&    24  	&     5  	&       44	&       44 		 \\   
 29  	&    14	 &   4		&    0		&    15  	&    42	&     7	&       70	&       72		  \\  
 35  	&    17	 &   5		&    1		&    17  	&    54  	&     12  	&       102	&       102		  \\  \hline
41	&    20	 &    6	&    0		&    21  	&    80 	&      15 	&      140	&      142		  \\  
47	&    23	 &    7	&    1		&    23  	&    96 	&      22 	&      184	&      184		  \\  
53	&    26	 &    8	&    0		&    27  	&    130 	&      26 	&      234	&      236		  \\  
59	&    29	 &    9	&    1		&    29  	&    150 	&      35 	&      290	&      290		  \\  
65	&    32	 &    10	&    0		&    33  	&    192 	&      40 	&      352	&      354		  \\  
71 	&    35	 &    11	&    1		&    35    	&    216    	&      51  	&      420	&      420		\\  \hline
77 	&    38	 &    12	&    0		&     39   	&     266   	&      57  	&      494	&      496  		\\ 
83 	&    41	 &    13	&    1		&     41   	&     294   	&      70  	&      574	&      574  		\\ 
89 	&    44	 &    14	&    0		&     45   	&     352   	&      77  	&      660	&      662  		\\ 
95 	&    47	 &    15	&    1		&     47   	&     384   	&      92  	&      752	&      752  		\\ 
101 	&    50	 &    16	&    0		&     51   	&     450   	&      100  	&      850	&      852  		\\ 
107 	&    53	 &    17	&    1		&     53   	&     486   	&      117  	&      954	&      954  		\\  \hline
113 	&    56	 &    18	&    0		&     57   	&     560   	&      126 	&      1064	&      1066  	\\ 	
119 	&    59	 &    19	&    1		&     59   	&     600   	&      145 	&      1180	&      1180  	\\ 	
125 	&    62	 &    20	&    0		&     63   	&     682   	&      155 	&      1302	&      1304  	\\ 	
131 	&    65	 &    21	&    1		&     65   	&     726   	&      176 	&      1430	&      1430  	\\ 	
137 	&    68	 &    22	&    0		&     69   	&     816   	&      187 	&      1564	&      1566  	\\ 	
143 	&    71	 &    23	&    1		&     71   	&     864   	&      210 	&      1704	&      1704  	\\  \hline	
 149 	&    74	 &    24	&    0		&     75   	&     962   	&      222 	&      1850	&      1852  	\\  	
 155 	&    77	 &    25	&    1		&     77   	&    1014   &      247 	&      2002	&      2002  	\\  	
 161 	&    80	 &    26	&    0		&     81   	&   1120   	&      260 	&      2160	&      2162  	\\  	
 167 	&    83	 &    27	&    1		&     83   	&   1176   	&      287 	&      2324	&      2324  	\\  	
 173 	&    86	 &    28	&    0		&     87   	&   1290   	&      301 	&      2494	&      2496  	\\  	
 179 	&    89	 &    29	&    1		&     89   	&   1350   	&      330 	&      2670	&      2670  	\\  \hline	
 185 	&    92	 &    30	&    0		&     93   	&   1472   	&      345 	&      2852	&      2854  	\\  	
 191 	&    95	 &    31	&    1		&     95   	&   1536   	&      376 	&      3040	&      3040  	\\  	
 197 	&    98	 &    32	&    0		&     99   	&   1666   	&      392 	&      3234	&      3236  	\\  	
 203 	&    101	 &    33	&    1		&    101   	&   1734   	&      425 	&      3434	&      3434  	\\  \hline	
\end{tabular}
\end{table}
%
%
%
%
%
%
%
\begin{table}[h]
\centering
\caption{Domino enumeration in the Diamond --  \texttt{Class} $ 00 $.  \newline --- 
$
           p \, \in  \, \mathbb{N} \,  ;  \hspace{2mm} m  \, = 3 \, p  \, ;  \hspace{2mm} n \, =  2 \, m  \, ;
$
\newline\newline\normalsize{
$
            p + 1 = 4 \, q_{ p + 1 } +  r_{ p + 1 }  \hspace{2mm}  ( 0  \le  r_{ p + 1 } < 4 )  \, ;  \hspace{8mm}
            f_{00} ( m ) = m  + 1 - r_{ p + 1 } \, ; 
$
\newline\newline
$
          W'_{00} ( p )  =  ( p - q_{ p + 1 } ) \, ( p - q_{ p + 1 } + 1 ) / 2   \, ;   \hspace{8mm}
           \xi_{ f_{00}( m ) }  :  
$
refer to 
Tab.\,\ref{Table: Dominos in the Square --  n even } ; 
\newline\newline
$
          W''_{00}( p )  =  \overline{b''}_{00} ( m )  \cdot
          ( ( p - q_{ p + 1 }  -  3 )  +  ( p - q_{ p + 1 }  -  3  \cdot  \overline{b''}_{00} ( m ) ) \,  / 2    \, ; 
$
\newline\newline
$
          \overline{b''}_{00} ( m )   =   ( q_{ p + 1 } - 1 )  +   \delta_{ r_{ p + 1 } }    \hspace{3mm} \mbox{with} \hspace{3mm} 
          \delta_{ r_{ p + 1 } }  =   \Bigl\lfloor  r_{ p + 1 } / 3  \Bigr\rceil  ; 
$
\newline\newline
$
         \psi_{n}  =   \xi_{ f_{00}( m ) }  + 4 \, ( W'_{00} ( p ) + W''_{00} ( p ) ) \, ;    \hspace{3mm} 
          \psi_{ 0 }  =   0 \, ;   \hspace{1mm}  \psi_n =  \psi_{ n - 24 } \ +  \ 4 \ ( n - 11 )  \hspace{2mm}  ( p  \ge  4 )  \, ;
$
\newline\newline
$
          \overline{\psi}_n =  \psi_n  +  2 \, \lfloor \, p / 4  \,  \rfloor 
$
}
\vspace{-1mm}
}
\label{Table: Domino enumeration -- Class 00 }
\begin{tabular}{|c|c|c|c|c||c|c|c||c|c|}
\hline
\rule{0pt}{11pt}
 $ \ n \ $   &   $  \ m \  $    &   $  \ p \  $   &	$  r_{ p + 1 } \ $   &   $ f_{00}( m ) $  &   $ \xi_{ f_{00}( m ) } $   &
 $ W'_{00} ( p ) $  & $ W''_{00}( p ) $  &   $  \  \psi_{n}  $ \ &  $ \      \overline{\psi}_n $ \  \\   [0.5ex]
\hline\hline	
%
%
 0 	&   	0  	&    0		&    1		&    0 	&    0		&     0	&     0	&        0	&       0		 \\   
 6 	&   	3  	&    1		&    2		&    2 	&    1 	&     1	&     0	&       5	&       5		 \\   	 	 
 12  	&      6	&    2		&    3		&    4 	&    4		&     3 	&     0 	&       16	&       16		 \\   \hline
18  	&      9	&     3	&    0		&    10 	&    20  	&     3  	&      0  	&       32	&       32 		 \\   
24  	&    12	 &    4	&    1		&    12  	&    28	&     6	&      0	&       52	&       54		  \\  
30	&    15	 &    5	&    2		&    14  	&    37 	&      10 	&      1 	&       81	&       83		  \\  
36	&    18	 &    6	&    3		&    16  	&    48 	&      15 	&      2 	&      116	&      118		  \\  \hline
42	&    21	 &    7	&    0		&    22  	&    88 	&      15 	&      2 	&      156	&      158		  \\  
48	&    24	 &    8	&    1		&    24  	&    104 	&      21 	&      3 	&      200	&      204		  \\  
54	&    27	 &    9	&    2		&    26  	&    121 	&      28 	&      5 	&      253	&      257		  \\  
60 	&    30	 &    10	&    3		&    28    	&    140    	&      36  	&      7  	&      312	&      316		  \\  \hline
66 	&    33	 &    11	&    0		&     34   	&     204   	&      36  	&      7  	&      376	&      380  		\\ 
72 	&    36	 &    12	&    1		&     36   	&     228   	&      45  	&      9  	&      444	&      450  		\\ 
78 	&    39	 &    13	&    2		&     38   	&     253   	&      55  	&      12  	&      521	&      527  		\\ 
84 	&    42	 &    14	&    3		&     40   	&     280   	&      66  	&      15  	&      604	&      610  		\\   \hline
90 	&    45	 &    15	&    0		&     46   	&     368   	&      66  	&      15  	&      692	&      698  		\\ 
96 	&    48	 &    16	&    1		&     48   	&     400   	&      78  	&      18  	&      784	&      792  		\\  	
102 	&    51	 &    17	&    2		&     50   	&     433   	&      91 	&      22 	&      885	&      893  		\\ 	
108 	&    54	 &    18	&    3		&     52   	&     468   	&      105 	&      26 	&      992	&      1000	  	\\  \hline
114 	&    57	 &    19	&    0		&     58   	&     580   	&      105 	&      26 	&      1104	&      1112  	\\ 	
120 	&    60	 &    20	&    1		&     60   	&     620   	&      120 	&      30 	&      1220	&      1230  	\\ 	
126 	&    63	 &    21	&    2		&     62   	&     661   	&      136 	&      35 	&      1345	&      1355  	\\ 	
132 	&    66	 &    22	&    3		&     64   	&     704   	&      153 	&      40 	&      1476	&      1486  	\\  \hline	
 138 	&    69	 &    23	&    0		&     70   	&     840   	&      153 	&      40 	&      1612	&      1622  	\\  	
 144 	&    72	 &    24	&    1		&     72   	&     888   	&      171 	&      45 	&      1752	&      1764  	\\  	
 150 	&    75	 &    25	&    2		&     74   	&     937   	&      190 	&      51 	&      1901	&      1913  	\\  	
 156 	&    78	 &    26	&    3		&     76   	&     988   	&      210 	&      57 	&      2056	&      2068  	\\  \hline
 162 	&    81	 &    27	&    0		&     82   	&   1148   	&      210 	&      57 	&      2216	&      2228  	\\  	
 168 	&    84	 &    28	&    1		&     84   	&   1204   	&      231 	&      63 	&      2380	&      2394  	\\  	
 174 	&    87	 &    29	&    2		&     86   	&   1261   	&      253 	&      70 	&      2553 &      2567  	\\  	
 180 	&    90	 &    30	&    3		&     88   	&   1320   	&      276 	&      77 	&      2732	&      2746  	\\  \hline
 186 	&    93	 &    31	&    0		&     94   	&   1504   	&      276 	&      77 	&      2916	&      2930  	\\  	
 192 	&    96	 &    32	&    1		&     96   	&   1568   	&      300 	&      84 	&      3104	&      3120  	\\  	
 198 	&    99	 &    33	&    2		&     98   	&   1633   	&      325 	&      92 	&      3301	&      3317  	\\  \hline	
\end{tabular}
\end{table}
%
%
%
%
%
%
\begin{table}[h]
\centering
\caption{Domino enumeration in the Diamond --  \texttt{Class} $ 01 $.  \newline --- 
$
           p \, \in  \, \mathbb{N} \,  ;  \hspace{2mm} m  \, = 3 \, p + 1  \, ;  \hspace{2mm} n \, =  2 \, m  \, ;
$
\newline\newline\normalsize{
$
            p  = 4 \, q_{ p  } +  r_{ p  }  \hspace{2mm}  ( 0  \le  r_{ p  } < 4 )  \, ;  \hspace{20mm}
            f_{01} ( m ) = m  + 1 - r_{ p } \, ; 
$
\newline\newline
$
          W'_{01} ( p )  =  ( p - q_{ p  } ) \, ( p - q_{ p  } + 1 ) / 2   \, ;    \hspace{10mm}
           \xi_{ f_{01}( m ) }  
$
refer to 
Tab.\,\ref{Table: Dominos in the Square --  n even } ; 
\newline\newline
$
          W''_{01}( p )  =  \overline{b''}_{01} ( m )  \cdot
          ( ( p - q_{ p  }  -  3 )  +  ( p - q_{ p  }  -  3  \cdot  \overline{b''}_{01} ( m ) ) \,  / 2    \, ; 
$
\newline\newline
$
          \overline{b''}_{01} ( m )   =   ( q_{ p } - 1 )  +   \delta_{ r_{ p } }    \hspace{3mm} \mbox{with} \hspace{3mm} 
          \delta_{ r_{ p } }  =   \Bigl\lfloor (  r_{ p }  + 1 ) / 3  \Bigr\rceil  ; 
$
\newline\newline
$
         \psi_{n}  =   \xi_{ f_{01}( m ) }  + 4 \, ( W'_{01} ( p ) + W''_{01} ( p ) ) \, ;    \hspace{3mm}
          \psi_{ 2 }  =   1 \, ;   \hspace{1mm}  \psi_n =  \psi_{ n - 24 } \ +  \ 4 \ ( n - 11 )  \hspace{2mm}  ( p  \ge  4 )  \, ;
$
\newline\newline
$
          \overline{\psi}_n =  \psi_n  +  2 \, \lfloor  ( p + 1 ) / 4  \rfloor 
$
}
\vspace{-1mm}
}
\label{Table: Domino enumeration -- Class 01 }
\begin{tabular}{|c|c|c|c|c||c|c|c||c|c|}
\hline
\rule{0pt}{11pt}
 $ \ n \ $   &   $  \ m \  $    &   $  \ p \  $   &	$  r_{ p } \ $   &   $ f_{01}( m ) $  &   $ \xi_{ f_{01}( m ) } $   &
 $ W'_{01} ( p ) $  & $ W''_{01}( p ) $  &   $  \  \psi_{n}  $ \ &  $ \   \overline{\psi}_n $ \  \\    [0.5ex]
\hline\hline	
%
%
 2 	&   	1  	&    0		&    0		&    2 	&    1		&     0	&     0	&        1	&      1		 \\   
 8 	&   	4  	&    1		&    1		&    4 	&    4 	&     1	&     0	&       8	&       8		 \\   	 	 
14  	&      7	&     2	&    2		&    6 	&    8		&     3 	&     0 	&       20	&      20		 \\   
20  	&      10	&     3	&    3		&    8 	&    13  	&     6  	&      0  	&       37 	&     39 		 \\   \hline
26  	&    13	 &    4	&    0		&    14  	&    37	&      6	&      0	&       61	&       63		  \\  
32	&    16	 &    5	&    1		&    16  	&    48 	&      10 	&      1 	&       92	&       94		  \\  
38	&    19	 &    6	&    2		&    18  	&    60 	&      15 	&      2 	&      128	&      130		  \\  
44	&    22	 &    7	&    3		&    20  	&    73 	&      21 	&      3 	&      169	&      173		  \\  \hline
50	&    25	 &    8	&    0		&    26  	&    121 	&      21 	&      3 	&      217	&      221		  \\  
56	&    28	 &    9	&    1		&    28  	&    140 	&      28 	&      5 	&      272	&      276		  \\  
62 	&    31	 &    10	&    2		&    30    	&    160    	&      36  	&      7  	&      332	&      336		  \\  
68 	&    34	 &    11	&    3		&    32   	&     181   	&      45  	&      9  	&      397	&      403  		 \\  \hline
74 	&    37	 &    12	&    0		&     38   	&     253   	&      45  	&      9  	&      469	&      475  		\\ 
80 	&    40	 &    13	&    1		&     40   	&     280   	&      55  	&      12  	&      548	&      554  		\\ 
86 	&    43	 &    14	&    2		&     42   	&     308   	&      66  	&      15  	&      632	&      638  		\\   
92 	&    46	 &    15	&    3		&     44   	&     337   	&      78  	&      18  	&      721	&      729  		\\  \hline
98 	&    49	 &    16	&    0		&     50   	&     433   	&      78  	&      18  	&      817 	&      825  		\\  	
104 	&    52	 &    17	&    1		&     52   	&     468   	&      91 	&      22 	&      920 	&      928  		\\ 	
110 	&    55	 &    18	&    2		&     54   	&     504   	&      105 	&      26 	&      1028	&      1036	  	\\  
116 	&    58	 &    19	&    3		&     56   	&     541   	&      120 	&      30 	&      1141	&      1151  	\\  \hline
122 	&    61	 &    20	&    0		&     62   	&     661   	&      120 	&      30 	&      1261 &      1271  	\\ 	
128 	&    64	 &    21	&    1		&     64   	&     704   	&      136 	&      35 	&      1388	 &      1398  	\\ 	
134 	&    67	 &    22	&    2		&     66   	&     748   	&      153 	&      40 	&      1520 &      1530  	\\  	
140 	&    70	 &    23	&    3		&     68   	&     793   	&      171 	&      45 	&      1657	&      1669  	\\  \hline	
 146 	&    73	 &    24	&    0		&     74   	&     937   	&      171 	&      45 	&      1801	&      1813  	\\  	
 152 	&    76	 &    25	&    1		&     76   	&     988   	&      190 	&      51 	&      1952	&      1964  	\\  	
 158 	&    79	 &    26	&    2		&     78   	&    1040  &      210 	&      57 	&      2108	&      2120  	\\   
 164 	&    82	 &    27	&    3		&     80   	&   1093   	&      231 	&      63 	&      2269	&      2283  	\\  \hline	
 170 	&    85	 &    28	&    0		&     86   	&   1261   	&      231 	&      63 	&      2437	&      2451  	\\  	
 176 	&    88	 &    29	&    1		&     88   	&   1320   	&      253 	&      70 	&      2612	&      2626  	\\  	
 182 	&    91	 &    30	&    2		&     90   	&   1380   	&      276 	&      77 	&      2792	&      2806  	\\  
 188 	&    94	 &    31	&    3		&     92   	&   1441   	&      300 	&      84 	&      2977	&      2993  	\\  \hline	
 194 	&    97	 &    32	&    0		&     98   	&   1633   	&      300 	&      84 	&      3169	&      3185  	\\  	
 200 	&    100	 &    33	&    1		&     100   	&   1700   	&      325 	&      92 	&      3368 &      3384  	\\  \hline	
\end{tabular}
\end{table}
%
%
%
%
%
%
\begin{table}[h]
\centering
\caption{Domino enumeration in the Diamond --  \texttt{Class} $ 02 $.  \newline --- 
$
           p \, \in  \, \mathbb{N} \,  ;  \hspace{2mm} m  \, = 3 \, p + 2  \, ;  \hspace{2mm} n \, =  2 \, m  \, ;
$
\newline\newline\normalsize{
$
             p - 1  = 4 \, q_{ p  - 1 } +  r_{ p  - 1 }  \hspace{2mm}  ( 0  \le  r_{ p - 1 } < 4 )  \, ;  \hspace{8mm}
            f_{02} ( m ) = m  + 1 - r_{ p - 1 } \, ; 
$
\newline\newline
$
          W'_{02} ( p )  =  ( p - q_{ p - 1 } ) \, ( p - q_{ p - 1 } + 1 ) / 2   \, ;   \hspace{8mm}
           \xi_{ f_{02}( m ) }  :  
$
refer to 
Tab.\,\ref{Table: Dominos in the Square --  n even } ; 
\newline\newline
$
          W''_{02}( p )  =  \overline{b''}_{02} ( m )  \cdot
          ( ( p - q_{ p - 1 }  -  3 )  +  ( p - q_{ p - 1 }  -  3  \cdot  \overline{b''}_{02} ( m ) ) \,  / 2    \, ; 
$
\newline\newline
$
          \overline{b''}_{02} ( m )  =   ( q_{ p - 1 } - 1 )  + \delta_{ r_{ p - 1 } }   \hspace{3mm} \mbox{with} \hspace{3mm} 
           \delta_{ r_{ p - 1 } }   =   \Bigl\lfloor (  r_{ p - 1 }  +  2 ) / 3  \Bigr\rceil  ; 
$
\newline\newline
$
         \psi_{n}  =   \xi_{ f_{02}( m ) }  + 4 \, ( W'_{02} ( p ) + W''_{02} ( p ) ) \, ;    \hspace{3mm}
          \psi_{ 4 }  =   2 \, ;   \hspace{1mm}  \psi_n =  \psi_{ n - 24 } \ +  \ 4 \ ( n - 11 )  \hspace{2mm}  ( p  >  4 )  \, ;
$
\newline\newline
$
          \overline{\psi}_n =  \psi_n  +  2 \, \lfloor  ( p + 2 ) / 4  \rfloor   -  \lambda_{ r_{ p - 1 } }  \, ;  \hspace{2mm}
           \lambda_{ r_{ p - 1 } } = 1 \hspace{2mm} \mbox{if} \hspace{2mm} r_{ p - 1 } = 1 
           \hspace{2mm}  (  \lambda_{ r_{ p - 1 } } = 0  \hspace{1mm}  \mbox{otherwise} )
$
}
\vspace{-3mm}
}
\label{Table: Domino enumeration -- Class 02 }
\begin{tabular}{|c|c|c|c|c||c|c|c||c|c|}
\hline
\rule{0pt}{10pt}
 $ \ n \ $   &   $  \ m \  $    &   $  \ p \  $   &	$  r_{ p - 1 } \ $   &   $ f_{02}( m ) $  &   $ \xi_{ f_{02}( m ) } $   &
 $ W'_{02} ( p ) $  & $ W''_{02}( p ) $  &   $  \  \psi_{n}  $ \ &  $ \  \overline{\psi}_n  $ \  \\     [0.4ex]
\hline\hline	
%
%
 4 	&   	2  	&    0		&    --	&    0 	&    0		&     --	&     0	&        2	&        2		 \\  \hline 
10 	&   	5  	&     1	&    0		&    6 	&    8		&     1	&     0	&       12	&       12		 \\   
16 	&   	8  	&     2	&    1		&    8 	&    13 	&     3	&     0	&       25	&       26		 \\   	 	 
22  	&      11	&     3	&    2		&    10 	&    20	&     6 	&     0 	&       44	&       46		 \\   
28  	&      14	&     4	&    3		&    12 	&    28  	&     10  	&      1  	&       72	&       74 		 \\   \hline
34  	&    17	 &    5	&    0		&    18  	&    60	&      10	&      1	&      104	&      106		  \\  
40	&    20	 &    6	&    1		&    20  	&    73 	&      15 	&      2 	&      141	&      144		  \\  
46	&    23	 &    7	&    2		&    22  	&    88 	&      21 	&      3 	&      184	&      188		  \\  
52	&    26	 &    8	&    3		&    24  	&    104 	&      28 	&      5 	&      236	&      240		  \\  \hline
58	&    29	 &    9	&    0		&    30  	&    160 	&      28 	&      5 	&      292	&      296		  \\  
64	&    32	 &    10	&    1		&    32  	&    181 	&      36 	&      7 	&      353	&      358		  \\  
70 	&    35	 &    11	&    2		&    34    	&    204    	&      45  	&      9  	&      420	&      426		  \\  
76 	&    38	 &    12	&    3		&    36   	&     228   	&      55  	&      12  	&      496	&      502  		 \\  \hline
82 	&    41	 &    13	&    0		&     42   	&     308   	&      55  	&      12  	&      576	&      582  		\\ 
88 	&    44	 &    14	&    1		&     44   	&     337   	&      66  	&      15  	&      661	&      668  		\\ 
94 	&    47	 &    15	&    2		&     46   	&     368   	&      78  	&      18  	&      752	&      760  		\\   
100 	&    50	 &    16	&    3		&     48   	&     400   	&      91  	&      22  	&      852	&      860  		\\  \hline
106 	&    53	 &    17	&    0		&     54   	&     504   	&      91  	&      22  	&      956	&      964  		\\  	
112 	&    56	 &    18	&    1		&     56   	&     541   	&      105 	&      26 	&      1065	&      1074  	\\ 	
118 	&    59	 &    19	&    2		&     58   	&     580   	&      120 	&      30 	&      1180	&      1190	  	\\  
124 	&    62	 &    20	&    3		&     60   	&     620   	&      136 	&      35 	&      1304	&      1314  	\\  \hline
130 	&    65	 &    21	&    0		&     66   	&     748   	&      136 	&      35 	&      1432	&      1442  	\\ 	
136 	&    68	 &    22	&    1		&     68   	&     793   	&      153 	&      40 	&      1565	&      1576  	\\ 	
142 	&    71	 &    23	&    2		&     70   	&     840   	&      171 	&      45 	&      1704 &      1716  	\\  	
148 	&    74	 &    24	&    3		&     72   	&     888   	&      190 	&      51 	&      1852	&      1864  	\\  \hline	
 154 	&    77	 &    25	&    0		&     78   	&   1040   &     190 	&      51 	&      2004	&      2016  	\\  	
 160 	&    80	 &    26	&    1		&     80   	&   1093   &     210 	&      57 	&      2161	&      2174  	\\  	
 166 	&    83	 &    27	&    2		&     82   	&   1148  	&     231 	&      63 	&      2324	&      2338  	\\   
 172 	&    86	 &    28	&    3		&     84   	&   1204   	&     253 	&      70 	&      2496	&      2510  	\\  \hline	
 178 	&    89	 &    29	&    0		&     90   	&   1380   	&      253 	&      70 	&      2672	&      2686  	\\  	
 184 	&    92	 &    30	&    1		&     92   	&   1441   	&      276 	&      77 	&      2853 &     2868  	\\  	
 190 	&    95	 &    31	&    2		&     94   	&   1504   	&      300 	&      84 	&      3040	&      3056  	\\  
 196 	&    98	 &    32	&    3		&     96   	&   1568   	&      325 	&      92 	&      3236	&      3252  	\\  \hline	
 202 	&    101	 &    33	&    0		&     102   	&   1768   	&      325 	&      92 	&      3436	&      3452  	\\  \hline	
\end{tabular}
\end{table}
%
%
\end{document}